\documentclass[10pt]{amsart}
\usepackage[T1]{fontenc}
\usepackage{geometry}
\usepackage[latin1] {inputenc}
\usepackage{amsmath}
\usepackage{amsfonts, amssymb, textcomp}
\usepackage[colorlinks=flase, linkcolor=red,urlcolor=green, citecolor=blue]{hyperref}
\usepackage{subeqnarray}

\usepackage{latexsym}
\usepackage{fancyhdr}
\usepackage{longtable}
\usepackage{amsmath, amssymb}
\usepackage{graphicx}
\usepackage{enumerate}

\numberwithin{equation}{section}

\theoremstyle{plain}
\newtheorem{proposition}{Proposition}[section]
\newtheorem{theorem}[proposition]{Theorem}
\newtheorem{lemma}[proposition]{Lemma}
\newtheorem{corollary}[proposition]{Corollary}
\newtheorem{definition}[proposition]{Definition}
\newtheorem{example}[proposition]{Example}

\newtheorem{remark}[proposition]{Remark}

\newcommand{\RR}{\mathbb{R}}

\newcommand{\NN}{\mathbb{N}}

\newcommand{\id}{\operatorname{id}}

\let\on=\operatorname

\newenvironment{demo}[1]{\par\smallskip\noindent{\bf #1.}}{\par\smallskip}

\makeatletter
\@namedef{subjclassname@2020}{%
  \textup{2020} Mathematics Subject Classification}
\makeatother

\newsavebox{\fmbox}
\newenvironment{fmpage}[1]
 {\begin{lrbox}{\fmbox}\begin{minipage}{#1}}
 {\end{minipage}\end{lrbox}\fbox{\usebox{\fmbox}}}

\title[Generalized upper and lower Legendre conjugates for weight functions]
{Generalized upper and lower Legendre conjugates for weight functions}

\author[G.~Schindl]{Gerhard Schindl}

\address{G.~Schindl: Fakult\"at f\"ur Mathematik, Universit\"at Wien, Oskar-Morgenstern-Platz~1, A-1090 Wien, Austria.}
\email{gerhard.schindl@univie.ac.at}

\begin{document}

\begin{abstract}
We introduce and study new transformations between two functions satisfying some basic growth properties and generalize the known lower and upper Legendre conjugate (or envelope). We also investigate how these transformations modify recently defined growth indices for weight functions. A special but important and useful situation, to which the knowledge is then applied, is when considering associated weight functions which are expressed in terms of an underlying weight sequence. In this case these transformations precisely correspond to the point-wise product resp. point-wise division of the given sequences. Therefore, the new approach studied in this work illustrates the genuineness and importance and suggests applications for weighted spaces in different directions.
\end{abstract}

\thanks{This research was funded in whole or in part by the Austrian Science Fund (FWF) 10.55776/PAT9445424}
\date{\today}

\maketitle

\section{Introduction}
Weighted spaces appear in different contexts and fields in Mathematics, like in Functional Analysis and the theory of partial differential operators but also in applications in various directions. We mention \emph{ultradifferentiable and ultraholomorphic classes,} see \cite{Komatsu73}, \cite{BraunMeiseTaylor90}, \cite{BonetMeiseMelikhov07}, \cite{compositionpaper}, and \cite{Thilliezdivision}, \cite{optimalflat23}; \emph{weighted spaces of sequences of complex numbers} or \emph{weighted spaces of formal power series,} see \cite{petzsche} and \cite{Borelmapalgebraity}; generalized \emph{Gelfand-Shilov classes,} see \cite{nuclearglobal1}, \cite{nuclearglobal2}, \cite{GelfandShilovincl24}; \emph{weighted spaces of entire functions,} see \cite{Bonet2022survey}, \cite{weightedentireinclusion1} and \cite{weightedentireinclusion2}; and \emph{Orlicz spaces,} see \cite{orliczpaper}. In each instance we also refer to the list of references in these works.\vspace{6pt}

Classically, one can find two approaches: either using a \emph{weight sequence} $\mathbf{M}=(M_p)_{p\in\NN}$ or a \emph{weight function} $\omega:[0,+\infty)\rightarrow[0,+\infty)$. More recently, even \emph{weight matrices} or \emph{weight function systems} have been considered; see e.g. \cite{compositionpaper}, \cite{GelfandShilovincl24}. In all these settings growth and regularity assumptions on the weights are required and unavoidable. In general, the weight sequence and weight function setting are mutually distinct which follows in the ultradifferentiable case by \cite{BonetMeiseMelikhov07}, see also \cite[Sect. 5]{compositionpaper}, and in the Gelfand-Shilov setting we refer to \cite[Sect. 5]{GelfandShilovincl}. Note that in the mentioned settings similar growth conditions appear; this has motivated the abstract study in \cite{index} where especially weight functions have been investigated under very general growth assumptions.\vspace{6pt}

Given a weight sequence $\mathbf{M}=(M_p)_{p\in\NN}\in\RR_{>0}^{\NN}$ satisfying basic growth properties one can consider the corresponding \emph{associated function} $\omega_{\mathbf{M}}$; see \cite[Chapitre I]{mandelbrojtbook}, \cite[Sect. 3]{Komatsu73}, the more recent work \cite{regularnew} and Section \ref{assofunctsect} below for more details. If $\mathbf{M}$ is \emph{log-convex} (see \eqref{logconv}), then there exists an intimate relation between $\mathbf{M}$ and $\omega_{\mathbf{M}}$. Note that $\omega_{\mathbf{M}}$ appears frequently in both the weight sequence and weight function setting and in the latter case it is also serving as a (counter-)example and is becoming relevant for comparison results between these settings; see \cite{BonetMeiseMelikhov07}.

On the other hand, for abstractly given weight functions, again satisfying basic growth restrictions, one can assign an \emph{associated weight sequence,} see \cite[Sect. 2.7]{weightedentireinclusion1} and \cite[Sect. 5]{orliczpaper}, or even a \emph{weight matrix,} see \cite[Sect. 5]{compositionpaper} and \cite{dissertation}.\vspace{6pt}

From an abstract point of view but also for concrete applications and when working with weighted structures one is interested in the following questions:
\begin{itemize}
\item[$(I)$] Let sequences $\mathbf{M}$, $\mathbf{N}$ be given and write $\mathbf{M}\cdot\mathbf{N}:=(M_pN_p)_{p\in\NN}$ and $\frac{\mathbf{M}}{\mathbf{N}}:=\left(\frac{M_p}{N_p}\right)_{p\in\NN}$ for the \emph{point-wise product resp. point-wise quotient sequence.} How is then $\omega_{\mathbf{M}\cdot\mathbf{N}}$, $\omega_{\frac{\mathbf{M}}{\mathbf{N}}}$ related to $\omega_{\mathbf{M}}$, $\omega_{\mathbf{N}}$?

\item[$(II)$] Study for more general resp. for a wider class of weight functions these analogous relations and transformations; more precisely one may ask:

    Given functions $\sigma,\tau:[0,+\infty)\rightarrow[0,+\infty)$, which function is then corresponding to the operation yielding $\omega_{\mathbf{M}\cdot\mathbf{N}}$ resp. $\omega_{\frac{\mathbf{M}}{\mathbf{N}}}$?
\end{itemize}
The content of this work is to answer and study $(I)$ and $(II)$ in detail: It turns out that taking the point-wise product is corresponding to
\begin{equation}\label{wedgeformula}
\sigma\check{\star}\tau(t):=\inf_{s>0}\{\sigma(s)+\tau(t/s)\},\;\;\;t\in[0,+\infty),
\end{equation}
whereas the point-wise quotient is corresponding to
\begin{equation}\label{widehatformula}
\sigma\widehat{\star}\tau(t):=\sup_{s\ge 0}\{\sigma(s)-\tau(s/t)\},\;\;\;t\in(0,+\infty).
\end{equation}
The choice $\tau=\id$ in \eqref{wedgeformula} gives, up to inversion of the variable, the known \emph{lower Legendre conjugate (or envelope) of $\sigma$} and \eqref{widehatformula} for $\tau=\id$ gives the \emph{upper Legendre conjugate (or envelope) of $\sigma$;} see e.g. \cite[Sect. 2.5]{index} and the references therein and confirm formulas \eqref{lowerenvelope} and \eqref{upperLegendre}. It is natural to study $(II)$ for weight functions in the sense of \cite{index}, see Definition \ref{weightfctdef}, and this notion encompasses in particular weight functions in the sense of \emph{Braun-Meise-Taylor} introduced in \cite{BraunMeiseTaylor90}. In \cite{index} for weight functions also new growth indices $\gamma(\cdot)$ and $\overline{\gamma}(\cdot)$ have been introduced, see Section \ref{growthindexsect}, and hence it is also natural to ask how these indices are modified under the actions $\check{\star}$ and $\widehat{\star}$.\vspace{6pt}

Indeed, the starting point of our investigations and interest has been the technical statements shown in \cite[Sect. 3]{sectorialextensions}. There we have revisited the lower/upper Legendre conjugate and established a connection between these conjugates and the associated functions $\omega_{\frac{\mathbf{M}}{\mathbf{G}^s}}$ resp. $\omega_{\mathbf{M}\cdot\mathbf{G}^1}$, $\mathbf{G}^s:=(p!^s)_{p\in\NN}$, $s>0$, see \cite[Lemmas 3.1, 3.4 \& 3.5]{sectorialextensions}. Note that $\omega_{\mathbf{G}^1}$ corresponds (up to equivalence) to the weight $\id: t\mapsto t$ and taking the power of a sequence is yielding a power substitution for the associated function (recall \cite[$(2.7)$]{sectorialextensions}). Moreover, in the ultraholomorphic setting (for proving extension results) the point-wise product resp. point-wise quotient of a sequence $\mathbf{M}$ and the \emph{Gevrey sequence} $\mathbf{G}^s$ appears naturally and is intimately related to the \emph{(formal) $s$-Laplace resp. $s$-Borel transform;} we refer to \cite[Thm. 3.5]{surjectivityregularsequence} and the citations and explanations stated just before this result. Then we have found \cite[Lemma 4]{Boman98} where an explicit relation between $\omega_{\mathbf{M}\cdot\mathbf{N}}$ and $\omega_{\mathbf{M}}$, $\omega_{\mathbf{N}}$ has been established in terms of \eqref{wedgeformula}. This operation is called \emph{``infimal convolution''} in \cite{Boman98} and in view of this $\widehat{\star}$ can be called \emph{``supremal anti-convolution''} but these terminologies will not be used in this work. To the best of our knowledge the analogous result for $\omega_{\frac{\mathbf{M}}{\mathbf{N}}}$, when using \eqref{widehatformula}, has not been shown so far and it seems that formally the proof given in \cite[Lemma 4]{Boman98} is not correct since the assumptions for $\mathbf{M}$, $\mathbf{N}$ there do not ensure necessarily the required property that the corresponding sequences of quotients tend to infinity. And when this standard behavior fails for $\mathbf{M}$ then it is known that $\omega_{\mathbf{M}}(t)=+\infty$ for all $t$ sufficiently large (see Lemma \ref{assofunctsectlemma}); in particular, $\omega_{\mathbf{M}}$ is not a weight function.

Moreover, in order to ensure that $\sigma\widehat{\star}\tau$ and $\omega_{\frac{\mathbf{M}}{\mathbf{N}}}$ are well defined some more technical investigations have to be taken into account and the weight functions resp. sequences have to satisfy a certain growth relation; see Sections \ref{arbuppertrafosect} and \ref{uppertrafosect} for details. Roughly speaking this is due to the fact that division destroys regularity; i.e. even if both $\mathbf{M}$ and $\mathbf{N}$ are very regular then the point-wise quotient $\frac{\mathbf{M}}{\mathbf{N}}$ will be not in general.\vspace{6pt}

For the sake of completeness in the weight sequence setting we are also investigating \emph{non-standard situations} and using the detailed knowledge from \cite{regularnew}; see Sections \ref{nonstandardlowersection}, \ref{nonstandardsection} and \ref{invnonstandardsection} in the Appendix \ref{nonstandardappendix}. More precisely, we are interested in the above mentioned situation when $\omega_{\mathbf{M}}(t)=+\infty$ holds for all sufficiently large $t$ and hence such associated functions are violating the standard assumptions of being a weight function. On the one hand, for applications in weighted settings this case is restricting resp. yielding trivial situations. On the other hand, from the point of view of studying growth properties for abstractly given sequences $\mathbf{M}\in\RR_{>0}^{\NN}$ it is natural to consider also these cases and to investigate the qualitative differences compared with the standard setting. The study of non-standard cases has also been motivated by the general formulation of \cite[Lemma 4]{Boman98}.\vspace{6pt}

Since we are dealing with weight functions in an abstract and general setting we expect that the obtained results in this work allow for applications in different weighted settings and in different contexts and we already mention two very recent uses: In \cite{genLegendreconjBMT} we apply the knowledge of this work to some more specific weight functions, more precisely we consider \emph{weights in the sense of Braun-Meise-Taylor} (see \cite{BraunMeiseTaylor90}) and study the effects of the conjugates in terms of the corresponding \emph{associated weight matrices $\mathcal{M}_{\omega}$.} This notion and method is introduced in \cite{compositionpaper} and \cite{dissertation}. And in \cite{modgrowthstrangeII} we are continuing our study from \cite{modgrowthstrange}. Via applying $\check{\star}$ and the main result Theorem \ref{propBomanLemma4} we are able to characterize a new crucial condition for weight sequences $\mathbf{M}$ in terms of the corresponding associated weight functions $\omega_{\mathbf{M}}$; this new condition characterizes the \emph{quotient-root comparison properties} of the weight matrices $\mathcal{M}_{\omega_{\mathbf{M}}}$ associated with $\omega_{\mathbf{M}}$. We refer to the main results \cite[Thm. 5.3, Cor. 5.4 \& 5.5]{modgrowthstrangeII}.\vspace{6pt}

The paper is structured as follows: In Section \ref{weightsequassofctsect} we revisit crucial definitions, relations and growth and regularity properties for weight sequences and (associated) weight functions. Then, in Sections \ref{arblowertrafosect} and \ref{arbuppertrafosect} we study the new operations $\check{\star}$ and $\widehat{\star}$ between weight functions and their effect on the growth indices $\gamma(\cdot)$, $\overline{\gamma}(\cdot)$; see Theorems \ref{lowertransformindexthm} and \ref{uppertransformindexthm}. In Section \ref{trafosect} we focus on associated weight functions. Here, since the functions are based on sequences this underlying additional information enables to verify Theorems \ref{propBomanLemma4} and \ref{propBomanLemma4inv} which solves $(I)$. In Section \ref{supplsection} we exploit more known results for weight sequences and prove, in particular, how $\omega_{\mathbf{M}}\check{\star}\omega_{\mathbf{N}}$ and $\omega_{\mathbf{M}}\widehat{\star}\omega_{\mathbf{N}}$ are related to $\omega_{\mathbf{M}}$, $\omega_{\mathbf{N}}$; see Propositions \ref{multlemma65} and \ref{divlemma65}. In Section \ref{inversection} we are showing that within the weight sequence setting the operations $\check{\star}$ and $\widehat{\star}$ are inverse to each other; see Theorem \ref{propBomanLemma4comb}. Finally, in the Appendix in Sections \ref{nonstandardlowersection}, \ref{nonstandardsection} and \ref{invnonstandardsection} we are commenting on non-standard situations in the weight sequence setting and see how the results from Sections \ref{trafosect} and \ref{inversection} are modified.\vspace{6pt}

\textbf{Acknowledgements.} The author of this article thanks the two anonymous referees for the careful reading and the valuable suggestions which have improved the presentation of the results.

\section{Weight sequences, weight functions and associated functions}\label{weightsequassofctsect}

\subsection{Basic notation}
We write $\NN:=\{0,1,2,\dots\}$, $\NN_{>0}:=\{1,2,\dots\}$, $\RR_{>0}:=(0,+\infty)$ and $\RR_{\ge 0}:=[0,+\infty)$. Occasionally, we use the conventions $\frac{0}{0}=0$ and $\frac{c}{0}=+\infty$ for $c>0$ and write that an inequality $a\le b$ holds formally if either $a=+\infty$ or $b=+\infty$. In the first case necessarily $b=+\infty$ follows then.

\subsection{On weight sequences}\label{preliminarysection}
Let $\mathbf{M}=(M_p)_{p\in\NN}\in\RR_{>0}^{\NN}$ be a sequence of positive real numbers. Moreover, let us set $\mu_p:=\frac{M_p}{M_{p-1}}$ for $p\ge 1$ and $\mu_0:=1$. Hence $\frac{M_p}{M_0}=\mu_p\cdots\mu_1$ for all $p\in\NN$; the case $p=0$ gives the empty product. $\mathbf{M}$ is called \emph{normalized} if $1=M_0\le M_1$.

Given $\mathbf{M},\mathbf{N}$, write $\mathbf{M}\cdot\mathbf{N}$ for the sequence obtained by $(M_p\cdot N_p)_{p\in\NN}$ (point-wise product) and $\frac{\mathbf{M}}{\mathbf{N}}$ for $\left(\frac{M_p}{N_p}\right)_{p\in\NN}$ (point-wise quotient). The corresponding sequences of quotients are denoted and given by $\mu\cdot\nu$ and $\frac{\mu}{\nu}$.

$\mathbf{M}$ is \emph{log-convex} if
\begin{equation}\label{logconv}
\forall\;p\in\NN_{>0}:\;\;\;M_p^2\le M_{p-1}M_{p+1},
\end{equation}
in the literature also denoted by $(M.1)$; see e.g \cite{Komatsu73}.

For concrete applications and examples it is convenient to introduce the following class of sequences:
$$\hypertarget{LCset}{\mathcal{LC}}:=\{\mathbf{M}\in\RR_{>0}^{\NN}:\;\mathbf{M}\;\text{is normalized, log-convex},\;\lim_{p\rightarrow+\infty}(M_p)^{1/p}=+\infty\}.$$
\emph{Note:} For each $\mathbf{M}\in\hyperlink{LCset}{\mathcal{LC}}$ we have $\mu_0=1\le\mu_1\le\mu_2\le\dots$ and $\lim_{p\rightarrow+\infty}\mu_p=+\infty$. In particular, $\mathbf{M}$ is non-decreasing and in fact there exists a one-to-one correspondence between sequences $\mathbf{M}$ belonging to the set \hyperlink{LCset}{$\mathcal{LC}$} and sequences $\mu_0=1\le\mu_1\le\mu_2\le\dots$ with $\lim_{p\rightarrow+\infty}\mu_p=+\infty$ via setting $M_p:=\prod_{i=0}^p\mu_i$.

Next let us introduce relations between given sequences $\mathbf{M}$, $\mathbf{N}$:

\begin{itemize}
\item[$(*)$] Write $\mathbf{M}\le\mathbf{N}$ if $M_p\le N_p$ for all $p\in\NN$,

\item[$(*)$] $\mathbf{M}\hypertarget{preceq}{\preceq}\mathbf{N}$ if
$$\sup_{p\in\NN_{>0}}\left(\frac{M_p}{N_p}\right)^{1/p}<+\infty;$$
i.e. if $M_p\le C^{p+1}N_p$ for some $C\ge 1$ and all $p\in\NN$.

\item[$(*)$] $\mathbf{M},\mathbf{N}\in\RR_{>0}^{\NN}$ are called \emph{equivalent} if $\mathbf{M}\hyperlink{preceq}{\preceq}\mathbf{N}$ and $\mathbf{N}\hyperlink{preceq}{\preceq}\mathbf{M}$; i.e. if
\begin{equation}\label{equivalencrelation}
0<\inf_{p\in\NN_{>0}}\left(\frac{M_p}{N_p}\right)^{1/p}\le\sup_{p\in\NN_{>0}}\left(\frac{M_p}{N_p}\right)^{1/p}<+\infty,
\end{equation}
equivalently
$$\exists\;C\ge 1\;\forall\;p\in\NN:\;\;\;\frac{1}{C^{p+1}}N_p\le M_p\le C^{p+1}N_p.$$
\item[$(*)$] Write $\mathbf{M}\hypertarget{mtriangle}{\vartriangleleft}\mathbf{N}$ if $\lim_{p\rightarrow+\infty}\left(\frac{M_p}{N_p}\right)^{1/p}=0$; i.e. if
\begin{equation}\label{triangleestim}
\forall\;h>0\;\exists\;C_h\ge 1\;\forall\;p\in\NN:\;\;\;M_p\le C_hh^pN_p.
\end{equation}
Obviously, $\mathbf{M}\hyperlink{mtriangle}{\vartriangleleft}\mathbf{N}$ implies $\mathbf{M}\hyperlink{preceq}{\preceq}\mathbf{N}$, but \hyperlink{mtriangle}{$\vartriangleleft$} neither is reflexive nor symmetric.
\end{itemize}

\begin{remark}\label{Cinfremark}
\emph{\hyperlink{preceq}{$\preceq$} and \hyperlink{mtriangle}{$\vartriangleleft$} are crucial for the inclusion relations of the corresponding weighted spaces and, in particular, equivalent sequences yield the same weighted class. Therefore, in concrete applications it is no restriction to assume always $M_0=1$ since otherwise $\mathbf{M}$ can be replaced by $\widetilde{\mathbf{M}}$ with $\widetilde{M}_p:=M_p/M_0$. Similarly, for each $\mathbf{M}\in\RR_{>0}^{\NN}$ one can find an equivalent sequence $\mathbf{N}$ which is even normalized.}

\emph{Finally, for any given log-convex sequence $\mathbf{M}$ satisfying $\lim_{p\rightarrow+\infty}(M_p)^{1/p}=+\infty$ there exists an equivalent sequence $\mathbf{N}\in\hyperlink{LCset}{\mathcal{LC}}$; see \cite[Rem. 2.5]{weightedentireinclusion1} and \cite[Prop. 3.5]{regularnew}.}
\end{remark}

A prominent example are the \emph{Gevrey sequences} $\mathbf{G}^s:=(p!^s)_{p\in\NN}$, $s>0$. Note that $\mathbf{G}^s$ is equivalent to $\overline{\mathbf{G}}^s$ with $\overline{\mathbf{G}}^s:=(p^{ps})_{p\in\NN}$ by \emph{Stirling's formula.}\vspace{6pt}

We revisit some notions from \cite[Sect. 2.5]{regularnew}. Let $\mathbf{M}\in\RR_{>0}^{\NN}$, then put
\begin{equation}\label{liminfcond}
\mathbf{M}_{\iota}:=\liminf_{p\rightarrow+\infty}\left(\frac{M_p}{M_0}\right)^{1/p}=\liminf_{p\rightarrow+\infty}(M_p)^{1/p},
\end{equation}
and
\begin{equation}\label{infcond}
\mathbf{M}_{\inf}:=\inf_{p\in\NN_{>0}}\left(\frac{M_p}{M_0}\right)^{1/p}.
\end{equation}

Consequently:

\begin{itemize}
\item[$(a)$] $\mathbf{M}_{\inf}\le\mathbf{M}_{\iota}$, $(\mathbf{M}\cdot\mathbf{N})_{\inf}\ge\mathbf{M}_{\inf}\cdot\mathbf{N}_{\inf}$ and $(\mathbf{M}\cdot\mathbf{N})_{\iota}\ge\mathbf{M}_{\iota}\cdot\mathbf{N}_{\iota}$.

\item[$(b)$] If $\mathbf{M}$ is also log-convex then, as mentioned in \cite[Sect.2.2]{regularnew} (see also \cite[Rem. 4.1]{conjugateweightfunction}), we get that $p\mapsto\left(\frac{M_p}{M_0}\right)^{1/p}$ is non-decreasing and
$$\mathbf{M}_{\inf}=\frac{M_1}{M_0}=\mu_1,\hspace{15pt}\lim_{p\rightarrow+\infty}\mu_p=\lim_{p\rightarrow+\infty}\left(\frac{M_p}{M_0}\right)^{1/p}=\mathbf{M}_{\iota}\in(0,+\infty].$$

\item[$(c)$] Note that $\mathbf{M}_{\iota}=0$ for a log-convex sequence corresponds to the ``extreme (limiting) case'' $M_0\ge 0$ and $M_p=0$ for all $p\in\NN_{>0}$ studied in \cite[Sect. 4.1]{regularnew}. However, formally this violates the (basic) assumption $\mathbf{M}\in\RR_{>0}^{\NN}$.
\end{itemize}

\emph{Based on the above comments, especially see Remark \ref{Cinfremark}, in the main results in this work we focus on the case $\mathbf{M}_{\iota}=+\infty$. Comments on non-standard situations are summarized in the separate Appendix \ref{nonstandardappendix}.}\vspace{6pt}

Let now $\mathbf{M},\mathbf{N}\in\RR_{>0}^{\NN}$ be given and remark that
\begin{equation}\label{MNrelations}
0\le\inf_{p\in\NN_{>0}}\left(\frac{M_p/M_0}{N_p/N_0}\right)^{1/p}\le\liminf_{p\rightarrow+\infty}\left(\frac{M_p/M_0}{N_p/N_0}\right)^{1/p}\le\limsup_{p\rightarrow+\infty}\left(\frac{M_p/M_0}{N_p/N_0}\right)^{1/p}\le\sup_{p\in\NN_{>0}}\left(\frac{M_p/M_0}{N_p/N_0}\right)^{1/p}\le+\infty,
\end{equation}
and
$$0<\inf_{p\in\NN_{>0}}\left(\frac{M_p/M_0}{N_p/N_0}\right)^{1/p}\Longleftrightarrow 0<\liminf_{p\rightarrow+\infty}\left(\frac{M_p/M_0}{N_p/N_0}\right)^{1/p},$$ $$\underline{C}_{\mathbf{M}\preceq\mathbf{N}}:=\limsup_{p\rightarrow+\infty}\left(\frac{M_p/M_0}{N_p/N_0}\right)^{1/p}<+\infty\Longleftrightarrow C_{\mathbf{M}\preceq\mathbf{N}}:=\sup_{p\in\NN_{>0}}\left(\frac{M_p/M_0}{N_p/N_0}\right)^{1/p}<+\infty.$$
Alternatively we get
$$C_{\mathbf{M}\preceq\mathbf{N}}=\inf\{C>0:\;\frac{M_p}{M_0}\le C^p\frac{N_p}{N_0},\;\forall\;p\in\NN\}$$
and summarize:

\begin{itemize}
\item[$(*)$] $C_{\mathbf{M}\preceq\mathbf{N}}\in(0,+\infty]$ whereas $\underline{C}_{\mathbf{M}\preceq\mathbf{N}}\in[0,+\infty]$ holds.

\item[$(*)$] $\underline{C}_{\mathbf{M}\preceq\mathbf{N}}=0$ if and only if $\mathbf{M}\hyperlink{mtriangle}{\vartriangleleft}\mathbf{N}$.

\item[$(*)$] $\mathbf{M}\hyperlink{preceq}{\preceq}\mathbf{N}$ if and only if $\underline{C}_{\mathbf{M}\preceq\mathbf{N}}<+\infty$ resp. if and only if $C_{\mathbf{M}\preceq\mathbf{N}}<+\infty$. Indeed, $\mathbf{M}\hyperlink{preceq}{\preceq}\mathbf{N}$ amounts to having
    \begin{equation}\label{preceqalternative}
    \forall\;p\in\NN:\;\;\;\frac{M_p}{M_0}\le C_{\mathbf{M}\preceq\mathbf{N}}^p\frac{N_p}{N_0},
    \end{equation}
and $\mathbf{M}\hyperlink{preceq}{\preceq}\mathbf{N}$ is equivalent to
\begin{equation}\label{preceqalternative1}
\forall\;h>\underline{C}_{\mathbf{M}\preceq\mathbf{N}}\;\exists\;C_h>0\;\forall\;p\in\NN:\;\;\;M_p\le C_hh^pN_p,
\end{equation}
which should be compared with \eqref{triangleestim}.

\item[$(*)$] When $\underline{C}_{\mathbf{M}\preceq\mathbf{N}}\in(0,+\infty)$, i.e. $\mathbf{M}\hyperlink{preceq}{\preceq}\mathbf{N}$ holds but $\mathbf{M}\hyperlink{vartriangle}{\vartriangleleft}\mathbf{N}$ fails, then there exists some $D\ge 1$ such that $C_{\mathbf{M}\preceq\mathbf{N}}\le D\underline{C}_{\mathbf{M}\preceq\mathbf{N}}$. (This estimate is ``meaningless/trivial'' if $\underline{C}_{\mathbf{M}\preceq\mathbf{N}}=+\infty=C_{\mathbf{M}\preceq\mathbf{N}}$.)

\item[$(*)$] When $\underline{C}_{\mathbf{M}\preceq\mathbf{N}}\in(0,+\infty)$, then by definition
\begin{equation}\label{MNiotarelation}
\mathbf{M}_{\iota}\le\underline{C}_{\mathbf{M}\preceq\mathbf{N}}\mathbf{N}_{\iota}.
\end{equation}
If $\mathbf{M}_{\iota}=+\infty$, then also $\mathbf{N}_{\iota}=+\infty$ and hence (formally) equality. However, in general equality in \eqref{MNiotarelation} fails resp. is unclear: There might even exist some subsequence of integers $(p_k)_{k\in\NN_{>0}}$ such that $\sup_{k\in\NN_{>0}}\left(\frac{M_{p_k}/M_0}{N_{p_k}/N_0}\right)^{1/p_k}<\underline{C}_{\mathbf{M}\preceq\mathbf{N}}$ and such that $\mathbf{M}_{\iota}=\lim_{k\rightarrow+\infty}\left(\frac{M_{p_k}}{M_0}\right)^{1/p_k}$.

\item[$(*)$] Concerning the ``extreme case'' $\underline{C}_{\mathbf{M}\preceq\mathbf{N}}=0$, i.e. if $\mathbf{M}\hyperlink{vartriangle}{\vartriangleleft}\mathbf{N}$, we remark: When $\mathbf{M}_{\iota}\in(0,+\infty]$, then necessarily $\mathbf{N}_{\iota}=\lim_{p\rightarrow+\infty}(N_p)^{1/p}=+\infty$ is valid.
\end{itemize}

Related to these observations and in view of the quotient sequence $\frac{\mathbf{M}}{\mathbf{N}}$  we comment on the behavior of $p\mapsto\left(\frac{M_p}{N_p}\right)^{1/p}$ appearing in $\mathbf{M}\hyperlink{mtriangle}{\vartriangleleft}\mathbf{N}$ and $\mathbf{M}\hyperlink{preceq}{\preceq}\mathbf{N}$.

\begin{itemize}
\item[$(i)$] Even if both $\mathbf{M}$ and $\mathbf{N}$ are log-convex, then in general $\frac{\mathbf{M}}{\mathbf{N}}$ \emph{will not be log-convex anymore:} Consider, e.g., $\mathbf{M}=\mathbf{G}^{s_1}$, $\mathbf{N}=\mathbf{G}^{s_2}$, with $0<s_1<s_2$, and so $\frac{M_p}{N_p}=p!^{s_1-s_2}$. The corresponding quotient sequence is given by $p^{s_1-s_2}$ which is even decreasing and tending to $0$ as $p\rightarrow+\infty$. Thus $\frac{\mathbf{M}}{\mathbf{N}}$ is \emph{log-concave} which means that the estimate in \eqref{logconv} is reversed. On the other hand, clearly $\frac{\mathbf{G}^{s_2}}{\mathbf{G}^{s_1}}=\mathbf{G}^{s_2-s_1}\in\hyperlink{LCset}{\mathcal{LC}}$.

\item[$(ii)$] Let $\mathbf{M}$ and $\mathbf{N}$ be log-convex.

Then $\lim_{p\rightarrow+\infty}(M_p/M_0)^{1/p}=\mathbf{M}_{\iota}\in(0,+\infty]$, $\lim_{p\rightarrow+\infty}(N_p/N_0)^{1/p}=\mathbf{N}_{\iota}\in(0,+\infty]$ and concerning $\frac{\mathbf{M}}{\mathbf{N}}_{\iota}$ we have to distinguish:
    \begin{itemize}
    \item[$(*)$] If $\mathbf{M}_{\iota},\mathbf{N}_{\iota}<+\infty$, then $\frac{\mathbf{M}}{\mathbf{N}}_{\iota}=\frac{\mathbf{M}_{\iota}}{\mathbf{N}_{\iota}}\in(0,+\infty)$. Indeed, in this case $\mathbf{M}$ and $\mathbf{N}$ are equivalent.

\item[$(*)$] If $\mathbf{M}_{\iota}=+\infty$ and $\mathbf{N}_{\iota}<+\infty$, then $\frac{\mathbf{M}}{\mathbf{N}}_{\iota}=+\infty$. In this case $\mathbf{N}\hyperlink{mtriangle}{\vartriangleleft}\mathbf{M}$ holds.

\item[$(*)$] If $\mathbf{M}_{\iota}<+\infty$ and $\mathbf{N}_{\iota}=+\infty$, then $\frac{\mathbf{M}}{\mathbf{N}}_{\iota}=0$ and  $\mathbf{M}\hyperlink{mtriangle}{\vartriangleleft}\mathbf{N}$ is valid.

\item[$(*)$] The case $\mathbf{M}_{\iota}=+\infty=\mathbf{N}_{\iota}$ is not clear; here $\frac{\mathbf{M}}{\mathbf{N}}_{\iota}\in[0,+\infty]$ is possible and for the extreme cases see $(i)$ before. But even the existence of $\lim_{p\rightarrow+\infty}\left(\frac{M_p/M_0}{N_p/N_0}\right)^{1/p}$ is not ensured: Indeed, it is not difficult to construct log-convex sequences $\mathbf{M}$, $\mathbf{N}$ (with $M_0=1=N_0$ and $\mathbf{M}_{\iota}=+\infty=\mathbf{N}_{\iota}$) satisfying the ``strongest possible oscillation''
\begin{equation}\label{notcomparable}
0=\inf_{p\ge 1}\left(\frac{M_p}{N_p}\right)^{1/p}<\sup_{p\ge 1}\left(\frac{M_p}{N_p}\right)^{1/p}=+\infty.
\end{equation}
So neither $\mathbf{M}\hyperlink{preceq}{\preceq}\mathbf{N}$ nor $\mathbf{N}\hyperlink{preceq}{\preceq}\mathbf{M}$ is valid, i.e. $\mathbf{M}$, $\mathbf{N}$ are \emph{not comparable,} $\lim_{p\rightarrow+\infty}\left(\frac{M_p}{N_p}\right)^{1/p}$ does not exist and $\frac{\mathbf{M}}{\mathbf{N}}_{\iota}=0$.

For an explicit but involved example we refer to \cite[Sect. 3]{GelfandShilovincl}: From given $\mathbf{N}\in\hyperlink{LCset}{\mathcal{LC}}$ one can construct $\mathbf{M}\in\hyperlink{LCset}{\mathcal{LC}}$ such that \eqref{notcomparable} holds. (In \cite{GelfandShilovincl} this question is related to the problem of non-triviality of Gelfand-Shilov-type spaces and so both $\mathbf{M}$ and $\mathbf{N}$ are required to have some more specific regularity properties.)
\end{itemize}

\item[$(iii)$] Note that the growth of $p\mapsto\left(\frac{M_p}{N_p}\right)^{1/p}$ is crucial for the inclusion relation of the corresponding weighted classes expressed in terms of $\mathbf{M}$ and $\mathbf{N}$; concerning the ultradifferentiable setting we refer e.g. to \cite[Prop. 2.12 \& 4.6]{compositionpaper}.
\end{itemize}

\subsection{Weight functions}\label{weightfctsection}
We start by recalling the following definition; see \cite[Sect. 2.3 \& 2.4]{index}.

\begin{definition}\label{weightfctdef}
$\omega:[0,+\infty)\rightarrow[0,+\infty)$ is called a \emph{weight function} if
\begin{itemize}
\item[$(*)$] $\omega$ is non-decreasing and

\item[$(*)$] $\lim_{t\rightarrow+\infty}\omega(t)=+\infty$.
\end{itemize}
\end{definition}
These growth requirements are sufficient in order to introduce and work with the crucial indices $\gamma(\omega)$ and $\overline{\gamma}(\omega)$; see Section \ref{growthindexsect} for details.

Let $\sigma,\tau: [0,+\infty)\rightarrow[0,+\infty)$ be weight functions.
\begin{itemize}
\item[$(*)$] Write $\sigma\hypertarget{ompreceq}{\preceq}\tau$ if
\begin{equation}\label{bigOrelation}
\tau(t)=O(\sigma(t))\;\text{as}\;t\rightarrow+\infty.	
\end{equation}
\item[$(*)$] $\sigma$ and $\tau$ are called \emph{equivalent,} denoted by $\sigma\hypertarget{sim}{\sim}\tau$, if $\sigma\hyperlink{ompreceq}{\preceq}\tau$ and $\tau\hyperlink{ompreceq}{\preceq}\sigma$.

\item[$(*)$] Write $\sigma\hypertarget{omvartriangle}{\vartriangleleft}\tau$ if
\begin{equation}\label{smallOrelation}
	\tau(t)=o(\sigma(t))\;\text{as}\;t\rightarrow+\infty,	
\end{equation}
and note that \hyperlink{omvartriangle}{$\vartriangleleft$} neither is reflexive nor symmetric.
\end{itemize}

We are going to use the following notation for any (weight) function $\omega:[0,+\infty)\rightarrow[0,+\infty)$ and arbitrary $\alpha>0$: set $\omega^{\iota}(t):=\omega(\frac{1}{t})$, $t>0$, and $\omega^{1/\alpha}(t):=\omega(t^{1/\alpha})$, $t\ge 0$. Note that $\omega^{1/\alpha}$ is obtained by a power-substitution and it is again a weight function. Moreover, $(\omega^{\iota})^{1/\alpha}(t)=(\omega^{1/\alpha})^{\iota}(t)$ for all $t>0$. Finally, let us write $\id^{1/\alpha}$ for $t\mapsto t^{1/\alpha}$, $\alpha>0$ arbitrary.

\subsection{The growth indices $\gamma(\omega)$ and $\overline{\gamma}(\omega)$}\label{growthindexsect}
We briefly recall the definitions of the growth indices $\gamma(\omega)$ and $\overline{\gamma}(\omega)$; see \cite[Sect. 2.3 $(6)$, Sect. 2.4 $(7)$]{index} and the references in these sections. Let $\omega:[0,+\infty)\rightarrow[0,+\infty)$ be a weight function and $\gamma>0$. We say that $\omega$ has property $(P_{\omega,\gamma})$ if
\begin{equation}\label{newindex1}
\exists\;K>1:\;\;\;\limsup_{t\rightarrow+\infty}\frac{\omega(K^{\gamma}t)}{\omega(t)}<K.
\end{equation}
If $(P_{\omega,\gamma})$ holds for some $K>1$, then $(P_{\omega,\gamma'})$ is satisfied for all $\gamma'\le\gamma$ with the same $K$ since $\omega$ is non-decreasing. Moreover we can restrict to $\gamma>0$, because for $\gamma\le 0$ condition $(P_{\omega,\gamma})$ is satisfied for any weight function $\omega$, again because $\omega$ is non-decreasing and $K>1$. Then put
\begin{equation}\label{newindex2}
\gamma(\omega):=\sup\{\gamma>0:\;\;(P_{\omega,\gamma})\;\;\text{is satisfied}\},
\end{equation}
and if none of the conditions $(P_{\omega,\gamma})$ (with $\gamma>0$) hold, then set $\gamma(\omega):=0$.

Analogously, for $\gamma>0$ we say that $\omega$ has property $(\overline{P}_{\omega,\gamma})$ if
\begin{equation}\label{newindex3}
\exists\;A>1:\;\;\;\liminf_{t\rightarrow+\infty}\frac{\omega(A^{\gamma}t)}{\omega(t)}>A.
\end{equation}
If $(\overline{P}_{\omega,\gamma})$ holds for some $A>1$, then $(\overline{P}_{\omega,\gamma'})$ is satisfied for all $\gamma'\ge\gamma$ with the same $A$ since $\omega$ is non-decreasing. Moreover, we can restrict to $\gamma>0$ because for $\gamma\le 0$ condition $(\overline{P}_{\omega,\gamma})$ is never satisfied for any weight function: $\omega$ is assumed to be non-decreasing and $A>1$. Then set
\begin{equation}\label{newindex4}
\overline{\gamma}(\omega):=\inf\{\gamma>0: \;\;(\overline{P}_{\omega,\gamma})\;\;\text{is satisfied}\}.
\end{equation}
We obtain for any weight function $(0\le)\gamma(\omega)\le\overline{\gamma}(\omega)$; see again \cite[Sect. 2.3 \& 2.4]{index}.

\begin{remark}\label{om1om6indexrem}
\emph{Let $\omega$ be a weight function and consider the following two growth conditions appearing for weight functions in the sense of Braun-Meise-Taylor, see e.g. \cite{BraunMeiseTaylor90} and \cite{BonetMeiseMelikhov07}, and the following abbreviations have already been used in \cite{dissertation}:}
\begin{itemize}
\item[\hypertarget{om1}{$(\omega_1)$}] $\omega(2t)=O(\omega(t))$ as $t\rightarrow+\infty$; i.e. $\exists\;L\ge 1\;\forall\;t\ge 0:\;\;\;\omega(2t)\le L\omega(t)+L$.
	
\item[\hypertarget{om6}{$(\omega_6)$}] $\exists\;H\ge 1\;\forall\;t\ge 0:\;\;\;2\omega(t)\le\omega(Ht)+H$.
\end{itemize}

\emph{$\omega$ has \hyperlink{om1}{$(\omega_1)$} if and only if $\gamma(\omega)>0$, see \cite[Thm. 2.11, Rem. 2.12, Cor. 2.14]{index}, and \hyperlink{om6}{$(\omega_6)$} holds if and only if $\overline{\gamma}(\omega)<+\infty$; see \cite[Thm. 2.16, $(7)$, Cor. 2.17]{index}.}
\end{remark}

\begin{remark}\label{Gevreyweightrem}
\emph{For any $\alpha>0$ we consider the weight function $\id^{1/\alpha}$ and obtain}
$$\gamma(\id^{1/\alpha})=\overline{\gamma}(\id^{1/\alpha})=\alpha.$$
\end{remark}

Next we verify the following useful consequence concerning the possible choices for $K$ and $A$ in \eqref{newindex1} and \eqref{newindex3}, respectively.

\begin{lemma}\label{bigKremark}
Let $\omega$ be a weight function and $\gamma>0$.
\begin{itemize}
\item[$(i)$] If $(P_{\omega,\gamma})$ is valid then there exists $K_0>1$ such that \eqref{newindex1} holds for all $K\ge K_0$.

\item[$(ii)$] If $(\overline{P}_{\omega,\gamma})$ is valid then there exists $A_0>1$ such that \eqref{newindex3} holds for all $A\ge A_0$.
\end{itemize}
\end{lemma}

\demo{Proof}
$(i)$ By \cite[Thm. 2.11 $(iii)\Leftrightarrow(iv)$]{index} and taking in $(iii)$ there $\alpha:=\gamma^{-1}$ we have that $\omega$ satisfies $(P_{\omega,\gamma})$ if and only if $\lim_{\epsilon\rightarrow 0}\limsup_{t\rightarrow+\infty}\frac{\epsilon^{1/\gamma}\omega(t)}{\omega(\epsilon t)}=0$. Set $\eta:=\frac{1}{\epsilon^{1/\gamma}}$ and $s:=\epsilon t=\frac{t}{\eta^{\gamma}}$ and hence the previous condition is equivalent to $\lim_{\eta\rightarrow+\infty}\limsup_{s\rightarrow+\infty}\frac{\omega(\eta^{\gamma}s)}{\eta\omega(s)}=0$. Thus we can find some $\eta_0>1$ such that for all $\eta\ge\eta_0$ one has $\limsup_{s\rightarrow+\infty}\frac{\omega(\eta^{\gamma}s)}{\eta\omega(s)}<1$ and this verifies, in particular, $(P_{\omega,\gamma})$ for each choice $K:=\eta\ge\eta_0$.\vspace{6pt}

$(ii)$ Similarly, by \cite[Thm. 2.16 $(ii)\Leftrightarrow(iii)$]{index} and taking in $(ii)$ there $\beta:=\gamma^{-1}$ we have that $\omega$ satisfies $(\overline{P}_{\omega,\gamma})$ if and only if $\lim_{k\rightarrow+\infty}\liminf_{t\rightarrow 0}\frac{\omega(kt)}{k^{1/\gamma}\omega(t)}=+\infty$. Take $\ell:=k^{1/\gamma}$ and hence the previous condition is equivalent to $\lim_{\ell\rightarrow+\infty}\liminf_{t\rightarrow+\infty}\frac{\omega(\ell^{\gamma}t)}{\ell\omega(t)}=+\infty$. Thus one can find some $\ell_0>1$ such that for all $\ell\ge\ell_0$ we get $\liminf_{t\rightarrow+\infty}\frac{\omega(\ell^{\gamma}t)}{\ell\omega(t)}>1$ and this verifies $(\overline{P}_{\omega,\gamma})$ for each choice $A:=\ell\ge\ell_0$.
\qed\enddemo

\emph{Note:} The conclusions of Lemma \ref{bigKremark} should be compared with the following direct observations. If \eqref{newindex1} holds for some $K>1$, then also for any $K'=K^n$, $n\in\NN_{>0}$, since by iteration
$$\limsup_{t\rightarrow+\infty}\frac{\omega(K^{n\gamma}t)}{\omega(t)}=\limsup_{t\rightarrow+\infty}\prod_{i=1}^n\frac{\omega(K^{i\gamma}t)}{\omega(K^{(i-1)\gamma}t)}\le\prod_{i=1}^n\limsup_{t\rightarrow+\infty}\frac{\omega(K^{i\gamma}t)}{\omega(K^{(i-1)\gamma}t)}<K^n;$$
see also the proof of \cite[Cor. 2.14]{index}. Similarly, if \eqref{newindex3} holds for some $A>1$, then also for all $A':=A^n$, $n\in\NN_{>0}$, because:
$$\liminf_{t\rightarrow+\infty}\frac{\omega(A^{n\gamma}t)}{\omega(t)}=\liminf_{t\rightarrow+\infty}\prod_{i=1}^n\frac{\omega(A^{i\gamma}t)}{\omega(A^{(i-1)\gamma}t)}\ge\prod_{i=1}^n\liminf_{t\rightarrow+\infty}\frac{\omega(A^{i\gamma}t)}{\omega(A^{(i-1)\gamma}t)}>A^n.$$

\subsection{Associated function}\label{assofunctsect}
For the following definition and the connection to the geometric construction of the \emph{log-convex minorant} we refer to \cite[Chapitre I]{mandelbrojtbook} and to the more recent work \cite{regularnew}. Let $\mathbf{M}=(M_p)_{p\in\NN}\in\RR_{>0}^{\NN}$ be given, then the \emph{associated function} $\omega_{\mathbf{M}}: [0,+\infty)\rightarrow[0,+\infty)\cup\{+\infty\}$ is defined as follows:
\begin{equation}\label{assofunc}
\omega_{\mathbf{M}}(t):=\sup_{p\in\NN}\log\left(\frac{M_0t^p}{M_p}\right),\qquad t\ge 0,
\end{equation}
with the conventions that $0^0:=1$ and $\log(0)=-\infty$. This ensures $\omega_{\mathbf{M}}(0)=0$ and $\omega_{\mathbf{M}}(t)\ge 0$ for any $t\ge 0$ since $\frac{t^0M_0}{M_0}=1$ for all $t\ge 0$. \eqref{assofunc} corresponds to \cite[$(3.1)$]{Komatsu73} and we immediately have that $\omega_{\mathbf{M}}$ is non-decreasing and satisfying $\lim_{t\rightarrow+\infty}\omega_{\mathbf{M}}(t)=+\infty$. Now recall \cite[Lemma 2.2]{regularnew}:

\begin{lemma}\label{assofunctsectlemma}
Let $\mathbf{M}\in\RR_{>0}^{\NN}$ be given.
\begin{itemize}
\item[$(i)$] $\mathbf{M}_{\iota}>0$ implies
$$\forall\;0\le t\le\mathbf{M}_{\inf}:\;\;\;\omega_{\mathbf{M}}(t)=0.$$
\item[$(ii)$] $\mathbf{M}_{\iota}<+\infty$ implies
$$\forall\;t>\mathbf{M}_{\iota}:\;\;\;\omega_{\mathbf{M}}(t)=+\infty.$$
\item[$(iii)$] If $\lim_{p\rightarrow+\infty}(M_p)^{1/p}=+\infty$, then $\omega_{\mathbf{M}}(t)<+\infty$ for all $t\ge 0$.
\end{itemize}
\end{lemma}

\emph{Note:} $(ii)$ and $(iii)$ together yield that $\lim_{p\rightarrow+\infty}(M_p)^{1/p}=+\infty$ is equivalent to $\omega_{\mathbf{M}}(t)<+\infty$ for any $t\ge 0$. In this case $\omega_{\mathbf{M}}$ is a weight function according to Definition \ref{weightfctdef}.

Let $\mathbf{M}\in\RR_{>0}^{\NN}$ be log-convex, hence $\mathbf{M}_{\iota}\in(0,+\infty]$, and we recall from \cite{regularnew} (and \cite[Chapitre I]{mandelbrojtbook}):

\begin{itemize}
\item[$(a)$] By the geometric approach presented and revisited in \cite[Sect. 3 \& 4]{regularnew}, in particular when involving the relation with the \emph{trace function} in \cite[$(3.11)$]{regularnew}, we know that $\omega_{\mathbf{M}}(t)<+\infty$ for all $t\in[0,\mathbf{M}_{\iota})$.

If $\mathbf{M}_{\iota}<+\infty$, then $\lim_{t\rightarrow\mathbf{M}_{\iota}}\omega_{\mathbf{M}}(t)$ is equal to the negative $y$-coordinate of the unique intersection point of the $y$-axis with a straight line being parallel to $t\mapsto t\cdot\mathbf{M}_{\iota}$; see \cite[Sect. 4.2, Ex. 4.3, Sect. 5.8]{regularnew}.

\item[$(b)$] When $\mathbf{M}_{\iota}<+\infty$, then in view of $(ii)$ in Lemma \ref{assofunctsectlemma} put $\omega_{\mathbf{M}}(t)=+\infty$ for all $t>\mathbf{M}_{\iota}$; see also \cite[Sect. 4.2]{regularnew}.

\item[$(c)$] Next, let us recall that by \cite[Sect. 3 \& 4]{regularnew} we have
\begin{equation}\label{reverseformula}
\forall\;p\in\NN:\;\;\;M_p=M_0\sup_{t\in(0,\mathbf{M}_{\iota})}\frac{t^p}{\exp(\omega_{\mathbf{M}}(t))};
\end{equation}
see also \cite[Chapitre I, 1.4, 1.8]{mandelbrojtbook} and \cite[Prop. 3.2]{Komatsu73}. If $\mathbf{M}$ is not log-convex, then \eqref{reverseformula} yields $M^{\on{lc}}_p$ with $\mathbf{M}^{\on{lc}}=(M^{\on{lc}}_p)_{p\in\NN}$ denoting the \emph{log-convex minorant} of $\mathbf{M}$. Moreover, even if $\mathbf{M}_{\iota}<+\infty$ by taking into account comment $(b)$ and using the convention $\frac{1}{+\infty}=0$, we can consider in \eqref{reverseformula} all $t\in(0,+\infty)$.

\item[$(d)$] Finally, in \cite[Sect. 4.1]{regularnew} we have given \eqref{reverseformula} a meaning even if $\mathbf{M}_{\iota}=0$: Here put
    $$\omega_{\mathbf{M}}(0)=0,\hspace{15pt}\omega_{\mathbf{M}}(t)=+\infty,\;\;\;\forall\;t>0.$$
    Note that for any log-convex $\mathbf{M}\in\RR_{>0}^{\NN}$ we cannot have $\mathbf{M}_{\iota}=0$ and when $\mathbf{M}\in\RR_{>0}^{\NN}$ satisfies $\mathbf{M}_{\iota}=0$, then $M^{\on{lc}}_p=0$ for any $p\in\NN_{>0}$ which contradicts $\mathbf{M}^{\on{lc}}\in\RR_{>0}^{\NN}$.
\end{itemize}

For concrete applications, when one is (only) interested in the associated function w.l.o.g. one can assume that the defining sequence is log-convex; this holds since the geometric procedure gives, in particular, that $\omega_{\mathbf{M}}=\omega_{\mathbf{M}^{\on{lc}}}$.

Next, let us recall \cite[Lemmas 2.4 \& 2.6]{regularnew}.

\begin{lemma}\label{lemma1}
Let $\mathbf{M}=(M_p)_{p\in\NN}$ be log-convex and assume that $(\mathbf{M}_{\iota}=)\lim_{p\rightarrow+\infty}(M_p)^{1/p}=\lim_{p\rightarrow+\infty}\mu_p=+\infty$. Then
\begin{equation}\label{lemma1equ}
\omega_{\mathbf{M}}(t)=0,\;\text{for}\;t\in[0,\mu_1],\;\;\;\;\omega_{\mathbf{M}}(t)=\log\left(\frac{M_0t^p}{M_p}\right)\;\text{for}\;t\in[\mu_p,\mu_{p+1}],\;p\ge 1,
\end{equation}
and when $(\mathbf{M}_{\iota}=)\lim_{p\rightarrow+\infty}(M_p)^{1/p}=\lim_{p\rightarrow+\infty}\mu_p=:C<+\infty$, then \eqref{lemma1equ} holds for all $t\in[0,C)$.
\end{lemma}

\emph{Note:} When taking into account \eqref{lemma1equ} we see that for any log-convex $\mathbf{M}$ satisfying $\mathbf{M}_{\iota}=+\infty$ the function $\omega_{\mathbf{M}}$ is a weight function in the sense of Definition \ref{weightfctdef}. In any case, Lemma \ref{lemma1} implies
$$\omega_{\mathbf{M}}(\mathbf{M}_{\iota})=\lim_{t\rightarrow\mathbf{M}_{\iota}}\omega_{\mathbf{M}}(t).$$

\begin{remark}\label{lemma1stabrem}
\emph{If $\mathbf{M}$ is log-convex with $\mathbf{M}_{\iota}=C<+\infty$, then there can exist $p_0\in\NN_{>0}$ such that $\mu_p=C$ for all $p\ge p_0$. In this case we say that the sequence of quotients eventually stabilizes and then $[\mu_p,C)=\emptyset$ for all $p\ge p_0$. However, \eqref{lemma1equ} is again valid for all $t\in[0,\mu_{p_0})=[0,\mathbf{M}_{\iota})$ and} $$\forall\;p\ge p_0:\;\;\;\omega_{\mathbf{M}}(\mathbf{M}_{\iota})=\log\left(\frac{M_0\mu_{p_0}^{p_0}}{M_{p_0}}\right)=\log\left(\frac{M_0\mu_{p}^{p}}{M_{p}}\right).$$
\end{remark}

\begin{example}\label{Gevreyweightfctexample}
\emph{For $s>0$ consider the \emph{Gevrey sequence} (of index $s$) defined by $\mathbf{G}^s=(p!^s)_{p\in\NN}$ and the corresponding quotients are given by $p^s$, $p\in\NN_{>0}$. We involve Lemma \ref{lemma1} and verify the known relation $\omega_{\mathbf{G}^s}\hyperlink{sim}{\sim}\id^{1/s}$ (which immediately gives $\gamma(\omega_{\mathbf{G}^s})=\overline{\gamma}(\omega_{\mathbf{G}^s})=s$ by Remark \ref{Gevreyweightrem}).}

\emph{First, $\omega_{\mathbf{G}^s}(t)=0$ for all $0\le t\le 1$ and when $p^s\le t\le(p+1)^s$ for some $p\in\NN_{>0}$, then $\omega_{\mathbf{G}^s}(t)=\log\left(\frac{t^p}{p!^s}\right)$ which gives}
$$\log\left(\frac{p^{sp}}{p!^s}\right)\le\omega_{\mathbf{G}^s}(t)\le\log\left(\frac{(p+1)^{sp}}{p!^s}\right).$$
\emph{Moreover, $\log\left(\frac{p^{sp}}{p!^s}\right)\ge\frac{p+1}{C}\Leftrightarrow p^{sp}\ge p!^s\exp(\frac{p+1}{C})$ for some $C\ge 1$ and all $p\in\NN_{>0}$ by \emph{Stirling's formula.} Similarly, $\log\left(\frac{(p+1)^{sp}}{p!^s}\right)\le pD\Leftrightarrow (p+1)^{sp}\le p!^s\exp(pD)$ for some $D\ge 1$ and all $p\in\NN_{>0}$. Finally, note that $\frac{p+1}{C}\ge\frac{t^{1/s}}{C}$ and $pD\le t^{1/s}D$.}
\end{example}

Let log-convex sequences $\mathbf{M},\mathbf{N}$ be given and assume that $\mathbf{M}\hyperlink{preceq}{\preceq}\mathbf{N}$; so $C_{\mathbf{M}\preceq\mathbf{N}}\in(0,+\infty)$. Then $\frac{M_p}{M_0}\le C_{\mathbf{M}\preceq\mathbf{N}}^p\frac{N_p}{N_0}$ for all $p\in\NN$ and so $\frac{N_0t^p}{N_p}\le\frac{M_0(tC_{\mathbf{M}\preceq\mathbf{N}})^p}{M_p}$ for all $t\ge 0$ and $p\in\NN$ is valid (recall \eqref{preceqalternative}). By definition this gives
\begin{equation}\label{assocfctrel}
\forall\;t\ge 0:\;\;\;\omega_{\mathbf{N}}(t)\le\omega_{\mathbf{M}}(tC_{\mathbf{M}\preceq\mathbf{N}}).
\end{equation}
Naturally one should restrict in \eqref{assocfctrel} to all $t\in[0,\mathbf{M}_{\iota}(\underline{C}_{\mathbf{M}\preceq\mathbf{N}})^{-1})$ (with the convention $\frac{1}{0}=+\infty$) since in this case \eqref{MNiotarelation} ensures $t\in[0,\mathbf{N}_{\iota})$ and hence finiteness for the left-hand side. However, even if either $\mathbf{N}_{\iota}<+\infty$ or $\mathbf{M}_{\iota}<+\infty$ formally the identity \eqref{assocfctrel} is valid for any $t\ge 0$.

\begin{itemize}
\item[$(i)$] $\mathbf{M}_{\iota}=+\infty$ and $\mathbf{M}\hyperlink{preceq}{\preceq}\mathbf{N}$ imply $\lim_{p\rightarrow+\infty}(N_p)^{1/p}=+\infty$ and by the above both sides are finite for any $t\ge 0$.

\item[$(ii)$] If $\mathbf{M}_{\iota}\in(0,+\infty)$ and $\mathbf{N}_{\iota}=+\infty$, then necessarily $\mathbf{M}\hyperlink{mtriangle}{\vartriangleleft}\mathbf{N}$ and so $\underline{C}_{\mathbf{M}\preceq\mathbf{N}}=0$. The left-hand side is finite for any $t\ge 0$, whereas the right-hand side yields the value $+\infty$ for all $t>\mathbf{M}_{\iota}(C_{\mathbf{M}\preceq\mathbf{N}})^{-1}$. Thus, formally \eqref{assocfctrel} holds again for all $t\ge 0$.

\item[$(iii)$] If $\mathbf{M}_{\iota},\mathbf{N}_{\iota}\in(0,+\infty)$, then for $t>\mathbf{N}_{\iota}$ the left-hand side is equal to $+\infty$ but by \eqref{MNiotarelation} we get that in this case $t>\mathbf{M}_{\iota}(\underline{C}_{\mathbf{M}\preceq\mathbf{N}})^{-1}\ge\mathbf{M}_{\iota}(C_{\mathbf{M}\preceq\mathbf{N}})^{-1}$ holds; therefore also the right-hand side is equal to $+\infty$ and again \eqref{assocfctrel} is formally valid for any $t\ge 0$. Note that here $\underline{C}_{\mathbf{M}\preceq\mathbf{N}}\in(0,+\infty)$.
\end{itemize}

\section{Generalized lower Legendre conjugate}\label{arblowertrafosect}
We introduce and study now the first generalized conjugate defined on the set of all weight functions.

\subsection{Preliminaries}
Let $\sigma$, $\tau$ be given weight functions and put
\begin{equation}\label{lowertransformgen}
\sigma\check{\star}\tau(t):=\inf_{s>0}\{\sigma(s)+\tau(t/s)\},\;\;\;t\in[0,+\infty).
\end{equation}

\begin{remark}\label{lowerLegendrepoint0rem}
\emph{By $\lim_{t\rightarrow+\infty}\tau(t)=+\infty$ one can use the convention $\tau(+\infty)=+\infty$ and thus consider in the definition all $s\ge 0$ since $s=0$ is then not effecting the value of $\sigma\check{\star}\tau(t)$ for any arbitrary $t\in(0,+\infty)$. And for $t=0$ we can either consider $\inf_{s\ge 0}\{\sigma(s)+\tau(0)\}$ directly or extend \eqref{lowertransformgen} to $s=0$ using the convention $\frac{0}{0}:=0$ in order to get $\sigma\check{\star}\tau(0)=\sigma(0)+\tau(0)(\ge 0)$ since $\sigma$ is non-decreasing. If $\sigma(0)\neq\lim_{t\rightarrow 0^{+}}\sigma(t)$, and this situation can occur since continuity is not assumed for being a weight function, then $\inf_{s>0}\{\sigma(s)+\tau(0)\}=\lim_{t\rightarrow 0^{+}}\sigma(t)+\tau(0)>\sigma(0)+\tau(0)$. In view of this inequality, the previous comments and since $\sigma\check{\star}\tau(0)=\sigma(0)+\tau(0)$ purely depends on $\sigma(0),\tau(0)$, this value of $\sigma\check{\star}\tau(0)$ is more natural and it automatically provides commutativity of $\check{\star}$ at $0$; see $(c)$ in Lemma \ref{wedgeproplemma}. Note that continuity at $0$ is sufficient to avoid this technical issue and (global) continuity of the weight functions under consideration is a mild extra condition which is natural for concrete examples; especially this is clear when treating associated weight functions.}
\end{remark}

We gather some more immediate consequences:

\begin{lemma}\label{wedgeproplemma}
Let $\sigma$, $\tau$ be weight functions.
\begin{itemize}
\item[$(a)$] It holds that
\begin{equation}\label{wedgeproplemmaequ}
\forall\;t,u\in[0,+\infty):\;\;\;\sigma\check{\star}\tau(tu)\le\min\{\sigma(t)+\tau(u),\sigma(u)+\tau(t)\},
\end{equation}
which gives $\tau,\sigma\hyperlink{ompreceq}{\preceq}\sigma\check{\star}\tau$;

\item[$(b)$] $\sigma\check{\star}\tau$ is a weight function;

\item[$(c)$] $\check{\star}$ is commutative.
\end{itemize}
\end{lemma}

\demo{Proof}
$(a)$ Since $\sigma\check{\star}\tau(0)=\sigma(0)+\tau(0)$ the estimate is valid if $t=0$ or $u=0$ since both $\sigma$ and $\tau$ are non-decreasing. For $t,u>0$ take $s=u$ or $s=t$ in \eqref{lowertransformgen} and hence \eqref{wedgeproplemmaequ} is verified.

Set (e.g.) $u=1$ in \eqref{wedgeproplemmaequ} and since $\sigma,\tau$ are weight functions one has $\sigma\check{\star}\tau(t)=O(\tau(t))$ and $\sigma\check{\star}\tau(t)=O(\sigma(t))$ as $t\rightarrow+\infty$; i.e. $\tau,\sigma\hyperlink{ompreceq}{\preceq}\sigma\check{\star}\tau$.\vspace{6pt}

$(b)$ By definition $\sigma\check{\star}\tau$ is non-decreasing and also $\lim_{t\rightarrow+\infty}\sigma\check{\star}\tau(t)=+\infty$ since both $\sigma$ and $\tau$ are weight functions. $(a)$ and Remark \ref{lowerLegendrepoint0rem} also give that $0\le\sigma\check{\star}\tau(t)<+\infty$ for any $t\in[0,+\infty)$.\vspace{6pt}

$(c)$ $\sigma\check{\star}\tau(0)=\tau\check{\star}\sigma(0)$ holds by Remark \ref{lowerLegendrepoint0rem}. For $t>0$ arbitrary and fixed we get with $u:=\frac{t}{s}$ that
$$\sigma\check{\star}\tau(t)=\inf_{s>0}\{\sigma(s)+\tau(t/s)\}=\inf_{u>0}\{\sigma(t/u)+\tau(u)\}=\tau\check{\star}\sigma(t).$$
Consequently, $\check{\star}$ is commutative.
\qed\enddemo

\begin{example}\label{logexamplelower}
\emph{Consider the weight function $\log_{+}(t):=\max\{0,\log(t)\}$ for $t>0$ and $\log_{+}(0):=0$; i.e. $\log_{+}(t)=0$ if $t\in[0,1)$ and $\log_{+}(t)=\log(t)$ if $t\ge 1$.}

\emph{We verify $\log_{+}\check{\star}\log_{+}=\log_{+}$ and distinguish:}

\emph{Let $t\in[0,1]$. Then for all $0<u\le 1$ one has $\log_{+}(u)+\log_{+}(t/u)=\log_{+}(t/u)$ and $\log_{+}(t/u)=0$ for all $u\ge t$. Summarizing, since both summands are non-negative and when choosing (e.g.) $u:=1$ we get $\log_{+}\check{\star}\log_{+}(t)=0=\log_{+}(t)$.}

\emph{Let $t>1$. Note that $u\mapsto\log_{+}(t/u)$ is non-increasing on $(0,1]$ and so the minimum of $u\mapsto\log_{+}(u)+\log_{+}(t/u)$ is equal to $\log_{+}(t)=\log(t)$ on $(0,1]$. For all $u\in(1,t)$ one has $\log_{+}(u)+\log_{+}(t/u)=\log(u)+\log(t/u)=\log(t)$ and finally, for all $u\in[t,+\infty)$, one has $\log_{+}(u)+\log_{+}(t/u)=\log(u)$ which is non-decreasing on this interval and so the minimum is attained at $u=t$ with the value $\log(t)$. Summarizing, for all values $t$ under consideration the infimum in the definition yields the value $\log(t)=\log_{+}(t)$.}
\end{example}

A special but important case is the following situation:
\begin{equation}\label{lowerenvelope}
\sigma\check{\star}\id(t)=\inf_{s>0}\{\sigma(s)+t/s\}=\inf_{u>0}\{\sigma(1/u)+tu\}=:(\sigma^{\iota})_{\star}(t),
\end{equation}
with $h_{\star}(t):=\inf_{u>0}\{h(u)+tu\}$ denoting the known \emph{lower Legendre conjugate (or envelope);} see \cite[Sect. 2.5]{index} for explanations and citations. More generally, for any $\alpha>0$ and $t\ge 0$ one has
$$\sigma\check{\star}\id^{1/\alpha}(t)=\inf_{s>0}\{\sigma(s)+(t/s)^{1/\alpha}\}=\inf_{u>0}\{\sigma\left(\left(\frac{1}{u}\right)^{\alpha}\right)+t^{1/\alpha}u\}
=((\sigma^{\iota})^{\alpha})_{\star}(t^{1/\alpha})=(((\sigma^{\iota})^{\alpha})_{\star})^{1/\alpha}(t);$$
hence
\begin{equation}\label{lowerLegendregeneral}
\forall\;\alpha>0\;\forall\;t\ge 0:\;\;\;\sigma\check{\star}\id^{1/\alpha}(t)=(((\sigma^{\iota})^{\alpha})_{\star})^{1/\alpha}(t).
\end{equation}

Next we verify that $\check{\star}$ preserves relations \hyperlink{ompreceq}{$\preceq$} and \hyperlink{omvartriangle}{$\vartriangleleft$} between weight functions.

\begin{lemma}\label{lowertransformrelationlemma}
Let $\sigma$, $\sigma_1$, $\tau$, and $\tau_1$ be weight functions.
\begin{itemize}
\item[$(i)$] If $\sigma\hyperlink{ompreceq}{\preceq}\sigma_1$ and $\tau\hyperlink{ompreceq}{\preceq}\tau_1$, then $\sigma\check{\star}\tau\hyperlink{ompreceq}{\preceq}\sigma_1\check{\star}\tau_1$ is satisfied.

Consequently, if $\sigma\hyperlink{sim}{\sim}\sigma_1$ and $\tau\hyperlink{sim}{\sim}\tau_1$, then $\sigma\check{\star}\tau\hyperlink{sim}{\sim}\sigma_1\check{\star}\tau_1$.

\item[$(ii)$] If $\sigma\hyperlink{omvartriangle}{\vartriangleleft}\sigma_1$ and $\tau\hyperlink{omvartriangle}{\vartriangleleft}\tau_1$, then $\sigma\check{\star}\tau\hyperlink{omvartriangle}{\vartriangleleft}\sigma_1\check{\star}\tau_1$ is satisfied.
\end{itemize}
\end{lemma}

\demo{Proof}
$(i)$ By assumption we can find $C_1,D_1,C_2,D_2\ge 1$ such that $\sigma_1(t)\le C_1\sigma(t)+D_1$ and $\tau_1(t)\le C_2\tau(t)+D_2$ for all $t\ge 0$ and so with $C:=\max\{C_1,C_2\}$:
$$\forall\;t\ge 0\;\forall\;s>0:\;\;\;\sigma_1(s)+\tau_1(t/s)\le C_1\sigma(s)+C_2\tau(t/s)+D_1+D_2\le C(\sigma(s)+\tau(t/s))+D_1+D_2.$$
By taking the infimum over all $s>0$ (see \eqref{lowertransformgen}) and since $\sigma\check{\star}\tau$ is also a weight function we have verified $\sigma\check{\star}\tau\hyperlink{ompreceq}{\preceq}\sigma_1\check{\star}\tau_1$.

$(ii)$ Here, even for any $0<c<1$ there exists $D_c\ge 1$ such that $\sigma_1(t)\le c\sigma(t)+D_c$ and $\tau_1(t)\le c\tau(t)+D_c$ for all $t\ge 0$ and the conclusion follows analogously.
\qed\enddemo

\subsection{Main results on the growth indices}
The next main statement shows how $\check{\star}$ is modifying the growth indices from Section \ref{growthindexsect}. For the sake of completeness we mention that in \cite[Rem. 5.11, p. 116]{modgrowthstrangeII} one should read in the citation Thm. 3.5 instead of Thm. 3.4.

\begin{theorem}\label{lowertransformindexthm}
Let $\sigma$, $\tau$ be weight functions.
\begin{itemize}
\item[$(i)$] If $\gamma(\sigma),\gamma(\tau)>0$, then
$$\gamma(\sigma)+\gamma(\tau)\le\gamma(\sigma\check{\star}\tau).$$

\item[$(ii)$] If $\overline{\gamma}(\sigma),\overline{\gamma}(\tau)<+\infty$, then
$$\overline{\gamma}(\sigma\check{\star}\tau)\le\overline{\gamma}(\sigma)+\overline{\gamma}(\tau).$$
\end{itemize}
\end{theorem}

\demo{Proof}
$(i)$ Let $0<\gamma<\gamma(\sigma)$, $0<\delta<\gamma(\tau)$, then by definition we get
$$\exists\;K_1>1\;\exists\;\epsilon_1\in(0,1)\;\exists\;H_1>1\;\forall\;t\ge 0:\;\;\;\sigma(K_1^{\gamma}t)\le K_1^{1-\epsilon_1}\sigma(t)+H_1,$$
and the analogous estimate holds for $\tau$ with constants $K_2,H_2>1$ and some $\epsilon_2\in(0,1)$. By $(i)$ in Lemma \ref{bigKremark} we can assume w.l.o.g. that $K_1=K_2$ and write $K$ for this constant. We also set $H:=\max\{H_1,H_2\}$ and $\epsilon:=\min\{\epsilon_1,\epsilon_2\}$ and by commutativity, see $(c)$ in Lemma \ref{wedgeproplemma}, we estimate as follows for all $t>0$:
\begin{align*}
\sigma\check{\star}\tau(K^{\gamma+\delta}t)&=\inf_{s>0}\{\sigma(s)+\tau((K^{\gamma+\delta}t)/s)\}=\inf_{s>0}\{\sigma(s)+\tau(K^{\gamma}(K^{\delta}t)/s)\}
\\&
\underbrace{=}_{u=K^{\delta}t/s}\inf_{u>0}\{\sigma(K^{\delta}t/u)+\tau(K^{\gamma}u)\}\le\inf_{u>0}\{K^{1-\epsilon}\sigma(t/u)+K^{1-\epsilon}\tau(u)\}+2H
\\&
=K^{1-\epsilon}\inf_{u>0}\{\sigma(t/u)+\tau(u)\}+2H=K^{1-\epsilon}\tau\check{\star}\sigma(t)+2H=K^{1-\epsilon}\sigma\check{\star}\tau(t)+2H.
\end{align*}
By definition this verifies $(P_{\sigma\check{\star}\tau,\gamma+\delta})$ and so $\gamma+\delta<\gamma(\sigma\check{\star}\tau)$ by \cite[Thm. 2.11]{index}. Since $\gamma<\gamma(\sigma)$ and $\delta<\gamma(\tau)$ can be chosen arbitrary we have shown $\gamma(\sigma)+\gamma(\tau)\le\gamma(\sigma\check{\star}\tau)$.\vspace{6pt}

$(ii)$ Let $\gamma>\overline{\gamma}(\sigma)$ and $\delta>\overline{\gamma}(\tau)$, then by definition
$$\exists\;A_1>1\;\exists\;\epsilon_1\in(0,1)\;\exists\;H_1\ge 1\;\forall\;t\ge 0:\;\;\;\sigma(A_1^{\gamma}t)\ge A_1^{1+\epsilon_1}\sigma(t)-H_1,$$
and
$$\exists\;A_2>1\;\exists\;\epsilon_2\in(0,1)\;\exists\;H_2\ge 1\;\forall\;t\ge 0:\;\;\;\tau(A_2^{\delta}t)\ge A_2^{1+\epsilon_2}\tau(t)-H_2.$$
By $(ii)$ in Lemma \ref{bigKremark} we can assume w.l.o.g. that $A_1=A_2$ and simply write $A$ for this constant. We set $H:=\max\{H_1,H_2\}$, $\epsilon:=\min\{\epsilon_1,\epsilon_2\}$ and estimate as follows for all $t>0$:
\begin{align*}
\sigma\check{\star}\tau(A^{\gamma+\delta}t)&=\inf_{u>0}\{\sigma((A^{\delta}t)/u)+\tau(A^{\gamma}u)\}\ge\inf_{u>0}\{A^{1+\epsilon}\sigma(t/u)+A^{1+\epsilon}\tau(u)\}-2H
\\&
=A^{1+\epsilon}\inf_{u>0}\{\sigma(t/u)+\tau(u)\}-2H=A^{1+\epsilon}\tau\check{\star}\sigma(t)-2H=A^{1+\epsilon}\sigma\check{\star}\tau(t)-2H.
\end{align*}
Therefore, $(\overline{P}_{\sigma\check{\star}\tau,\gamma+\delta})$ is verified and in view of \cite[Thm. 2.16, $(7)$]{index} this implies $\gamma+\delta>\overline{\gamma}(\sigma\check{\star}\tau)$. Since $\gamma>\overline{\gamma}(\sigma)$ and $\delta>\overline{\gamma}(\tau)$ are arbitrary we have shown $\overline{\gamma}(\sigma)+\overline{\gamma}(\tau)\ge\overline{\gamma}(\sigma\check{\star}\tau)$.
\qed\enddemo

Let us apply Theorem \ref{lowertransformindexthm} to some special cases:

\begin{corollary}\label{lowertransformindexcor}
Let $\sigma$, $\tau$ be weight functions.
\begin{itemize}
\item[$(i)$] When $\gamma(\sigma)>0$, then
$$\forall\;\alpha>0:\;\;\;\gamma(\sigma)+\alpha\le\gamma((((\sigma^{\iota})^{\alpha})_{\star})^{1/\alpha})$$
and $\overline{\gamma}(\sigma)<+\infty$ implies
$$\forall\;\alpha>0:\;\;\;\overline{\gamma}((((\sigma^{\iota})^{\alpha})_{\star})^{1/\alpha})\le\overline{\gamma}(\sigma)+\alpha.$$

\item[$(ii)$] When $0<\gamma(\sigma)=\overline{\gamma}(\sigma)<+\infty$, $0<\gamma(\tau)=\overline{\gamma}(\tau)<+\infty$, then
$$\gamma(\sigma)+\gamma(\tau)=\gamma(\sigma\check{\star}\tau)=\overline{\gamma}(\sigma\check{\star}\tau)=\overline{\gamma}(\sigma)+\overline{\gamma}(\tau).$$
\end{itemize}
\end{corollary}

\emph{Note:} $(i)$ with $\alpha=1$ should be compared with \cite[Prop. 2.22]{index}; there equality has been established in \cite[$(10)$]{index} under the additional assumption that $\sigma^{\iota}$ is equivalent to its \emph{largest convex minorant.} Moreover, when $0<\gamma(\sigma)=\overline{\gamma}(\sigma)<+\infty$, then
\begin{equation}\label{lowertransformindexcorequ}
\forall\;\alpha>0:\;\;\;\gamma(\sigma)+\alpha=\gamma((((\sigma^{\iota})^{\alpha})_{\star})^{1/\alpha})=\overline{\gamma}((((\sigma^{\iota})^{\alpha})_{\star})^{1/\alpha}),
\end{equation}
thus, in particular, the following identity:
\begin{equation}\label{Gevreyadd}
\forall\;\alpha,\beta>0:\;\;\;\alpha+\beta=\gamma(\id^{1/\alpha}\check{\star}\id^{1/\beta})=\overline{\gamma}(\id^{1/\alpha}\check{\star}\id^{1/\beta}).
\end{equation}

\demo{Proof}
$(i)$ By using Remark \ref{Gevreyweightrem} and Theorem \ref{lowertransformindexthm} applied to $\sigma$ and $\tau:=\id^{1/\alpha}$, $\alpha>0$ arbitrary, we see
$$\gamma(\sigma)+\alpha\le\gamma(\sigma\check{\star}\id^{1/\alpha})\le\overline{\gamma}(\sigma\check{\star}\id^{1/\alpha})\le\overline{\gamma}(\sigma)+\alpha.$$
And then \eqref{lowerLegendregeneral} yields the conclusion.

$(ii)$ is an immediate consequence of Theorem \ref{lowertransformindexthm}.
\qed\enddemo

\emph{Note:} Both indices $\gamma(\cdot)$ and $\overline{\gamma}(\cdot)$ are preserved under the relation \hyperlink{sim}{$\sim$}; see \cite[Rem. 2.12, Thm. 2.11 \& 2.16]{index}. Hence, in view of Lemma \ref{lowertransformrelationlemma} all inequalities and identities listed in Theorem \ref{lowertransformindexthm} and Corollary \ref{lowertransformindexcor} are preserved under equivalence of weight functions, too.

In view of Example \ref{Gevreyweightfctexample}, identity \eqref{Gevreyadd} is then equivalent to
$$\forall\;\alpha,\beta>0:\;\;\;\alpha+\beta=\gamma(\omega_{\mathbf{G}^{\alpha}}\check{\star}\omega_{\mathbf{G}^{\beta}})=\overline{\gamma}(\omega_{\mathbf{G}^{\alpha}}\check{\star}\omega_{\mathbf{G}^{\beta}}).$$

\begin{remark}\label{indexnaturallowerrem}
\emph{The assumptions on the indices in Theorem \ref{lowertransformindexthm} and Corollary \ref{lowertransformindexcor} are somehow natural in this context; however we remark:}

\begin{itemize}
\item[$(*)$] \emph{In $(i)$ in Theorem \ref{lowertransformindexthm}, if either $\gamma(\sigma)=+\infty$ or $\gamma(\tau)=+\infty$, then by the given proof $\gamma(\sigma\check{\star}\tau)=+\infty$ as well. And, indeed, $(i)$ and the given proof also holds when either $\gamma(\sigma)=0$ or $\gamma(\tau)=0$: Take $\gamma\le 0$ resp. $\delta\le 0$ and recall that in this case $(P_{\sigma,\gamma})$ resp. $(P_{\tau,\delta})$ is trivial. And note that $\gamma(\sigma\check{\star}\tau)\ge 0$ because $\sigma\check{\star}\tau$ is a weight function.}

\item[$(*)$] \emph{In $(ii)$ in Theorem \ref{lowertransformindexthm}, when either $\overline{\gamma}(\sigma)=+\infty$ or $\overline{\gamma}(\tau)=+\infty$, then we cannot follow the proof but the desired inequality is (formally) trivial.}

\item[$(*)$] \emph{Thus $(i)$ in Corollary \ref{lowertransformindexcor} still holds when $\gamma(\sigma)=0$ or $\overline{\gamma}(\sigma)=+\infty$ and in the second situation the desired estimate is formally trivial.}

\item[$(*)$] \emph{Finally, also the equalities stated in $(ii)$ in Corollary \ref{lowertransformindexcor} are (formally) valid for the remaining cases.}
\end{itemize}
\end{remark}

\section{Generalized upper Legendre conjugate}\label{arbuppertrafosect}
Now we proceed with the second conjugate; this approach is more technical.

\subsection{Preliminaries}\label{arbuppertrafosectprelim}
Let $\sigma$, $\tau$ be given weight functions and put
\begin{equation}\label{uppertransformgen}
\sigma\widehat{\star}\tau(t):=\sup_{s\ge 0}\{\sigma(s)-\tau(s/t)\},\;\;\;t\in(0,+\infty).
\end{equation}
By definition $\sigma\widehat{\star}\tau(t)\ge\sigma(0)-\tau(0)$ (formally) for all $t\in(0,+\infty)$ but in general it is not clear that $0\le\sigma\widehat{\star}\tau(t)<+\infty$ for some/any $t\in[0,+\infty)$ and so $\sigma\widehat{\star}\tau$ is \emph{not automatically a weight function} according to Definition \ref{weightfctdef}. To see this, and this example has been suggested by one of the referees, we consider e.g. $\sigma(t):=t^2$ and $\tau(t):=t$ and then even $\sigma\widehat{\star}\tau(t)=+\infty$ for all $t\in(0,+\infty)$. Moreover, the definition resp. the value of $\sigma\widehat{\star}\tau(0)$ is not clear. We start with the following observation.

\begin{lemma}\label{widehatproplemma}
Let $\sigma$, $\tau$ be weight functions.
\begin{itemize}
\item[$(a)$] $\sigma\widehat{\star}\tau$ is non-decreasing, $\lim_{t\rightarrow+\infty}\sigma\widehat{\star}\tau(t)=+\infty$ and (formally)
\begin{equation}\label{widehatproplemmaequ}
\forall\;t,u>0:\;\;\;\sigma\widehat{\star}\tau(tu)\ge\max\{\sigma(t)-\tau(1/u),\sigma(u)-\tau(1/t)\},
\end{equation}
which gives $\sigma\widehat{\star}\tau\hyperlink{ompreceq}{\preceq}\sigma$.

\item[$(b)$] The conventions $\frac{0}{0}=0$, $\frac{1}{0}=+\infty$, $\tau(+\infty)=+\infty$, imply $\sigma\widehat{\star}\tau(0)=\sigma(0)-\tau(0)$.
\end{itemize}
\end{lemma}

\demo{Proof}
$(a)$ Since $\tau$ is non-decreasing also $\sigma\widehat{\star}\tau$ is non-decreasing and take $s:=t$ resp. $s:=u$ in \eqref{uppertransformgen} to obtain \eqref{widehatproplemmaequ}. Choosing $u:=1$ this estimate verifies $\lim_{t\rightarrow+\infty}\sigma\widehat{\star}\tau(t)=+\infty$ and $\sigma(t)=O(\sigma\widehat{\star}\tau(t))$ as $t\rightarrow+\infty$.\vspace{6pt}

$(b)$ This is immediate by the conventions which give $\sigma(0)-\tau(0/0)=\sigma(0)-\tau(0)$ and $\sigma(s)-\tau(s/0)=-\infty$ for any $s>0$.
\qed\enddemo

Lemma \ref{widehatproplemma} implies the next consequence:

\begin{lemma}\label{uppertransformnonnegative}
Let $\sigma$ and $\tau$ be weight functions. In order to ensure (formally) $0\le\sigma\widehat{\star}\tau(t)$ for all $t\in[0,+\infty)$ it suffices to assume $\tau(0)=0$.
\end{lemma}

\demo{Proof}
By definition (set $s:=0$) one has $\sigma\widehat{\star}\tau(t)\ge\sigma(0)-\tau(0)=\sigma\widehat{\star}\tau(0)$ for any $t>0$; see $(b)$ in Lemma \ref{widehatproplemma}. If $\tau(0)=0$, then the conclusion follows since $\sigma(0)\ge 0$.
\qed\enddemo

\emph{Note:} $(b)$ in Lemma \ref{widehatproplemma} is in accordance with Remark \ref{lowerLegendrepoint0rem} for the lower conjugate and in view of the fact that $\sigma\widehat{\star}\tau$ is non-decreasing and since (formally) $\sigma\widehat{\star}\tau(t)\ge\sigma(0)-\tau(0)$ for all $t\in(0,+\infty)$, the equality $\sigma\widehat{\star}\tau(0)=\sigma(0)-\tau(0)$ is somehow natural. However, note that the expected equality $\sigma(0)-\tau(0)=\lim_{t\rightarrow 0^{+}}\sigma\widehat{\star}\tau(t)$ fails in general as it can be seen when considering e.g. $\sigma(t):=t^2$ and $\tau(t):=t$; recall the comments above. This example also motivates the study in which cases we can ensure $\sigma\widehat{\star}\tau(t)<+\infty$, at least for all sufficiently small $t>0$. But by Corollary \ref{inversethm1cor}, even if $\sigma\widehat{\star}\tau$ is a weight function, the aforementioned expected equality at $0$ fails in general. And similarly as in Remark \ref{lowerLegendrepoint0rem} also $\sigma(0)-\tau(0)=\lim_{t\rightarrow 0^{+}}\sigma(t)-\tau(t)$ is violated when considering e.g. any globally continuous weight function $\tau$ and $\sigma$ arbitrary such that $\sigma(0)\neq\lim_{t\rightarrow 0^{+}}\sigma(t)$ (cf. again with Corollary \ref{inversethm1cor}).

\subsection{On the well-definedness of the generalized upper Legendre conjugate}
We start with the first main technical result:

\begin{lemma}\label{uppertransformfinite}
Let $\sigma$ and $\tau$ be weight functions and let $t_0\in(0,+\infty]$. Consider the following assertions:
\begin{itemize}
\item[$(i)$] We have
\begin{equation}\label{uppertransformfiniteequ}
\sup_{0<t<t_0}\limsup_{u\rightarrow+\infty}\frac{\sigma(tu)}{\tau(u)}<1.
\end{equation}

\item[$(ii)$] $\sigma\widehat{\star}\tau(t)<+\infty$ holds for all $t\in(0,t_0)$; i.e.
\begin{equation}\label{uppetransfromwelldefequ}
\forall\;t\in(0,t_0)\;\exists\;D_t>0\;\forall\;s\ge 0:\;\;\;\sigma(s)-\tau(s/t)\le D_t.
\end{equation}

\item[$(iii)$] We have
\begin{equation}\label{uppertransformfiniteequweak}
\sup_{0<t<t_0}\limsup_{u\rightarrow+\infty}\frac{\sigma(tu)}{\tau(u)}\le 1.
\end{equation}
\end{itemize}
Then $(i)\Rightarrow(ii)\Rightarrow(iii)$ is valid.
\begin{itemize}
\item[$(iv)$] On the other hand, if
\begin{equation}\label{uppertransformfiniteequweakviol}
\exists\;t_0\in(0,+\infty):\;\;\;\liminf_{u\rightarrow+\infty}\frac{\sigma(t_0u)}{\tau(u)}>1,
\end{equation}
then $\sigma\widehat{\star}\tau(t)=+\infty$ for all $t\ge t_0$.
\end{itemize}
\end{lemma}

For the quotients in the $\limsup$-expressions note that $\tau$ is non-decreasing and $\lim_{u\rightarrow+\infty}\tau(u)=+\infty$; so, in particular, $\tau(u)>0$ for all $u$ sufficiently large. On the other hand note that in general it is not clear if any of the conditions \eqref{uppertransformfiniteequ}, \eqref{uppertransformfiniteequweak}, \eqref{uppertransformfiniteequweakviol} is preserved under equivalence of weight functions.

\demo{Proof}
$(i)\Rightarrow(ii)$ By assumption $\limsup_{u\rightarrow+\infty}\frac{\sigma(tu)}{\tau(u)}<1$ for all $t\in(0,t_0)$, therefore
$$\exists\;\epsilon>0\;\forall\;t\in(0,t_0)\;\exists\;u_t>0\;\forall\;u\ge u_t:\;\;\;\sigma(tu)\le(1-\epsilon)\tau(u),$$
equivalently
$$\exists\;\epsilon>0\;\forall\;t\in(0,t_0)\;\exists\;D_t>0\;\forall\;u\ge 0:\;\;\;\sigma(tu)\le(1-\epsilon)\tau(u)+D_t.$$
Set $s:=tu$, then this estimate turns into
$$\exists\;\epsilon>0\;\forall\;t\in(0,t_0)\;\exists\;D_t>0\;\forall\;s\ge 0:\;\;\;\sigma(s)\le(1-\epsilon)\tau(s/t)+D_t<\tau(s/t)+D_t,$$
hence by definition $\sigma\widehat{\star}\tau(t)\le D_t<+\infty$ for all $t\in(0,t_0)$ follows immediately.\vspace{6pt}

$(ii)\Rightarrow(iii)$ In \eqref{uppetransfromwelldefequ} put $u:=\frac{s}{t}$ and so
$$\forall\;t\in(0,t_0)\;\exists\;D_t>0\;\forall\;u\ge 0:\;\;\;\sigma(ut)\le\tau(u)+D_t.$$
Thus \eqref{uppertransformfiniteequweak} follows by taking into account $\lim_{u\rightarrow+\infty}\tau(u)=+\infty$.\vspace{6pt}

$(iv)$ It holds that \eqref{uppertransformfiniteequweakviol} implies $\liminf_{u\rightarrow+\infty}\frac{\sigma(tu)}{\tau(u)}>1$ for all $t\ge t_0$ and by assumption with $A:=\liminf_{u\rightarrow+\infty}\frac{\sigma(t_0u)}{\tau(u)}\in(1,+\infty]$ we get
$$\exists\;\epsilon_1>0\;\forall\;\epsilon\in(\epsilon_1,A-1)\;\exists\;u_{\epsilon}>0:\;\;\;\sigma(t_0u_{\epsilon})>(1+\epsilon)\tau(u_{\epsilon})>(1+\epsilon_1)\tau(u_{\epsilon}).$$
Hence $\sigma(s_{\epsilon})-\tau(s_{\epsilon}/t_0)>\epsilon_1\tau(s_{\epsilon}/t_0)$ with $s_{\epsilon}:=t_0u_{\epsilon}$ which immediately gives $\sigma\widehat{\star}\tau(t_0)=+\infty$: As $\epsilon\rightarrow A-1$ one has $s_{\epsilon}\rightarrow+\infty\Leftrightarrow u_{\epsilon}\rightarrow+\infty$ and so $\tau(s_{\epsilon}/t_0)\rightarrow+\infty$.
\qed\enddemo

\emph{Note:}

\begin{itemize}
\item[$(i)$] In general it is not clear if $(iii)\Rightarrow(ii)$ in Lemma \ref{uppertransformfinite} holds; for this note that \eqref{uppertransformfiniteequweak} precisely means
    $$\forall\;t\in(0,t_0)\;\forall\;\epsilon>0\;\exists\;D_{t,\epsilon}>0\;\forall\;u\ge 0:\;\;\;\sigma(tu)\le(1+\epsilon)\tau(u)+D_{t,\epsilon},$$
    but which does not imply necessarily \eqref{uppetransfromwelldefequ}.

\item[$(ii)$] If
\begin{equation}\label{uppertransformfiniteequ1}
\forall\;t>0:\;\;\;\liminf_{u\rightarrow+\infty}\frac{\sigma(tu)}{\tau(u)}>1,
\end{equation}
i.e. \eqref{uppertransformfiniteequweakviol} for any $t_0\in(0,+\infty)$, then
$$\forall\;t>0\;\exists\;\epsilon_t>0\;\exists\;u_t>0\;\forall\;u\ge u_t:\;\;\;\sigma(tu)\ge(1+\epsilon_t)\tau(u),$$
and hence $\sigma\widehat{\star}\tau(t)=+\infty$ for any $t\in(0,+\infty)$. As an example take, e.g., as mentioned above $\sigma(t)=t^2$ and $\tau(t)=t$.
\end{itemize}

\begin{remark}\label{slowlyvarrem}
\emph{Lemma \ref{uppertransformfinite} for $t_0=+\infty$ implies that, under appropriate (and mild) assumptions on $\sigma$ and $\tau$, we get $\sigma\widehat{\star}\tau(t)<+\infty$ for all $t\in(0,+\infty)$; i.e. $\sigma\widehat{\star}\tau$ is well defined.}

\emph{If $t_0<+\infty$, then $\sigma\widehat{\star}\tau$ is well-defined for all $t<t_0$. If in addition \eqref{uppertransformfiniteequweakviol} holds then $\sigma\widehat{\star}\tau(t)=+\infty$ for all $t\ge t_0$ and this situation is becoming relevant for some cases in the weight sequence setting studied in Section \ref{uppertrafosect}.}

\emph{Moreover, for $\sigma\widehat{\star}\sigma$ we have the following comments for which we use that $\sigma$ is non-decreasing:}

\begin{itemize}
\item[$(i)$] \emph{The $\limsup$-expression satisfies the estimate in \eqref{uppertransformfiniteequweak} for all $0<t\le 1$.}

\emph{\eqref{uppertransformfiniteequ} has to fail whenever $t_0>1$: Consider any $t\in(1,t_0)$ and then $\frac{\sigma(tu)}{\sigma(u)}\ge 1$ (for all large $u$ such that $\sigma(u)>0$).}

\item[$(ii)$] \emph{Since $\sigma$ is non-decreasing, \eqref{uppertransformfiniteequweak} for $t_0=+\infty$ amounts to having $$\forall\;t>0:\;\;\;\lim_{u\rightarrow+\infty}\frac{\sigma(tu)}{\sigma(u)}=1,$$
    which precisely means that $\sigma$ is \emph{slowly varying;} see e.g. \cite[Sect. 2.2, $(3)$]{index}. Consequently, if $\sigma\widehat{\star}\sigma(t)<+\infty$ for all $t>0$, then $\sigma$ has to be slowly varying.}
\end{itemize}
\end{remark}

The next example illustrates that implication $(i)\Rightarrow(ii)$ in Lemma \ref{uppertransformfinite} is strict in general; see also comment $(i)$ in the previous Remark.

\begin{example}\label{logexampleupper}
\emph{Consider again the weight function $\log_{+}$ and let us verify $\log_{+}\widehat{\star}\log_{+}=\log_{+}$. We distinguish:}

\emph{For $t=0$ one has $\log_{+}\widehat{\star}\log_{+}(0)=\log_{+}(0)-\log_{+}(0)=0-0=0=\log_{+}(0)$.}

\emph{Let $t\in(0,1]$. Then for all $s\in[0,t)$ one has $\log_{+}(s)-\log_{+}(s/t)=0-0=0$, for $s\in[t,1)$ one has $\log_{+}(s)-\log_{+}(s/t)=0-\log(s/t)=\log(t/s)$ and so this expression attains its maximum at $s=t$ which yields the value $\log(1)=0$. And if $s\in[1,+\infty)$, then $\log_{+}(s)-\log_{+}(s/t)=\log(s)-\log(s/t)=\log(t)\le 0$. Therefore, $\log_{+}\widehat{\star}\log_{+}(t)=0=\log_{+}(t)$ is verified.}

\emph{Let $t>1$. Then for $s\in[0,1)$ one gets $\log_{+}(s)-\log_{+}(s/t)=0$, for $s\in[1,t)$ one gets $\log_{+}(s)-\log_{+}(s/t)=\log(s)\le\log(t)$ and for $s\in[t,+\infty)$ one has $\log_{+}(s)-\log_{+}(s/t)=\log(s)-\log(s/t)=\log(t)$.}

\emph{Summarizing, in this case the expression in the supremum is equal to $\log(t)=\log_{+}(t)$ which establishes the claim. However, as seen before, \eqref{uppertransformfiniteequ} fails (for all $t_0>1$) and $\log_{+}$ is obviously slowly varying.}
\end{example}

By the comments and results obtained so far, in order to ensure that $\sigma\widehat{\star}\tau$ is a weight function one shall naturally assume that $\sigma$ and $\tau$ satisfy

\begin{itemize}
\item[$(A)$] $\tau(0)=0$ and

\item[$(B)$] $\sigma\widehat{\star}\tau$ is \emph{well defined on $[0,+\infty)$} i.e. \eqref{uppetransfromwelldefequ} holds with $t_0=+\infty$ and so $\sigma\widehat{\star}\tau(t)<+\infty$ for any $t>0$.
\end{itemize}
Indeed, $(B)$ is ensured when \eqref{uppertransformfiniteequ} is satisfied for $t_0=+\infty$:
\begin{equation}\label{conditionB}
\sup_{t>0}\limsup_{u\rightarrow+\infty}\frac{\sigma(tu)}{\tau(u)}<1.
\end{equation}
For convenience, frequently we say that $\sigma\widehat{\star}\tau$ is well defined if \eqref{uppetransfromwelldefequ} holds with $t_0=+\infty$ and that \emph{$\sigma$ is subordinate to $\tau$ via $\widehat{\star}$} if \eqref{conditionB} holds.

Again, as special cases we consider
\begin{equation}\label{upperLegendre}
\sigma\widehat{\star}\id(t)=\sup_{s\ge 0}\{\sigma(s)-s/t\}=:\sigma^{\star}(1/t)=(\sigma^{\star})^{\iota}(t),\;\;\;t\in(0,+\infty),
\end{equation}
with $\sigma^{\star}(t):=\sup_{s\ge 0}\{\sigma(s)-ts\}$ denoting the known \emph{upper Legendre conjugate (or envelope);} see again \cite[Sect. 2.5]{index} for detailed explanations and citations. More generally, for any $\alpha>0$ and $t>0$ consider
\begin{align*}
\sigma\widehat{\star}\id^{1/\alpha}(t)&=\sup_{s\ge 0}\{\sigma(s)-(s/t)^{1/\alpha}\}=\sup_{u\ge 0}\{\sigma(u^{\alpha})-ut^{-1/\alpha}\}
\\&
=(\sigma^{\alpha})^{\star}(t^{-1/\alpha})=(((\sigma^{\alpha})^{\star})^{\iota})^{1/\alpha}(t).
\end{align*}
Hence, since $\id^{1/\alpha}(0)=0$ in view of $(b)$ in Lemma \ref{widehatproplemma} this identity is also valid for $t=0$ and so
\begin{equation}\label{upperLegendregeneral}
\forall\;\alpha>0\;\forall\;t\ge 0:\;\;\;\sigma\widehat{\star}\id^{1/\alpha}(t)=(((\sigma^{\alpha})^{\star})^{\iota})^{1/\alpha}(t).
\end{equation}

Note that $(A)$ for $\sigma\widehat{\star}\id^{1/\alpha}$ is trivial since $\id^{1/\alpha}(0)=0$ and Lemma \ref{uppertransformfinite} takes the following characterizing form:

\begin{corollary}\label{uppertransformfinitecor}
Let $\sigma$ be a weight function and $\alpha>0$. Then the following are equivalent:
\begin{itemize}
\item[$(i)$] We have
$$\forall\;t>0:\;\;\;\lim_{u\rightarrow+\infty}\frac{\sigma(tu)}{\id^{1/\alpha}(u)}=\lim_{u\rightarrow+\infty}\frac{\sigma(tu)}{u^{1/\alpha}}=0.$$

\item[$(ii)$] $\sigma\widehat{\star}\id^{1/\alpha}(t)<+\infty$ for all $t\in(0,+\infty)$; i.e. \eqref{uppetransfromwelldefequ} holds for $t_0=+\infty$.

\item[$(iii)$] We have
$$\sup_{t>0}\limsup_{u\rightarrow+\infty}\frac{\sigma(tu)}{\id^{1/\alpha}(u)}\le 1.$$

\item[$(iv)$] We have $\sigma(s)=o(s^{1/\alpha})$ as $s\rightarrow+\infty$; i.e. $\id^{1/\alpha}\hyperlink{omvartriangle}{\vartriangleleft}\sigma$.
\end{itemize}
\end{corollary}

\demo{Proof}
$(i)\Rightarrow(ii)$ follows by $(i)\Rightarrow(ii)$ in Lemma \ref{uppertransformfinite} applied to $\tau=\id^{1/\alpha}$ and $t_0=+\infty$ and $(ii)\Rightarrow(iii)$ is $(ii)\Rightarrow(iii)$ in Lemma \ref{uppertransformfinite}.\vspace{6pt}

$(iii)\Rightarrow(iv)$ For any fixed $t>0$ it follows that
$$\limsup_{v\rightarrow+\infty}\frac{\sigma(v)}{(v/t)^{1/\alpha}}=\limsup_{u\rightarrow+\infty}\frac{\sigma(tu)}{u^{1/\alpha}}=\limsup_{u\rightarrow+\infty}\frac{\sigma(tu)}{\id^{1/\alpha}(u)}\le 1.$$
So, when we let $t\rightarrow+\infty$ then $\sigma(v)=o(v^{1/\alpha})$ as $v\rightarrow+\infty$ is valid.\vspace{6pt}

$(iv)\Rightarrow(i)$ Let $t>0$ be arbitrary but fixed. $\sigma(s)=o(s^{1/\alpha})$ as $s\rightarrow+\infty$ implies
$$\forall\;\epsilon>0\;\exists\;D_{\epsilon}\ge 1\;\forall\;u\ge 0:\;\;\;\sigma(tu)\le\epsilon(tu)^{1/\alpha}+D_{\epsilon}.$$
Now, let $\epsilon_1>0$ be small and choose $\epsilon:=\frac{\epsilon_1}{t^{1/\alpha}}$ and thus $\limsup_{u\rightarrow+\infty}\frac{\sigma(tu)}{u^{1/\alpha}}\le\epsilon t^{1/\alpha}=\epsilon_1$. Since $\epsilon_1$ is arbitrary we have verified $\limsup_{u\rightarrow+\infty}\frac{\sigma(tu)}{u^{1/\alpha}}=0$ and since $t>0$ was arbitrary we are done.
\qed\enddemo

The previous Corollary suggests that for more particular weights one can expect in Lemma \ref{uppertransformfinite} a characterization. The next technical result relates the introduced $\limsup$-conditions with the known growth relation \eqref{smallOrelation}.

\begin{lemma}\label{om1om6lemma}
Let $\sigma$ and $\tau$ be weight functions. Consider the following assertions:
\begin{itemize}
\item[$(i)$] \eqref{uppertransformfiniteequweak} holds with $t_0=+\infty$; i.e.
$$\sup_{t>0}\limsup_{u\rightarrow+\infty}\frac{\sigma(tu)}{\tau(u)}\le 1.$$

\item[$(ii)$] The weights are related by $\tau\hyperlink{omvartriangle}{\vartriangleleft}\sigma$; i.e. $\sigma(t)=o(\tau(t))$ as $t\rightarrow+\infty$.

\item[$(iii)$] We have
\begin{equation}\label{om1om6lemmaequ}
\sup_{t>0}\lim_{u\rightarrow+\infty}\frac{\sigma(tu)}{\tau(u)}=0,
\end{equation}
which implies \eqref{conditionB} and so $\sigma$ is subordinate to $\tau$ via $\widehat{\star}$ and $\sigma\widehat{\star}\tau$ is well defined.
\end{itemize}

If either $\sigma$ or $\tau$ satisfies \hyperlink{om6}{$(\omega_6)$}, then $(i)\Rightarrow(ii)$; if either $\sigma$ or $\tau$ satisfies \hyperlink{om1}{$(\omega_1)$}, then $(ii)\Rightarrow(iii)$.
\end{lemma}



\demo{Proof}
$(i)\Rightarrow(ii)$ \eqref{uppertransformfiniteequweak} with $t_0=+\infty$ precisely means
\begin{equation}\label{om1om6lemmaequ1}
\forall\;t>0\;\forall\;\epsilon>0\;\exists\;D_{t,\epsilon}>0\;\forall\;u\ge 0:\;\;\;\sigma(tu)\le(1+\epsilon)\tau(u)+D_{t,\epsilon}.
\end{equation}
Fix now $\epsilon$ and write $D_t$ instead of $D_{t,\epsilon}$. Let $c>0$ small be given and fixed, choose $n\in\NN_{>0}$ (minimal) such that $\frac{1}{c}\le 2^n$. When $\sigma$ satisfies \hyperlink{om6}{$(\omega_6)$}, then by iterating this property $n$-times we get $\frac{1}{c}\sigma(u)\le 2^n\sigma(u)\le\sigma(Hu)+H$ for some $H\ge 1$ depending on $c$ and all $u\ge 0$. And so, by applying \eqref{om1om6lemmaequ1} to $t:=H$, we have for all $u\ge 0$ the estimate $\frac{1}{c}\sigma(u)\le(1+\epsilon)\tau(u)+H+D_H$. Note that $n$ and $H$ only depend on chosen $c$. Since $c>0$ is arbitrary, this verifies $\sigma(t)=o(\tau(t))$ as $t\rightarrow+\infty$ as desired.

When $\tau$ satisfies \hyperlink{om6}{$(\omega_6)$}, then similarly
$$\frac{1}{c}\sigma(u)\le 2^n\sigma(u)\le 2^n(1+\epsilon)\tau(uH^{-1})+2^nD_H\le(1+\epsilon)\tau(u)+(1+\epsilon)H+2^nD_H,$$
for some $H\ge 1$ depending on $c$ and all $u\ge 0$. So $\sigma(t)=o(\tau(t))$ as $t\rightarrow+\infty$ follows again.\vspace{6pt}

$(ii)\Rightarrow(iii)$ Assertion $(ii)$ precisely means
\begin{equation}\label{om1om6lemmaequ2}
\forall\;c>0\;\exists\;A_c\ge 1\;\forall\;u\ge 0:\;\;\;\sigma(u)\le c\tau(u)+A_c.
\end{equation}
Let now $t>0$ be given, arbitrary large but fixed. Choose $n\in\NN_{>0}$ (minimal) such that $t\le 2^n$. When $\sigma$ satisfies \hyperlink{om1}{$(\omega_1)$}, then by iterating this property $n$-times we get $\sigma(2^nu)\le L\sigma(u)+L$ for some $L\ge 1$ depending on $t$ and all $u\ge 0$. Let $\epsilon>0$ be arbitrary (small) but fixed and apply \eqref{om1om6lemmaequ2} to $c:=\frac{\epsilon}{L}$, depending on given $t$ and $\epsilon$, to get for all $u\ge 0$:
$$\sigma(tu)\le\sigma(2^nu)\le L\sigma(u)+L\le\epsilon\tau(u)+LA_{\epsilon L^{-1}}+L.$$
Since $\tau$ is a weight function, this estimate implies $\limsup_{u\rightarrow+\infty}\frac{\sigma(tu)}{\tau(u)}\le\epsilon$ and since $\epsilon>0$ was arbitrary this gives $\limsup_{u\rightarrow+\infty}\frac{\sigma(tu)}{\tau(u)}=0$. Since $t>0$ was arbitrary, \eqref{om1om6lemmaequ} is verified.

When $\tau$ satisfies \hyperlink{om1}{$(\omega_1)$}, then similarly by using again \eqref{om1om6lemmaequ2} for $c:=\frac{\epsilon}{L}$ one has for all $u\ge 0$:
$$\sigma(tu)\le\sigma(2^nu)\le\frac{\epsilon}{L}\tau(2^nu)+A_{\epsilon L^{-1}}\le\epsilon\tau(u)+\epsilon+A_{\epsilon L^{-1}}.$$
So again \eqref{om1om6lemmaequ} is verified.
\qed\enddemo

\begin{remark}\label{om6notslowlyvarying}
\emph{Since $\tau\hyperlink{omvartriangle}{\vartriangleleft}\tau$ cannot hold, Lemma \ref{om1om6lemma} for the special case $\sigma=\tau$ yields the immediate fact that any slowly varying weight function cannot have \hyperlink{om6}{$(\omega_6)$}; recall $(ii)$ in Remark \ref{slowlyvarrem} and Example \ref{logexampleupper}.}

\emph{This observation is consistent with \cite[Rem. 2.18]{index}; see also the statements \cite[Lemma 2.10, Thm. 2.11 \& 2.16, Cor. 2.17]{index}.}
\end{remark}

\begin{remark}\label{om1om6lemmaremnew}
\emph{In the proof of Lemma \ref{om1om6lemma} one requires growth properties for weights, here \hyperlink{om1}{$(\omega_1)$} and  \hyperlink{om6}{$(\omega_6)$}, only for \emph{one of the weights and not necessarily for both or for a particular one.} This fact is occurring in different contexts dealing with weights and their growth properties and is an universal principle; we refer to \cite[Prop. 3.22, Thm. 3.25]{weightedentireinclusion1} and \cite[Lemma 3.3, Thm. 4.4, Cor. 4.5, Rem. 4.6]{orliczpaper}. Indeed, this observation is useful in the proof of $(ii)$ in Theorem \ref{upperconjequivcor}.}
\end{remark}

\subsection{Main results on the growth indices}
We study how $\widehat{\star}$ modifies the indices from Section \ref{growthindexsect}.

\begin{theorem}\label{uppertransformindexthm}
Let $\sigma$, $\tau$ be weight functions such that $\tau(0)=0$ and $\sigma\widehat{\star}\tau$ is well defined.
\begin{itemize}
\item[$(i)$] If $\gamma(\sigma)>0$ and $\overline{\gamma}(\tau)<+\infty$, then
\begin{equation}\label{uppertransformindexthmequ}
\gamma(\sigma)\le\gamma(\sigma\widehat{\star}\tau)+\overline{\gamma}(\tau).
\end{equation}
\item[$(ii)$] If $0<\gamma(\tau)\le\overline{\gamma}(\sigma)<+\infty$, then
\begin{equation}\label{uppertransformindexthmequ1}
\overline{\gamma}(\sigma\widehat{\star}\tau)+\gamma(\tau)\le\overline{\gamma}(\sigma).
\end{equation}
\end{itemize}
\end{theorem}

\demo{Proof}
$(i)$ Let $0<\gamma<\gamma(\sigma)$, $\delta>\overline{\gamma}(\tau)$, then by definition
$$\exists\;K>1\;\exists\;\epsilon_1\in(0,1)\;\exists\;H_1\ge 1\;\forall\;t\ge 0:\;\;\;\sigma(K^{\gamma}t)\le K^{1-\epsilon_1}\sigma(t)+H_1,$$
and
$$\exists\;A>1\;\exists\;\epsilon_2\in(0,1)\;\exists\;H_2\ge 1\;\forall\;t\ge 0:\;\;\;\tau(A^{\delta}t)\ge A^{1+\epsilon_2}\tau(t)-H_2.$$
By Lemma \ref{bigKremark} we can assume w.l.o.g. that $K=A$ and write from now on $K$ for this constant. Moreover, set $H:=\max\{H_1,H_2\}$ and $\epsilon:=\min\{\epsilon_1,\epsilon_2\}$ and estimate as follows for all $t>0$:
\begin{align*}
\sigma\widehat{\star}\tau(K^{\gamma-\delta}t)&=\sup_{s\ge 0}\{\sigma(s)-\tau(s/(K^{\gamma-\delta}t))\}=\sup_{s\ge 0}\{\sigma(s)-\tau(K^{\delta}s/(K^{\gamma}t))\}
\\&
\underbrace{=}_{u=K^{-\gamma}s}\sup_{u\ge 0}\{\sigma(K^{\gamma}u)-\tau(K^{\delta}u/t)\}\le\sup_{u\ge 0}\{K^{1-\epsilon}\sigma(u)-\tau(K^{\delta}u/t)\}+H
\\&
\le\sup_{u\ge 0}\{K^{1-\epsilon}\sigma(u)-K^{1+\epsilon}\tau(u/t)\}+2H\le K^{1-\epsilon}\sup_{u\ge 0}\{\sigma(u)-\tau(u/t)\}+2H
\\&
=K^{1-\epsilon}\sigma\widehat{\star}\tau(t)+2H.
\end{align*}
Since $\lim_{t\rightarrow+\infty}\sigma\widehat{\star}\tau(t)=+\infty$ (see $(a)$ in Lemma \ref{widehatproplemma}) we have verified $(P_{\sigma\widehat{\star}\tau,\gamma-\delta})$, hence $\gamma-\delta<\gamma(\sigma\widehat{\star}\tau)$. And since $\gamma<\gamma(\sigma)$ and $\delta>\overline{\gamma}(\tau)$ are arbitrary the relation $\gamma(\sigma)\le\gamma(\sigma\widehat{\star}\tau)+\overline{\gamma}(\tau)$ follows.\vspace{6pt}

$(ii)$ Analogously, for any $0<\delta<\gamma$ satisfying $\gamma>\overline{\gamma}(\sigma)$ and $\delta<\gamma(\tau)$ we have for all $t>0$:
\begin{align*}
\sigma\widehat{\star}\tau(K^{\gamma-\delta}t)&=\sup_{u\ge 0}\{\sigma(K^{\gamma}u)-\tau(K^{\delta}u/t)\}\ge\sup_{u\ge 0}\{K^{1+\epsilon}\sigma(u)-\tau(K^{\delta}u/t)\}-H
\\&
\ge\sup_{u\ge 0}\{K^{1+\epsilon}\sigma(u)-K^{1-\epsilon}\tau(u/t)\}-2H\ge K^{1+\epsilon}\sup_{u\ge 0}\{\sigma(u)-\tau(u/t)\}-2H
\\&
=K^{1+\epsilon}\sigma\widehat{\star}\tau(t)-2H.
\end{align*}
Since $\lim_{t\rightarrow+\infty}\sigma\widehat{\star}\tau(t)=+\infty$ we have verified $(\overline{P}_{\sigma\widehat{\star}\tau,\gamma-\delta})$, hence $\gamma-\delta>\overline{\gamma}(\sigma\widehat{\star}\tau)$ for any $\gamma>\overline{\gamma}(\sigma)$ and $\delta<\gamma(\tau)$. Thus $\overline{\gamma}(\sigma)\ge\overline{\gamma}(\sigma\widehat{\star}\tau)+\gamma(\tau)$ holds as desired.
\qed\enddemo

Again, let us apply Theorem \ref{uppertransformindexthm} to some special cases.

\begin{corollary}\label{uppertransformindexcor}
Let $\sigma$ be a weight function.
\begin{itemize}
\item[$(i)$] Let $\alpha>0$ and assume that $\sigma(s)=o(s^{1/\alpha})$ as $s\rightarrow+\infty$. When $\gamma(\sigma)>0$, then $$\gamma(\sigma)\le\gamma((((\sigma^{\alpha})^{\star})^{\iota})^{1/\alpha})+\alpha$$
    and $\alpha\le\overline{\gamma}(\sigma)<+\infty$ implies
$$\overline{\gamma}((((\sigma^{\alpha})^{\star})^{\iota})^{1/\alpha})+\alpha\le\overline{\gamma}(\sigma).$$

\item[$(ii)$] Assume that $\tau$ is a weight function satisfying $\tau(0)=0$ and such that $\sigma\widehat{\star}\tau$ is well-defined. If $0<\gamma(\sigma)$, $0<\gamma(\tau)=\overline{\gamma}(\tau)\le\overline{\gamma}(\sigma)<+\infty$, then
$$\gamma(\sigma)\le\gamma(\sigma\widehat{\star}\tau)+\overline{\gamma}(\tau)\le\overline{\gamma}(\sigma\widehat{\star}\tau)+\gamma(\tau)\le\overline{\gamma}(\sigma).$$
\end{itemize}
\end{corollary}

\demo{Proof}
$(i)$ We use Remark \ref{Gevreyweightrem} and Theorem \ref{uppertransformindexthm} applied to $\tau:=\id^{1/\alpha}$, $\alpha>0$ arbitrary, and take into account Corollary \ref{uppertransformfinitecor}. Consequently, $\gamma(\sigma)\le\gamma(\sigma\widehat{\star}\id^{1/\alpha})+\alpha$ and $\overline{\gamma}(\sigma\widehat{\star}\id^{1/\alpha})+\alpha\le\overline{\gamma}(\sigma)$,
thus the assertion follows by recalling \eqref{upperLegendregeneral}.\vspace{6pt}

$(ii)$ This is an immediate consequence of Theorem \ref{uppertransformindexthm} and the assumptions.
\qed\enddemo

In particular, one infers:
\begin{itemize}
\item[$(*)$] When $\sigma$ satisfies $\sigma(s)=o(s^{1/\alpha})$ and $\gamma(\sigma)=\overline{\gamma}(\sigma)\ge\alpha>0$ then
\begin{equation}\label{uppertransformindexcorequ1}
\gamma(\sigma)=\gamma((((\sigma^{\alpha})^{\star})^{\iota})^{1/\alpha})+\alpha=\overline{\gamma}((((\sigma^{\alpha})^{\star})^{\iota})^{1/\alpha})+\alpha=\overline{\gamma}(\sigma);
\end{equation}

\item[$(*)$] more specially,
\begin{equation}\label{uppertransformindexcorequ}
\forall\;\beta>\alpha>0:\;\;\;\beta=\gamma(\id^{1/\beta}\widehat{\star}\id^{1/\alpha})+\alpha=\overline{\gamma}(\id^{1/\beta}\widehat{\star}\id^{1/\alpha})+\alpha.
\end{equation}
\end{itemize}

\begin{remark}\label{indexnaturalupperrem}
\emph{Analogously to Remark \ref{indexnaturallowerrem} we comment on the extreme cases for the growth indices:}
\begin{itemize}
\item[$(*)$] \emph{If in $(i)$ in Theorem \ref{uppertransformindexthm} one has $\overline{\gamma}(\tau)<\gamma(\sigma)=+\infty$, then $\gamma(\sigma\widehat{\star}\tau)=+\infty$ as well. Since by assumption $\sigma\widehat{\star}\tau$ is a weight function and so $\gamma(\sigma\widehat{\star}\tau)\ge 0$, \eqref{uppertransformindexthmequ} becomes trivial if $\gamma(\sigma)\le\overline{\gamma}(\tau)$ and so formally, in particular, when $\overline{\gamma}(\tau)=+\infty$. However, in this case the precise value of $\gamma(\sigma\widehat{\star}\tau)$ becomes unclear.}

\item[$(*)$] \emph{\eqref{uppertransformindexthmequ1} is formally trivial when $\overline{\gamma}(\sigma)=+\infty$ but again, in this case, the value of $\overline{\gamma}(\sigma\widehat{\star}\tau)$ is unclear. Moreover, the given proof is also valid if $0=\gamma(\tau)\le\overline{\gamma}(\sigma)<+\infty$.}

\item[$(*)$] \emph{The relation $\gamma(\tau)\le\overline{\gamma}(\sigma)$ in \eqref{uppertransformindexthmequ1} is natural for having $\overline{\gamma}(\sigma\widehat{\star}\tau)+\gamma(\tau)\le\overline{\gamma}(\sigma)$: Therefore, note and recall that $\overline{\gamma}(\omega)\ge 0$ for any weight function $\omega$.}

\item[$(*)$] \emph{In $(i)$ in Corollary \ref{uppertransformindexcor} the first estimate is still valid (but trivial) if $\gamma(\sigma)=0$ and the second estimate is formally trivial when $\overline{\gamma}(\sigma)=+\infty$.}

\item[$(*)$] \emph{Finally, also $(ii)$ in Corollary \ref{uppertransformindexcor} is valid for the extreme cases as long as $\gamma(\tau)=\overline{\gamma}(\tau)\le\overline{\gamma}(\sigma)$ is (formally) still valid.}
\end{itemize}
\end{remark}

For the generalized upper Legendre conjugate the previous remark becomes relevant since $\overline{\gamma}(\sigma)=+\infty$ appears naturally:

\begin{remark}\label{slowlyvaryingom6problemrem}
\emph{As stated in $(ii)$ in Remark \ref{slowlyvarrem} the assumption that $\sigma\widehat{\star}\sigma$ is well-defined implies the fact that $\sigma$ has to be slowly varying. Therefore, $\gamma(\sigma)=+\infty$ has to be valid, see \cite[Sect. 2.2 \& 2.3, Thm. 2.11]{index} and then, by \cite[Thm. 2.16, Cor. 2.17, Rem. 2.18]{index}, also $\overline{\gamma}(\sigma)=+\infty$. Hence the formal assumptions on the indices above fail, the precise values for $\gamma(\sigma\widehat{\star}\sigma)$ and $\overline{\gamma}(\sigma\widehat{\star}\sigma)$ are becoming unclear, but, however, the desired estimates are (formally) trivial.}

\emph{This comment should be compared with \eqref{uppertransformindexcorequ} where we naturally require the strict inequality $\beta>\alpha$. And for the weight from Examples \ref{logexamplelower} and \ref{logexampleupper} we have $\gamma(\log_{+})=+\infty=\overline{\gamma}(\log_{+})$ and hence all estimates in Theorems \ref{lowertransformindexthm} and \ref{uppertransformindexthm} are reduced to the formal equality $+\infty=+\infty$.}
\end{remark}

\subsection{On the relations between weight functions}
In general, the analogue of Lemma \ref{lowertransformrelationlemma} is not clear for the generalized upper conjugate; therefore note that the property that $\sigma\widehat{\star}\tau$ is well defined and also \eqref{conditionB} are not automatically preserved under equivalences of weight functions.

On the other hand, when also involving additional growth conditions on the weights, we have the following comments:

\begin{itemize}
\item[$(i)$] First, the indices $\gamma(\cdot)$ and $\overline{\gamma}(\cdot)$ from Section \ref{growthindexsect} and the relation $\sigma(s)=o(s^{1/\alpha})$ are preserved w.r.t. \hyperlink{sim}{$\sim$}. Hence $(i)$ in Corollary \ref{uppertransformindexcor} remains valid if $\sigma$ is replaced by some/any equivalent weight $\sigma_1$.

\item[$(ii)$] In the weight sequence setting we infer in Section \ref{uppertrafosect} that in \eqref{uppertransformindexcorequ} one can replace $\id^{1/\beta}$ by $\omega_{\mathbf{G}^{\beta}}$ and $\id^{1/\alpha}$ by $\omega_{\mathbf{G}^{\alpha}}$; see Example \ref{uppertrafoGevreyex}.

    This comment is also consistent with the following observations which are based on Remark \ref{om1om6indexrem}.

\item[$(iii)$] Let $\sigma$, $\tau$ be weight functions, a-priori not necessarily satisfying $\tau(0)=0$ and/or such that $\sigma\widehat{\star}\tau$ is well defined. Let $\sigma_1$ be such that $\sigma\hyperlink{ompreceq}{\preceq}\sigma_1$, i.e. $\sigma_1(t)=O(\sigma(t))$ as $t\rightarrow+\infty$ and so $\sigma_1(t)\le C_1\sigma(t)+C_1$ for some $C_1\ge 1$ and all $t\ge 0$. Now let us set $\tau_1:=\frac{1}{C_1}\tau$ and estimate as follows for all $t>0$:
\begin{align*}
\sigma_1\widehat{\star}\tau(t)&=\sup_{s\ge 0}\{\sigma_1(s)-\tau(s/t)\}\le\sup_{s\ge 0}\{C_1\sigma(s)-\tau(s/t)\}+C_1
\\&
=\sup_{s\ge 0}\{C_1\sigma(s)-C_1\tau_1(s/t)\}+C_1=C_1\sigma\widehat{\star}\tau_1(t)+C_1.
\end{align*}
This computation also implies that $\sigma_1\widehat{\star}\tau(t)<+\infty$ provided that $\sigma\widehat{\star}\tau_1(t)<+\infty$ hence, in particular, if $\sigma\widehat{\star}\tau_1$ is well defined then $\sigma_1\widehat{\star}\tau$, too. By enlarging the constant $C_1$ if necessary we also get this estimate for $t=0$; see $(b)$ in Lemma \ref{widehatproplemma}.

\item[$(iv)$] If $\tau$ satisfies \hyperlink{om6}{$(\omega_6)$} resp. $\overline{\gamma}(\tau)<+\infty$, then by iterating this condition $n$ times with $n\in\NN_{>0}$ minimal such that $C_1\le 2^n$, we have that $C_1\tau(s/(Ht))\le\tau(s/t)+H$ holds for some $H\ge 1$ which is depending on given $C_1$ via \hyperlink{om6}{$(\omega_6)$} and for all $s\ge 0$, $t>0$. Following the above computations we get $\sigma_1\widehat{\star}\tau(t)\le C_1\sigma\widehat{\star}\tau(Ht)+H+C_1$ for all $t>0$. In case $\sigma\widehat{\star}\tau$ is well defined, then $\sigma_1\widehat{\star}\tau$ too and in order to have $\sigma\widehat{\star}\tau\hyperlink{ompreceq}{\preceq}\sigma_1\widehat{\star}\tau$ one requires \hyperlink{om1}{$(\omega_1)$} either for $\sigma\widehat{\star}\tau$ or for $\sigma_1\widehat{\star}\tau$  similarly like in the proof of $(ii)\Rightarrow(iii)$ in Lemma \ref{om1om6lemma}; see also Remark \ref{om1om6lemmaremnew}.

\item[$(v)$] Analogously, let $\sigma$, $\tau$, $\tau_1$ be given weight functions such that $\tau_1\hyperlink{ompreceq}{\preceq}\tau$ and assume that either $\tau$ or $\tau_1$ satisfies \hyperlink{om6}{$(\omega_6)$}. So $\tau(t)\le C_1\tau_1(t)+C_1$ holds for some $C_1\ge 1$ and all $t\ge 0$. If $\tau$ has \hyperlink{om6}{$(\omega_6)$}, then by iteration we get $\tau(t)\le C_1\tau_1(t)+C_1\le\tau_1(Ht)+H+C_1$ for some $H\ge 1$ and all $t\ge 0$. If $\tau_1$ has \hyperlink{om6}{$(\omega_6)$}, then similarly $\tau(t)\le\frac{1}{C_1}\tau(Ht)+\frac{H}{C_1}\le\tau_1(Ht)+1+\frac{H}{C_1}\le\tau_1(Ht)+C_1+H$. Hence, in any case, for all $t>0$:
\begin{align*}
\sigma\widehat{\star}\tau(t)&=\sup_{s\ge 0}\{\sigma(s)-\tau(s/t)\}\ge\sup_{s\ge 0}\{\sigma(s)-\tau_1((Hs)/t)\}-C_1-H
\\&
=\sigma\widehat{\star}\tau_1(t/H)-C_1-H.
\end{align*}
Consequently, $\sigma\widehat{\star}\tau_1(t)<+\infty$ provided that $\sigma\widehat{\star}\tau(Ht)<+\infty$ and so, in particular, if $\sigma\widehat{\star}\tau$ is well defined then $\sigma\widehat{\star}\tau_1$, too.

Finally, when $\sigma\widehat{\star}\tau$ is well defined, then in order to have $\sigma\widehat{\star}\tau\hyperlink{ompreceq}{\preceq}\sigma\widehat{\star}\tau_1$ one requires \hyperlink{om1}{$(\omega_1)$} either for $\sigma\widehat{\star}\tau$ or for $\sigma\widehat{\star}\tau_1$.
\end{itemize}

We summarize this information in the next statement.

\begin{theorem}\label{upperconjequivcor}
Let $\sigma$, $\sigma_1$, $\tau$, $\tau_1$ be weight functions with $\tau(0)=0=\tau_1(0)$ and assume that $\sigma\widehat{\star}\tau$ is well defined.

\begin{itemize}
\item[$(i)$] Let $\sigma\hyperlink{ompreceq}{\preceq}\sigma_1$ be valid and $\overline{\gamma}(\tau)<+\infty$; i.e. \hyperlink{om6}{$(\omega_6)$} for $\tau$.
\begin{itemize}
\item[$(*)$] Then $\sigma_1\widehat{\star}\tau$ is a weight function as well.

\item[$(*)$] When either $\gamma(\sigma\widehat{\star}\tau)>0$ or $\gamma(\sigma_1\widehat{\star}\tau)>0$ is valid, then $\sigma\widehat{\star}\tau\hyperlink{ompreceq}{\preceq}\sigma_1\widehat{\star}\tau$ holds.

\item[$(*)$] If in addition even $\sigma_1\hyperlink{sim}{\sim}\sigma$, then $\sigma\widehat{\star}\tau\hyperlink{sim}{\sim}\sigma_1\widehat{\star}\tau$ and \eqref{uppertransformindexthmequ} holds for $\sigma_1\widehat{\star}\tau$.
\end{itemize}

\item[$(ii)$] Let $\tau_1\hyperlink{ompreceq}{\preceq}\tau$ be valid and let either $\overline{\gamma}(\tau)<+\infty$ or $\overline{\gamma}(\tau_1)<+\infty$.
\begin{itemize}
\item[$(*)$] Then $\sigma\widehat{\star}\tau_1$ is a weight function as well.

\item[$(*)$] When either $\gamma(\sigma\widehat{\star}\tau)>0$ or $\gamma(\sigma\widehat{\star}\tau_1)>0$, then $\sigma\widehat{\star}\tau\hyperlink{ompreceq}{\preceq}\sigma\widehat{\star}\tau_1$.

\item[$(*)$] If in addition even $\tau_1\hyperlink{sim}{\sim}\tau$, then $\sigma\widehat{\star}\tau\hyperlink{sim}{\sim}\sigma\widehat{\star}\tau_1$ and \eqref{uppertransformindexthmequ} holds for $\sigma\widehat{\star}\tau_1$.
\end{itemize}
\end{itemize}
Both $(i)$ and $(ii)$ is ensured when, in addition to the particular growth relations between the given weights, we assume
\begin{equation}\label{upperconjequivcorequ}
0\le\overline{\gamma}(\tau)<\gamma(\sigma)\le+\infty.
\end{equation}
\end{theorem}

The formulation of this result should be compared with Remark \ref{om1om6lemmaremnew} and this Theorem applies to Corollary \ref{uppertransformindexcor}. Note that the second strict inequality in \eqref{upperconjequivcorequ} has to fail for the case $\tau=\sigma$ and for $\sigma\widehat{\star}\sigma$ being well defined necessarily $\gamma(\sigma)=+\infty=\overline{\gamma}(\sigma)$; recall Remark \ref{slowlyvaryingom6problemrem}.

\demo{Proof}
$(i)$ By the assumptions $\sigma\widehat{\star}\tau$ is a weight function. Then apply the previous comments $(iii)$ and $(iv)$. Necessarily, $\sigma_1\widehat{\star}\tau(t)<+\infty$ for all $t>0$ and so $\sigma_1\widehat{\star}\tau$ is a weight function, too.

When $\sigma$ and $\sigma_1$ are even equivalent, then the already shown implications hold. Therefore, repeating the arguments in $(iii)$ and $(iv)$ but starting with the weight function $\sigma_1\widehat{\star}\tau$, gives $\sigma_1\widehat{\star}\tau\hyperlink{ompreceq}{\preceq}\sigma\widehat{\star}\tau$ and hence the desired equivalence.\vspace{6pt}

Concerning the supplement, if $0\le\overline{\gamma}(\tau)<\gamma(\sigma)\le+\infty$, then $(i)$ in Theorem \ref{uppertransformindexthm} yields $\gamma(\sigma\widehat{\star}\tau)\ge\gamma(\sigma)-\overline{\gamma}(\tau)>0$ and $\gamma(\sigma\widehat{\star}\tau)=+\infty$ provided that $\gamma(\sigma)=+\infty$, see \eqref{uppertransformindexthmequ}. Thus \hyperlink{om1}{$(\omega_1)$} follows for $\sigma\widehat{\star}\tau$ and hence the conclusions (and in case of equivalence also $\gamma(\sigma_1)=\gamma(\sigma)$).\vspace{6pt}

$(ii)$ Again, $\sigma\widehat{\star}\tau$ is a weight function and recall that $\overline{\gamma}(\tau)<+\infty$ resp. $\overline{\gamma}(\tau_1)<+\infty$ is precisely \hyperlink{om6}{$(\omega_6)$} for $\tau$ resp. for $\tau_1$. Thus we apply comment $(v)$ and hence the conclusion follows again.

In case $\tau_1\hyperlink{sim}{\sim}\tau$, then by the above $\sigma\widehat{\star}\tau_1$ is a weight function and we apply $(v)$ to $\sigma\widehat{\star}\tau_1$.\vspace{6pt}

Concerning the supplement note that again $(i)$ in Theorem \ref{uppertransformindexthm} ensures  $\gamma(\sigma\widehat{\star}\tau)>0$ (and $\gamma(\sigma\widehat{\star}\tau)=+\infty$ provided that $\gamma(\sigma)=+\infty$ holds). In case of equivalence even $\overline{\gamma}(\tau_1)=\overline{\gamma}(\tau)<+\infty$ is true.
\qed\enddemo

Finally, we comment on the stronger relation \hyperlink{omvartriangle}{$\vartriangleleft$}:
Analogously to $(iii)$ above, let $\sigma$, $\tau$ be weight functions and consider $\sigma_1$ such that $\sigma\hyperlink{omvartriangle}{\vartriangleleft}\sigma_1$. Thus for all $0<c<1$ there exists $D_c\ge 1$ such that $\sigma_1(t)\le c\sigma(t)+D_c$ for all $t\ge 0$ and so
\begin{align*}
\sigma_1\widehat{\star}\tau(t)&=\sup_{s\ge 0}\{\sigma_1(s)-\tau(s/t)\}\le\sup_{s\ge 0}\{c\sigma(s)-\tau(s/t)\}+D_c
\\&
\le c\sup_{s\ge 0}\{\sigma(s)-\tau(s/t)\}+D_c=c\sigma\widehat{\star}\tau(t)+D_c.
\end{align*}
Analogously to $(v)$ above let $\sigma$, $\tau$, $\tau_1$ be given such that $\tau_1\hyperlink{omvartriangle}{\vartriangleleft}\tau$, then for all $0<c<1$ there exists $D_c\ge 1$ such that $\tau(t)\le c\tau_1(t)+D_c$ for all $t\ge 0$ and so
\begin{align*}
\sigma\widehat{\star}\tau(t)&=\sup_{s\ge 0}\{\sigma(s)-\tau(s/t)\}\ge\sup_{s\ge 0}\{\sigma(s)-c{\tau_1}(s/t)\}-D_c
\\&
\ge c\sup_{s\ge 0}\{\sigma(s)-\tau_1(s/t)\}-D_c=c\sigma\widehat{\star}\tau_1(t)-D_c.
\end{align*}

These computations directly yield the next main statement.

\begin{theorem}\label{upperconjequivcor1}
Let $\sigma$, $\sigma_1$, $\tau$, $\tau_1$ be weight functions with $\tau(0)=0=\tau_1(0)$ and assume that $\sigma\widehat{\star}\tau$ is well defined.

\begin{itemize}
\item[$(i)$] Let $\sigma\hyperlink{omvartriangle}{\vartriangleleft}\sigma_1$ be valid, then $\sigma\widehat{\star}\tau\hyperlink{omvartriangle}{\vartriangleleft}\sigma_1\widehat{\star}\tau$ holds and $\sigma_1\widehat{\star}\tau$ is a weight function as well.

\item[$(ii)$] Let $\tau_1\hyperlink{omvartriangle}{\vartriangleleft}\tau$ be valid, then $\sigma\widehat{\star}\tau\hyperlink{ompreceq}{\preceq}\sigma\widehat{\star}\tau_1$ and $\sigma\widehat{\star}\tau_1$ is a weight function as well.
\end{itemize}
\end{theorem}

\section{Generalized Legendre conjugates of associated weight functions}\label{trafosect}
The aim is to study in detail the above operations within the weight sequence setting; i.e. when $\sigma=\omega_{\mathbf{M}}$ and $\tau=\omega_{\mathbf{N}}$ for weight sequences $\mathbf{M},\mathbf{N}\in\RR_{>0}^{\NN}$. In order to deal with weight functions (in the sense of Definition \ref{weightfctdef}), it is natural to focus on (log-convex) sequences satisfying the standard requirement $\mathbf{M}_{\iota}=+\infty=\mathbf{N}_{\iota}$. Concerning non-standard cases we refer to Appendix \ref{nonstandardappendix}.\vspace{6pt}

First, note that the weight function $\log_{+}$ from Examples \ref{logexamplelower} and \ref{logexampleupper} corresponds to the ``exotic situation'' when one allows resp. considers ``sequences'' $\mathbf{M}$ such that $M_p=+\infty$ for all sufficiently large $p$; see \cite[Lemma 7.2]{solidassociatedweight}. Indeed, this Lemma completely illustrates this particular and limiting example.

\subsection{Generalized lower Legendre conjugate}\label{lowertrafosect}
Let $\mathbf{M}=(M_p)_{p\in\NN},\mathbf{N}=(N_p)_{p\in\NN}\in\RR_{>0}^{\NN}$ be given such that $\mathbf{M}_{\iota},\mathbf{N}_{\iota}>0$ and let $\omega_{\mathbf{M}}$, $\omega_{\mathbf{N}}$ denote the corresponding associated weight functions. Set
\begin{equation}\label{Itset}
\mathcal{I}_t:=(t\mathbf{N}_{\iota}^{-1},\mathbf{M}_{\iota}),\;\;\;t\in[0,\mathbf{M}_{\iota}\cdot\mathbf{N}_{\iota}),
\end{equation}
with the conventions $\mathbf{N}_{\iota}^{-1}=0$ if $\mathbf{N}_{\iota}=+\infty$, and $\mathbf{M}_{\iota}\cdot\mathbf{N}_{\iota}=+\infty$ if either $\mathbf{M}_{\iota}=+\infty$ or $\mathbf{N}_{\iota}=+\infty$. Then, analogously as in \eqref{lowertransformgen} let us introduce
\begin{equation}\label{lowertransform}
\omega_{\mathbf{M}}\check{\star}\omega_{\mathbf{N}}(t):=\inf_{s\in\mathcal{I}_t}\{\omega_{\mathbf{M}}(s)+\omega_{\mathbf{N}}(t/s)\},\;\;\;t\in[0,\mathbf{M}_{\iota}\cdot\mathbf{N}_{\iota}).
\end{equation}
In case $\mathbf{M}_{\iota}=+\infty=\mathbf{N}_{\iota}$ with the above conventions $\mathcal{I}_t=(0,+\infty)$ and \eqref{lowertransform} precisely corresponds to \eqref{lowertransformgen}. For this expression we refer also to \cite[$(3.6)$]{Boman98}. Indeed, in \cite{Boman98} \emph{log-convex} sequences have been considered and so $\mathbf{M}_{\iota}>0$, $\mathbf{N}_{\iota}>0$ follows but, however, $\mathbf{M}_{\iota}=+\infty$ resp. $\mathbf{N}_{\iota}=+\infty$ is not clear in general; see $(b)$ in Section \ref{preliminarysection}. Also in \cite{Boman98} in \eqref{lowertransform} all $t>0$ and $\inf_{s>0}$ have been considered but for non-standard cases $\mathcal{I}_t$ is the canonical natural set of definition; we refer to Section \ref{nonstandardlowersection} for more details.\vspace{6pt}

In the main result we establish now the connection between $\omega_{\mathbf{M}}\check{\star}\omega_{\mathbf{N}}$ and $\omega_{\mathbf{M}\cdot\mathbf{N}}$; the proof basically repeats the arguments given in \cite[Lemma 4]{Boman98}.

\begin{theorem}\label{propBomanLemma4}
Let $\mathbf{M},\mathbf{N}\in\RR_{>0}^{\NN}$ be log-convex with $\mathbf{M}_{\iota}=+\infty=\mathbf{N}_{\iota}$. Then
\begin{equation}\label{propBomanLemma4equ}
\forall\;t\in[0,+\infty):\;\;\;\omega_{\mathbf{M}\cdot\mathbf{N}}(t)=\omega_{\mathbf{M}}\check{\star}\omega_{\mathbf{N}}(t).
\end{equation}
\end{theorem}

Recall that in \cite[Lemma 4]{Boman98} by the assumptions on the sequences one only has $\mathbf{M}_{\iota}>0$ and $\mathbf{N}_{\iota}>0$. If now $\mathbf{M}_{\iota},\mathbf{N}_{\iota}<+\infty$, then the part where it is claimed \emph{``Fix an arbitrary $t>0$ and choose $k$ such that $t\in[\mu_k\nu_k,\mu_{k+1}\nu_{k+1}]$''} is not satisfied for $t>\mathbf{M}_{\iota}\mathbf{N}_{\iota}$; see $(b)$ in Section \ref{preliminarysection}. In Section \ref{nonstandardlowersection} we study this non-standard setting in more detail. Finally, note that Theorem \ref{propBomanLemma4} implies as a consequence the commutativity of $\check{\star}$ in the weight sequence setting.

\demo{Proof}
First, since both $\mathbf{M}$ and $\mathbf{N}$ are log-convex we have that $\mathbf{M}\cdot\mathbf{N}$ is log-convex. Then $\lim_{p\rightarrow+\infty}\mu_p=\mathbf{M}_{\iota}=+\infty=\mathbf{N}_{\iota}=\lim_{p\rightarrow+\infty}\nu_p$ and, moreover,
\begin{equation}\label{propBomanLemma4equ1}
\lim_{p\rightarrow+\infty}\mu_p\nu_p=(\mathbf{M}\cdot\mathbf{N})_{\iota}=\mathbf{M}_{\iota}\cdot\mathbf{N}_{\iota}=+\infty.
\end{equation}
For any given $t\ge 0$, $p\in\NN$ and $s>0$ we write
\begin{equation}\label{stdecompose}
\frac{t^pM_0\cdot N_0}{M_p\cdot N_p}=\frac{s^pM_0}{M_p}\frac{(t/s)^pN_0}{N_p}.
\end{equation}
For fixed $t\ge 0$ take in this expression the supremum over all $p\in\NN$ and the infimum over all $s>0$, hence by definition $\exp(\omega_{\mathbf{M}\cdot\mathbf{N}}(t))\le\exp(\omega_{\mathbf{M}}\check{\star}\omega_{\mathbf{N}}(t))$ for all $t\ge 0$.

For the converse, first consider $t\ge 0$ such that $\mu_p\nu_p\le t\le\mu_{p+1}\nu_{p+1}$ for some $p\in\NN_{>0}$. Then choose $\theta\in[0,1]$ such that $t=(\mu_p\nu_p)^{1-\theta}(\mu_{p+1}\nu_{p+1})^{\theta}$ and set $s_t:=\mu_p^{1-\theta}\mu_{p+1}^{\theta}(>0)$. Hence $\mu_p\le s_t\le\mu_{p+1}$, $\frac{t}{s_t}=\nu_p^{1-\theta}\nu_{p+1}^{\theta}$ and so $\nu_p\le\frac{t}{s_t}\le\nu_{p+1}$. In view of \eqref{lemma1equ} we have $$\omega_{\mathbf{M}\cdot\mathbf{N}}(t)=\log\left(\frac{t^pM_0N_0}{M_pN_p}\right)=\log\left(\frac{s_t^pM_0}{M_p}\right)+\log\left(\frac{(t/s_t)^pN_0}{N_p}\right)=\omega_{\mathbf{M}}(s_t)+\omega_{\mathbf{N}}(t/s_t),$$ hence $\omega_{\mathbf{M}\cdot\mathbf{N}}(t)\ge\omega_{\mathbf{M}}\check{\star}\omega_{\mathbf{N}}(t)$ is verified for all $t\ge\mu_1\nu_1$.

Second, for $0\le t<\mu_1\nu_1$ there exists $0<s_t<\mu_1$ such that $\frac{t}{s_t}<\nu_1$. Consequently, again in view of \eqref{lemma1equ}, we obtain $\omega_{\mathbf{M}\cdot\mathbf{N}}(t)=0=\omega_{\mathbf{M}}(s_t)=\omega_{\mathbf{N}}(t/s_t)$ and so the desired estimate is established for all $t\ge 0$.
\qed\enddemo

\subsection{Generalized upper Legendre conjugate}\label{uppertrafosect}
Let $\mathbf{M},\mathbf{N}\in\RR_{>0}^{\NN}$ be given with associated weight functions $\omega_{\mathbf{M}}$, $\omega_{\mathbf{N}}$. In order to introduce $\omega_{\mathbf{M}}\widehat{\star}\omega_{\mathbf{N}}$, first note that assumption $(A)$ in Section \ref{arbuppertrafosect} is clear because $\omega_{\mathbf{N}}(0)=0$.

But concerning the fact that $\omega_{\mathbf{M}}\widehat{\star}\omega_{\mathbf{N}}$ is well defined recall Lemma \ref{uppertransformfinite}; we show that $\mathbf{M},\mathbf{N}$ have to satisfy a particular growth relation and formulate the next technical result in a very general framework since it purely deals with abstractly given weight sequences and might be of interest in a wider context as well.

\begin{proposition}\label{sequencenotslowly}
Let $\mathbf{M},\mathbf{N}\in\RR_{>0}^{\NN}$ be given and assume that $\mathbf{N}$ is log-convex. Then $\lim_{p\rightarrow+\infty}(N_p)^{1/p}=\mathbf{N}_{\iota}\in(0,+\infty]$ and we get:
\begin{itemize}
\item[$(I)$] The following are equivalent:
\begin{itemize}
\item[$(i)$] The associated weight functions satisfy
\begin{equation}\label{sequencenotslowlyequ}
\forall\;t\ge 0\;\exists\;D_t\ge 1\;\forall\;u\in[0,\mathbf{N}_{\iota}):\;\;\;\omega_{\mathbf{M}}(tu)\le\omega_{\mathbf{N}}(u)+D_t.
\end{equation}

\item[$(ii)$] $\mathbf{N}\hyperlink{mtriangle}{\vartriangleleft}\mathbf{M}$ holds.
\end{itemize}

\item[$(II)$]
\begin{itemize}
\item[$(i)$] The associated weight functions satisfy
\begin{equation}\label{sequencenotslowlyequvar}
\exists\;t_0\in(0,+\infty)\;\forall\;t\in[0,t_0)\;\exists\;D_t\ge 1\;\forall\;u\in[0,\mathbf{N}_{\iota}):\;\;\;\omega_{\mathbf{M}}(tu)\le\omega_{\mathbf{N}}(u)+D_t.
\end{equation}

\item[$(ii)$] $\mathbf{N}\hyperlink{preceq}{\preceq}\mathbf{M}$ holds.
\end{itemize}
\end{itemize}
\end{proposition}

\emph{Note:} Any of the equivalent assertions in $(I)$ implies $\mathbf{M}_{\iota}=+\infty$ whereas in $(II)$ this is only clear when $\mathbf{N}_{\iota}=+\infty$. If $\mathbf{N}_{\iota}\in(0,+\infty)$, then formally in \eqref{sequencenotslowlyequ} one can even take any $u\ge 0$ by setting $\omega_{\mathbf{N}}(u)=+\infty$ for all $u>\mathbf{N}_{\iota}$. If $\mathbf{N}_{\iota}=+\infty$, then \eqref{sequencenotslowlyequ} is precisely \eqref{uppetransfromwelldefequ} with $t_0=+\infty$ between $\omega_{\mathbf{M}}$ and $\omega_{\mathbf{N}}$; see $(B)$ in Section \ref{arbuppertrafosect}. Hence, for log-convex sequences $\mathbf{M}$, $\mathbf{N}$ such that $\mathbf{N}_{\iota}=+\infty$ the function $\omega_{\mathbf{M}}\widehat{\star}\omega_{\mathbf{N}}$ is well defined if and only if $\mathbf{N}\hyperlink{mtriangle}{\vartriangleleft}\mathbf{M}$.

Moreover, in $(II)$ the equality $t_0=(\underline{C}_{\mathbf{N}\preceq\mathbf{M}})^{-1}$ holds and therefore, via \eqref{MNiotarelation}, we (formally) get $\mathbf{M}_{\iota}\ge\mathbf{N}_{\iota}(\underline{C}_{\mathbf{N}\preceq\mathbf{M}})^{-1}=\mathbf{N}_{\iota}t_0$. Here, $\underline{C}_{\mathbf{N}\preceq\mathbf{M}}\in(0,+\infty)$ and so it follows that $t_0\in(0,+\infty)$ too and \eqref{sequencenotslowlyequvar} is precisely \eqref{uppetransfromwelldefequ} with $t_0$; i.e. $\omega_{\mathbf{M}}\widehat{\star}\omega_{\mathbf{N}}(t)<+\infty$ for all $t\in(0,t_0)$. (In $(I)$ one has the same correspondence by recognizing $t_0=+\infty=(\underline{C}_{\mathbf{N}\preceq\mathbf{M}})^{-1}$ since $\mathbf{N}\hyperlink{mtriangle}{\vartriangleleft}\mathbf{M}$.)

\demo{Proof}
First recall that log-convexity of $\mathbf{N}$ implies $\lim_{p\rightarrow+\infty}(N_p)^{1/p}=\mathbf{N}_{\iota}\in(0,+\infty]$.\vspace{6pt}

Using different notation and when $\mathbf{N}_{\iota}=+\infty=\mathbf{M}_{\iota}$, then $(I)$ has already been shown in \cite[Lemma 4.2]{weightedentireinclusion2}.

$(I)(i)\Rightarrow(ii)$ \eqref{sequencenotslowlyequ} gives, in particular, that $\mathbf{M}_{\iota}=+\infty$: If $\mathbf{M}_{\iota}<+\infty$, then for any fixed $u\in(0,\mathbf{N}_{\iota})$ and all $t>\mathbf{M}_{\iota}u^{-1}$ one has $\omega_{\mathbf{M}}(tu)=+\infty$ and $\omega_{\mathbf{N}}(u)<+\infty$, a contradiction. Then \eqref{reverseformula} implies for all $p\in\NN$ and $t>0$:
\begin{align*}
N_p&=N_0\sup_{u\in(0,\mathbf{N}_{\iota})}\frac{u^p}{\exp(\omega_{\mathbf{N}}(u))}\le N_0e^{D_t}\sup_{u\in(0,\mathbf{N}_{\iota})}\frac{u^p}{\exp(\omega_{\mathbf{M}}(tu))}
\\&
=N_0e^{D_t}\sup_{v\in(0,t\mathbf{N}_{\iota})}\frac{(v/t)^p}{\exp(\omega_{\mathbf{M}}(v))}\le N_0e^{D_t}\sup_{v\in(0,+\infty)}\frac{(v/t)^p}{\exp(\omega_{\mathbf{M}}(v))}
\\&
=\frac{N_0e^{D_t}}{M_0}\left(\frac{1}{t}\right)^pM^{\on{lc}}_p\le\frac{N_0e^{D_t}}{M_0}\left(\frac{1}{t}\right)^pM_p.
\end{align*}
This verifies \eqref{triangleestim} with $h:=\frac{1}{t}$ and $C_h=C_t:=\frac{N_0e^{D_t}}{M_0}$. Then let $t\rightarrow+\infty$ and so $\mathbf{N}\hyperlink{mtriangle}{\vartriangleleft}\mathbf{M}$ is shown.\vspace{6pt}

$(I)(ii)\Rightarrow(i)$ First, $\mathbf{N}_{\iota}\in(0,+\infty]$ and $\mathbf{N}\hyperlink{mtriangle}{\vartriangleleft}\mathbf{M}$ imply $\mathbf{M}_{\iota}=+\infty$. Then, by assumption
$$\forall\;h>0\;\exists\;C_h\ge 1\;\forall\;p\in\NN\;\forall\;u\ge 0:\;\;\;\frac{M_0(u/h)^p}{M_p}\le C_h\frac{M_0}{N_0}\frac{N_0u^p}{N_p},$$
hence
$$\forall\;h>0\;\exists\;C_h\ge 1\;\forall\;u\ge 0:\;\;\;\exp(\omega_{\mathbf{M}}(u/h))\le\frac{C_hM_0}{N_0}\exp(\omega_{\mathbf{N}}(u)),$$
which implies the estimate in \eqref{sequencenotslowlyequ} with $t:=\frac{1}{h}$, $D_t=D_h:=\log\left(\frac{C_hM_0}{N_0}\right)$ and note that naturally we restrict to $u\in[0,\mathbf{N}_{\iota})$. Finally, since $h>0$ is arbitrary estimate \eqref{sequencenotslowlyequ} is verified.\vspace{6pt}

$(II)(i)\Rightarrow(ii)$ First note that for any $t\in[0,t_0)$ we have $t\mathbf{N}_{\iota}\le\mathbf{M}_{\iota}$: If there exists $t_1\in(0,t_0)$ such that $t_1\mathbf{N}_{\iota}>\mathbf{M}_{\iota}$ (which is only possible if $\mathbf{M}_{\iota}<+\infty$), then for all $t$ sufficiently close to $t_0$ and $u$ sufficiently close to $\mathbf{N}_{\iota}$ we get $tu>\mathbf{M}_{\iota}$ hence contradicting \eqref{sequencenotslowlyequvar} since the left-hand side is equal to infinity for all such values whereas the right-hand side is finite. Then let $t\in[0,t_0)$, follow the estimate in $(I)(i)\Rightarrow(ii)$ and with $+\infty$ being replaced by $\mathbf{M}_{\iota}$ to obtain $\mathbf{N}\hyperlink{preceq}{\preceq}\mathbf{M}$; more precisely \eqref{preceqalternative1} with $h:=\frac{1}{t}$ and $C_h=C_t:=\frac{N_0e^{D_t}}{M_0}$. Then let $t\rightarrow t_0$ and so $\underline{C}_{\mathbf{N}\preceq\mathbf{M}}\le t_0^{-1}$. \vspace{6pt}

$(II)(ii)\Rightarrow(i)$ By \eqref{preceqalternative1} we have
$$\forall\;h>\underline{C}_{\mathbf{N}\preceq\mathbf{M}}\;\exists\;C_h>0\;\forall\;p\in\NN\;\forall\;u\ge 0:\;\;\;\frac{M_0(u/h)^p}{M_p}\le\frac{M_0}{N_0}C_h\frac{N_0u^p}{N_p}.$$
We can restrict to $u\in[0,\mathbf{N}_{\iota})$ and in this case $\frac{u}{h}\in[0,\mathbf{M}_{\iota})$ since $\mathbf{N}_{\iota}\le \underline{C}_{\mathbf{N}\preceq\mathbf{M}}\mathbf{M}_{\iota}$; see $(i)-(iii)$ in Section \ref{assofunctsect} and \eqref{MNiotarelation}. Then the above estimate gives
$$\forall\;h>\underline{C}_{\mathbf{N}\preceq\mathbf{M}}\;\exists\;C_h> 0\;\forall\;u\in[0,\mathbf{N}_{\iota}):\;\;\;\exp(\omega_{\mathbf{M}}(u/h))\le\frac{C_hM_0}{N_0}\exp(\omega_{\mathbf{N}}(u));$$
thus the estimate in \eqref{sequencenotslowlyequvar} is shown with $t:=\frac{1}{h}$ and $D_t=D_h:=\log\left(\frac{C_hM_0}{N_0}\right)$ and note that for all values $t,u$ under consideration both sides of the estimate are finite by the comments above. Since $h>\underline{C}_{\mathbf{N}\preceq\mathbf{M}}$ is arbitrary, \eqref{sequencenotslowlyequvar} is verified with $t_0\ge(\underline{C}_{\mathbf{N}\preceq\mathbf{M}})^{-1}$.
\qed\enddemo

Since \hyperlink{mtriangle}{$\vartriangleleft$} is not reflexive, by combining $(I)$ in Proposition \ref{sequencenotslowly}, Lemma \ref{uppertransformfinite} and finally Remark \ref{slowlyvarrem} we get:

\begin{corollary}\label{sequencenotslowlycor}
Let $\mathbf{M}\in\RR_{>0}^{\NN}$ be log-convex and such that $\lim_{p\rightarrow+\infty}(M_p)^{1/p}=\mathbf{M}_{\iota}=+\infty$. Then \eqref{sequencenotslowlyequ} fails for $\omega_{\mathbf{M}}$ and hence $\omega_{\mathbf{M}}\widehat{\star}\omega_{\mathbf{M}}$ is not a weight function.
\end{corollary}

\begin{remark}
\emph{On the other hand it is not excluded that $\omega_{\mathbf{M}}$ is slowly varying. The quotient in $(ii)$ in Remark \ref{slowlyvarrem} for $\omega_{\mathbf{M}}$ should be compared with the quotient studied in the weighted entire setting in \cite[Lemma 2.9]{weightedentireinclusion1} (see also \cite[Lemma 2.3 $(ii)$]{weightedentireinclusion2}): This result shows, in particular, that the function $v_{\mathbf{M}}:=\exp\circ-\omega_{\mathbf{M}}$ is not slowly varying for any log-convex sequence $\mathbf{M}$ with $\mathbf{M}_{\iota}=+\infty$ (and $M_0=1$); see \cite[Def. 2.4]{weightedentireinclusion1}, \cite[Def. 2.2]{weightedentireinclusion2}. However, via \cite[Sect. 2.5, $(2.14)$]{weightedentireinclusion1} in this setting the weights $v_{\mathbf{M}}$ are defined differently and show different growth behavior compared with $\omega_{\mathbf{M}}$.}
\end{remark}

By taking into account the information from Lemma \ref{uppertransformfinite} and Proposition \ref{sequencenotslowly}, analogously as in \eqref{uppertransformgen} we give the following definition.

\begin{definition}\label{uppertransformdef}
Let log-convex sequences $\mathbf{M},\mathbf{N}\in\RR_{>0}^{\NN}$ be given. Assume that
\begin{itemize}
\item[$(a)$] $\mathbf{M}_{\iota}=+\infty=\mathbf{N}_{\iota}$ and

\item[$(b)$] $\mathbf{N}\hyperlink{mtriangle}{\vartriangleleft}\mathbf{M}$.
\end{itemize}
Then set
\begin{equation}\label{uppertransform}
\omega_{\mathbf{M}}\widehat{\star}\omega_{\mathbf{N}}(0)=0,\hspace{15pt}\omega_{\mathbf{M}}\widehat{\star}\omega_{\mathbf{N}}(t):=\sup_{s\ge 0}\{\omega_{\mathbf{M}}(s)-\omega_{\mathbf{N}}(s/t)\},\;\;\;t\in(0,+\infty).
\end{equation}
If instead of $(b)$ one has the weaker relation
\begin{itemize}
\item[$(c)$] $\mathbf{N}\hyperlink{preceq}{\preceq}\mathbf{M}$,
\end{itemize}
then
\begin{equation*}
\omega_{\mathbf{M}}\widehat{\star}\omega_{\mathbf{N}}(0)=0,\hspace{15pt}\omega_{\mathbf{M}}\widehat{\star}\omega_{\mathbf{N}}(t):=\sup_{s\ge 0}\{\omega_{\mathbf{M}}(s)-\omega_{\mathbf{N}}(s/t)\},\;\;\;t\in(0,(\underline{C}_{\mathbf{N}\preceq\mathbf{M}})^{-1}).
\end{equation*}
\end{definition}

\emph{Note:}

\begin{itemize}
\item[$(*)$] Assumptions $(a)$ and $(b)$ guarantee that $\omega_{\mathbf{M}}\widehat{\star}\omega_{\mathbf{N}}$ is a weight function and $\omega_{\mathbf{M}}\widehat{\star}\omega_{\mathbf{N}}(0)=\omega_{\mathbf{M}}(0)-\omega_{\mathbf{N}}(0)=0$ follows.

\item[$(*)$] If $(c)$ instead of $(b)$, then $\omega_{\mathbf{M}}\widehat{\star}\omega_{\mathbf{N}}$ can naturally only be defined on the interval $[0,(\underline{C}_{\mathbf{N}\preceq\mathbf{M}})^{-1})$ and taking $\omega_{\mathbf{M}}\widehat{\star}\omega_{\mathbf{N}}((\underline{C}_{\mathbf{N}\preceq\mathbf{M}})^{-1})=\lim_{t\rightarrow(\underline{C}_{\mathbf{N}\preceq\mathbf{M}})^{-1}}\omega_{\mathbf{M}}\widehat{\star}\omega_{\mathbf{N}}(t)$ since $\omega_{\mathbf{M}}\widehat{\star}\omega_{\mathbf{N}}(t)=+\infty$ for all $t>(\underline{C}_{\mathbf{N}\preceq\mathbf{M}})^{-1}$.

\item[$(*)$] In particular, this last comment applies to the case $\mathbf{M}=\mathbf{N}$ which gives that one should restrict to $t\in(0,1)$; compare also with $(i)$ in Remark \ref{slowlyvarrem}.
\end{itemize}

\begin{example}\label{uppertrafoGevreyex}
\emph{Let us apply Definition \eqref{uppertransformdef} to $\mathbf{N}=\mathbf{G}^{\alpha}$, $\mathbf{M}=\mathbf{G}^{\beta}$ with $\beta>\alpha>0$. Then clearly $(a)$ and $(b)$ are valid and so $\omega_{\mathbf{G^{\beta}}}\widehat{\star}\omega_{\mathbf{G}^{\alpha}}$ is a weight function. Indeed, when repeating the proof of Theorem \ref{uppertransformindexthm} for the associated weight functions and recalling that the growth indices $\gamma(\cdot)$ and $\overline{\gamma}(\cdot)$ are preserved under equivalence of weight functions we get}
\begin{equation}\label{uppertransformindexcorequvar}
\forall\;\beta>\alpha>0:\;\;\;\beta=\gamma(\omega_{\mathbf{G^{\beta}}}\widehat{\star}\omega_{\mathbf{G}^{\alpha}})+\alpha=\overline{\gamma}(\omega_{\mathbf{G^{\beta}}}\widehat{\star}\omega_{\mathbf{G}^{\alpha}})+\alpha.
\end{equation}
\end{example}

The next main result is the analogue of Theorem \ref{propBomanLemma4} and is also motivated by the previous example and the following observation: The product of two Gevrey-sequences yields again a Gevrey-sequence; however for the point-wise quotient sequence the indices have to be well-related (as in \eqref{uppertransformindexcorequvar}).

\begin{theorem}\label{propBomanLemma4inv}
Let log-convex sequences $\mathbf{M},\mathbf{N}\in\RR_{>0}^{\NN}$ be given such that
\begin{itemize}
\item[$(a)$] $\mathbf{M}_{\iota}=+\infty=\mathbf{N}_{\iota}$,

\item[$(b)$] $\mathbf{N}\hyperlink{mtriangle}{\vartriangleleft}\mathbf{M}$.
\end{itemize}
Then
\begin{equation}\label{propBomanLemma4invequ}
\forall\;t\in[0,+\infty):\;\;\;\omega_{\mathbf{M}}\widehat{\star}\omega_{\mathbf{N}}(t)\le\omega_{\frac{\mathbf{M}}{\mathbf{N}}}(t)
\end{equation}
holds and, if $\frac{\mathbf{M}}{\mathbf{N}}$ is even log-convex, then in \eqref{propBomanLemma4invequ} equality holds for all $t\in[0,+\infty)=[0,\frac{\mathbf{M}}{\mathbf{N}}_{\iota})$.\vspace{6pt}

If $(b)$ is replaced by the weaker relation
\begin{itemize}
\item[$(c)$] $\mathbf{N}\hyperlink{preceq}{\preceq}\mathbf{M}$,
\end{itemize}
then in \eqref{propBomanLemma4invequ} naturally one has to restrict to $t\in[0,(\underline{C}_{\mathbf{N}\preceq\mathbf{M}})^{-1})=[0,\frac{\mathbf{M}}{\mathbf{N}}_{\iota})$.
\end{theorem}

Concerning the second part of this result we point out:

\begin{itemize}
\item[$(*)$] In view of Lemma \ref{uppertransformfinite} and $(II)$ in Proposition \ref{sequencenotslowly} for the left-hand side, and in view of $(ii)$ in Lemma \ref{assofunctsectlemma} and comment $(b)$ in Section \ref{assofunctsect} applied to the sequence $\frac{\mathbf{M}}{\mathbf{N}}$ for the right-hand side, we have that \eqref{propBomanLemma4invequ} can be extended to all $t\in[0,+\infty)$ when putting formally both sides equal to $+\infty$ for all $t>\frac{\mathbf{M}}{\mathbf{N}}_{\iota}$ and for $t=\frac{\mathbf{M}}{\mathbf{N}}_{\iota}$ take the limit on both sides. The given arguments verify, in particular, that both limits exist and coincide.

\item[$(*)$] Taking $\mathbf{M}=\mathbf{N}$ gives $\frac{\mathbf{M}}{\mathbf{N}}=\mathbf{1}$, $\mathbf{1}:=(1)_{p\in\NN}$. Thus $\mathbf{1}_{\iota}=1$ and this special and degenerate situation corresponds to $\omega_{\frac{\mathbf{M}}{\mathbf{N}}}(t)=0$ for all $t\in[0,1]$ and $\omega_{\frac{\mathbf{M}}{\mathbf{N}}}(t)=+\infty$ for all $t>1$; see again $(ii)$ in Lemma \ref{assofunctsectlemma} and Remark \ref{lemma1stabrem} and the comments above related to Definition \ref{uppertransformdef}.
\end{itemize}

\demo{Proof}
First, concerning the equality $\frac{\mathbf{M}}{\mathbf{N}}_{\iota}=+\infty$ remark that by log-convexity of $\frac{\mathbf{M}}{\mathbf{N}}$ and $\mathbf{N}\hyperlink{mtriangle}{\vartriangleleft}\mathbf{M}$ we have $\frac{\mathbf{M}}{\mathbf{N}}_{\iota}=\lim_{p\rightarrow+\infty}\left(\frac{M_p}{N_p}\right)^{1/p}=+\infty$.\vspace{6pt}

For $t=0$ in \eqref{propBomanLemma4invequ} one has the equality $0=0$ and so consider from now on $t>0$. Then, by definition, $\exp(\omega_{\mathbf{M}}\widehat{\star}\omega_{\mathbf{N}}(t))=\sup_{s\ge 0}\frac{\sup_{p\in\NN}s^pM_0/M_p}{\sup_{q\in\NN}(s/t)^qN_0/N_q}$ and $\exp(\omega_{\frac{\mathbf{M}}{\mathbf{N}}}(t))=\sup_{\ell\in\NN}\frac{t^{\ell}M_0/N_0}{M_{\ell}/N_{\ell}}$.

For any $t>0$, $s\ge 0$ and $p\in\NN$ write $\frac{t^pM_0/N_0}{M_p/N_p}=\frac{s^pM_0}{M_p}\left(\frac{(s/t)^pN_0}{N_p}\right)^{-1}$ and note: When $(a_p)_{p\in\NN}$ and $(b_p)_{p\in\NN}$ are sequences of positive real numbers write $a_p=\frac{a_p}{b_p}b_p$ and we have $\sup_{p\in\NN}a_p\le\sup_{q\in\NN}\frac{a_q}{b_q}\sup_{r\in\NN}b_r$. Then, by using this estimate and taking the supremum over all $s\ge 0$, we get by definition $\exp(\omega_{\mathbf{M}}\widehat{\star}\omega_{\mathbf{N}}(t))\le\exp(\omega_{\frac{\mathbf{M}}{\mathbf{N}}}(t))$ for all $t>0$.\vspace{6pt}

For the converse, let $t>0$ be such that $\frac{\mu_p}{\nu_p}\le t\le\frac{\mu_{p+1}}{\nu_{p+1}}$ for some $p\in\NN_{>0}$ and note that $\frac{\mathbf{M}}{\mathbf{N}}$ is log-convex by assumption now. Then choose $\theta\in[0,1]$ such that $t=(\frac{\mu_p}{\nu_p})^{1-\theta}(\frac{\mu_{p+1}}{\nu_{p+1}})^{\theta}$ and set $s_t:=\mu_p^{1-\theta}\mu_{p+1}^{\theta}$. Hence $\mu_p\le s_t\le\mu_{p+1}$, $\frac{s_t}{t}=\nu_p^{1-\theta}\nu_{p+1}^{\theta}$ and so $\nu_p\le\frac{s_t}{t}\le\nu_{p+1}$. Then write $$\log\left(\frac{t^pM_0/N_0}{M_p/N_p}\right)=\log\left(\frac{s_t^pM_0}{M_p}\right)-\log\left(\frac{(s_t/t)^pN_0}{N_p}\right),$$
and by \eqref{lemma1equ} we have $$\omega_{\frac{\mathbf{M}}{\mathbf{N}}}(t)=\omega_{\mathbf{M}}(s_t)-\omega_{\mathbf{N}}(s_t/t)\le\sup_{s\ge 0}\{\omega_{\mathbf{M}}(s)-\omega_{\mathbf{N}}(s/t)\}.$$
Thus $\omega_{\frac{\mathbf{M}}{\mathbf{N}}}(t)\le\omega_{\mathbf{M}}\widehat{\star}\omega_{\mathbf{N}}(t)$ and hence equality is verified for all $t\ge\frac{\mu_1}{\nu_1}$.

Finally, for $0<t<\frac{\mu_1}{\nu_1}$ again in view of \eqref{lemma1equ} applied to the log-convex sequence $\frac{\mathbf{M}}{\mathbf{N}}$ we obtain $\omega_{\frac{\mathbf{M}}{\mathbf{N}}}(t)=0$ and $0\le\omega_{\mathbf{M}}\widehat{\star}\omega_{\mathbf{N}}(t)$ for all $t>0$ because $\omega_{\mathbf{N}}(0)=0$; recall Lemma \ref{uppertransformnonnegative}.\vspace{6pt}

When considering assumption $(c)$ instead of $(b)$, we use the information from $(II)$ in Proposition \ref{sequencenotslowly} and $(c)$ in Definition \ref{uppertransformdef}. Everything follows analogously by restricting to all $t\in[0,(\underline{C}_{\mathbf{N}\preceq\mathbf{M}})^{-1})$ and recall that log-convexity for $\frac{\mathbf{M}}{\mathbf{N}}$ yields $\frac{\mathbf{M}}{\mathbf{N}}_{\iota}=\liminf_{p\rightarrow+\infty}\left(\frac{M_p/M_0}{N_p/N_0}\right)^{1/p}=\lim_{p\rightarrow+\infty}\left(\frac{M_p/M_0}{N_p/N_0}\right)^{1/p}=(\underline{C}_{\mathbf{N}\preceq\mathbf{M}})^{-1}$.
\qed\enddemo

\subsection{Supplementary results}\label{supplsection}
In this section we are studying how crucial and useful growth properties for weight sequences and relations between them are preserved under $(\mathbf{M},\mathbf{N})\mapsto\mathbf{M}\cdot\mathbf{N}$ and $(\mathbf{M},\mathbf{N})\mapsto\frac{\mathbf{M}}{\mathbf{N}}$. The relevance of this goal follows from the shown equalities in the main results Theorems \ref{propBomanLemma4} and \ref{propBomanLemma4inv}.

From now on, in this Section, let us assume that $\mathbf{M},\mathbf{N}\in\hyperlink{LCset}{\mathcal{LC}}$ which is only a subtle technical restriction when dealing with log-convex sequences satisfying $\mathbf{M}_{\iota}=+\infty=\mathbf{N}_{\iota}$ and not effecting the definition of weighted spaces; recall Remark \ref{Cinfremark}. Consequently, one has $\mathbf{M}\cdot\mathbf{N}\in\hyperlink{LCset}{\mathcal{LC}}$ too.

\begin{itemize}
\item[$(i)$] In \cite[Thm. 3.1]{subaddlike} it has been shown that $\omega_{\mathbf{M}}$ satisfies \hyperlink{om1}{$(\omega_1)$} (i.e. $\gamma(\omega_{\mathbf{M}})>0$) if and only if
\begin{equation}\label{liminfom1}
\exists\;Q\in\NN_{\ge 2}:\;\;\;\liminf_{p\rightarrow+\infty}\frac{(M_{Qp})^{1/(Qp)}}{(M_p)^{1/p}}>1.
\end{equation}
\item[$(ii)$] Concerning \hyperlink{om6}{$(\omega_6)$} for $\omega_{\mathbf{M}}$ (i.e. $\overline{\gamma}(\omega_{\mathbf{M}})<+\infty$) several equivalent reformulations exist; we refer to \cite[Thm. 3.1]{modgrowthstrange} and the citations given in this paper. Indeed, \hyperlink{om6}{$(\omega_6)$} for $\omega_{\mathbf{M}}$ holds if and only if $\mathbf{M}$ satisfies \emph{moderate growth (mg)}
    $$\exists\;C\ge 1\;\forall\;p,q\in\NN:\;\;\;M_{p+q}\le C^{p+q}M_pM_q,$$
    equivalently
    \begin{equation}\label{mgquotient}
    \sup_{p\in\NN}\frac{\mu_{2p}}{\mu_p}<+\infty.
    \end{equation}
    In the literature (mg) is also denoted by (M.2); this notion is due to H. Komatsu in \cite{Komatsu73}.
\end{itemize}

\begin{lemma}\label{om1om6preservelemma}
Let $\mathbf{M},\mathbf{N}\in\hyperlink{LCset}{\mathcal{LC}}$ be given.
\begin{itemize}
\item[$(i)$] If $\mathbf{M},\mathbf{N}$ satisfy \eqref{liminfom1}, then also $\mathbf{M}\cdot\mathbf{N}$ and if $\mathbf{M},\mathbf{N}$ satisfy \emph{(mg),} then also $\mathbf{M}\cdot\mathbf{N}$.

\item[$(ii)$] If $\mathbf{M}$ satisfies
\begin{equation}\label{om1om6preservelemmaequ}
\exists\;Q\in\NN_{\ge 2}:\;\;\;\lim_{p\rightarrow+\infty}\frac{(M_{Qp})^{1/(Qp)}}{(M_p)^{1/p}}=+\infty,
\end{equation}
and $\mathbf{N}$ satisfies \emph{(mg),} then $\frac{\mathbf{M}}{\mathbf{N}}$ has \eqref{liminfom1}.

\item[$(iii)$] If $\mathbf{M}$ satisfies \emph{(mg),} $\mathbf{N}$ satisfies \eqref{liminfom1} and $\frac{\mathbf{M}}{\mathbf{N}}$ is log-convex, then $\frac{\mathbf{M}}{\mathbf{N}}$ has \emph{(mg).}
\end{itemize}
\end{lemma}

\demo{Proof}
$(i)$ Concerning \eqref{liminfom1} note that if this condition holds for some $Q_1\in\NN_{\ge 2}$, then also for all $Q\ge Q_1$ because $p\mapsto(M_p)^{1/p}$ is non-decreasing by $M_0=1$ and log-convexity. If now $\mathbf{M}$ satisfies \eqref{liminfom1} with $Q_1$ and $\mathbf{N}$ with $Q_2$, then with $Q:=\max\{Q_1,Q_2\}$ one has
$$\liminf_{p\rightarrow+\infty}\frac{(M_{Qp}\cdot N_{Qp})^{1/(Qp)}}{(M_p\cdot N_p)^{1/p}}\ge\liminf_{p\rightarrow+\infty}\frac{(M_{Qp})^{1/(Qp)}}{(M_p)^{1/p}}\liminf_{p\rightarrow+\infty}\frac{(N_{Qp})^{1/(Qp)}}{(N_p)^{1/p}},$$
and the conclusion follows from this estimate.

Concerning (mg), everything is immediate by $\sup_{p\in\NN}\frac{\mu_{2p}\nu_{2p}}{\mu_p\nu_p}\le\sup_{p\in\NN}\frac{\mu_{2p}}{\mu_p}\sup_{p\in\NN}\frac{\nu_{2p}}{\nu_p}$ and \eqref{mgquotient}.\vspace{6pt}

$(ii)$ Let \eqref{om1om6preservelemmaequ} for $\mathbf{M}$ be valid with some $Q\in\NN_{>0}$ (and hence for all $Q_1\ge Q$, too). Since $\mathbf{N}$ satisfies (mg) one has, by iteration,
\begin{equation}\label{mgiterate}
\exists\;C\ge 1\;\forall\;Q\in\NN_{>0}\;\forall\;p\in\NN:\;\;\;N_{Qp}\le C^{pQ(Q+1)/2-1}(N_p)^Q\le C^{2pQ^2}(N_p)^Q.
\end{equation}
By \eqref{om1om6preservelemmaequ} for $\mathbf{M}$ we can find $p_0$ depending on $Q$ and $C$ such that for all $p\ge p_0$: $$\frac{(M_{Qp}/N_{Qp})^{1/(Qp)}}{(M_p/N_p)^{1/p}}\ge\frac{(M_{Qp})^{1/(Qp)}}{(M_p)^{1/p}}\frac{1}{C^{2Q}}>1.$$

$(iii)$ Since $\mathbf{N}$ satisfies \eqref{liminfom1} we get
$$\exists\;Q\in\NN_{\ge 2}\;\exists\;\epsilon>0\;\exists\;p_{\epsilon}\in\NN_{>0}\;\forall\;p\ge p_{\epsilon}:\;\;\;\frac{(N_{Qp})^{1/(Qp)}}{(N_p)^{1/p}}\ge 1+\epsilon.$$
Hence, when applying \eqref{mgiterate} for $\mathbf{M}$ to this $Q$ we estimate by $\frac{(M_{Qp}/N_{Qp})^{1/(Qp)}}{(M_p/N_p)^{1/p}}\le C^{2Q}(1+\epsilon)^{-1}$ for all $p\ge p_{\epsilon}$ and when choosing a sufficiently large constant $D$ this verifies $\sup_{p\in\NN_{>0}}\frac{(M_{Qp}/N_{Qp})^{1/(Qp)}}{(M_p/N_p)^{1/p}}\le D$. By assumption $\frac{\mathbf{M}}{\mathbf{N}}$ is log-convex and then the technical result \cite[Lemma 2.2 $(ii)$]{subaddlike} applied to $\frac{\mathbf{M}}{\mathbf{N}}$ yields the conclusion. (Indeed, formally in this result one shall assume that $\frac{\mathbf{M}}{\mathbf{N}}\in\hyperlink{LCset}{\mathcal{LC}}$ but a switch to an equivalent and normalized sequence $\widetilde{\frac{\mathbf{M}}{\mathbf{N}}}\in\hyperlink{LCset}{\mathcal{LC}}$ yields the conclusion since the boundedness of the above expression in the supremum and (mg) are preserved under equivalence.)
\qed\enddemo

Combining the above comments and Lemma \ref{om1om6preservelemma} with Theorems \ref{propBomanLemma4} and \ref{propBomanLemma4inv} we immediately obtain:

\begin{corollary}\label{mainthmsindexcor}
Let $\mathbf{M},\mathbf{N}\in\hyperlink{LCset}{\mathcal{LC}}$ be given.
\begin{itemize}
\item[$(i)$] If $\gamma(\omega_{\mathbf{M}})>0$ and $\gamma(\omega_{\mathbf{N}})>0$, then $\gamma(\omega_{\mathbf{M}\cdot\mathbf{N}})=\gamma(\omega_{\mathbf{M}}\check{\star}\omega_{\mathbf{N}})>0$ and if $\overline{\gamma}(\omega_{\mathbf{M}})<+\infty$ and $\overline{\gamma}(\omega_{\mathbf{N}})<+\infty$, then $\overline{\gamma}(\omega_{\mathbf{M}\cdot\mathbf{N}})=\overline{\gamma}(\omega_{\mathbf{M}}\check{\star}\omega_{\mathbf{N}})<+\infty$.

\item[$(ii)$] Assume that $\frac{\mathbf{M}}{\mathbf{N}}$ is log-convex and $\mathbf{N}\hyperlink{mtriangle}{\vartriangleleft}\mathbf{M}$. If $\mathbf{M}$ satisfies \eqref{om1om6preservelemmaequ} and $\overline{\gamma}(\omega_{\mathbf{N}})<+\infty$, then $\gamma(\omega_{\mathbf{M}}\widehat{\star}\omega_{\mathbf{N}})=\gamma(\omega_{\frac{\mathbf{M}}{\mathbf{N}}})>0$.

\item[$(iii)$] Assume that $\frac{\mathbf{M}}{\mathbf{N}}$ is log-convex and $\mathbf{N}\hyperlink{mtriangle}{\vartriangleleft}\mathbf{M}$. If $\overline{\gamma}(\omega_{\mathbf{M}})<+\infty$ and $\gamma(\omega_{\mathbf{N}})>0$, then $\overline{\gamma}(\omega_{\mathbf{M}}\widehat{\star}\omega_{\mathbf{N}})=\overline{\gamma}(\omega_{\frac{\mathbf{M}}{\mathbf{N}}})<+\infty$.
\end{itemize}
\end{corollary}

\emph{Note:} $(i)$ in Corollary \ref{mainthmsindexcor} gives in the weight sequence setting an independent proof of the positivity resp. finiteness of the growth indices and hence a partial conclusion of  Theorem \ref{lowertransformindexthm}. (But there more precise information on the relations between the indices is obtained as well.)\vspace{6pt}

Let $\sigma,\tau$ be weight functions and recall that $\sigma,\tau\hyperlink{ompreceq}{\preceq}\sigma\check{\star}\tau$ holds, see $(a)$ in Lemma \ref{wedgeproplemma}, and $\sigma\widehat{\star}\tau\hyperlink{ompreceq}{\preceq}\sigma$ is valid by $(a)$ in Lemma \ref{widehatproplemma}. In view of these observations and Example(s) \ref{logexamplelower}, \ref{logexampleupper} one is interested in the following question:

\begin{itemize}
\item[$(*)$] For which weight functions $\sigma,\tau$ does one have $\sigma\check{\star}\tau\hyperlink{ompreceq}{\preceq}\sigma,\tau$ resp. $\sigma\hyperlink{ompreceq}{\preceq}\sigma\widehat{\star}\tau$; in particular for which weight functions it holds that $\sigma\check{\star}\sigma\hyperlink{sim}{\sim}\sigma$, $\sigma\widehat{\star}\sigma\hyperlink{sim}{\sim}\sigma$?
\end{itemize}

In the weight sequence setting, by the main results Theorem \ref{propBomanLemma4} and Theorem \ref{propBomanLemma4inv} and by taking into account the auxiliary result \cite[Lemma 6.5]{PTTvsmatrix} we get some information. Recall that the weight $\log_{+}$ treated in Examples \ref{logexamplelower}, \ref{logexampleupper} is an extreme case for the weight sequence setting and hence it is not considered in the following.

For any $\mathbf{M}\in\hyperlink{LCset}{\mathcal{LC}}$ one has that $p\mapsto M_p$ is non-decreasing, in particular $M_p\ge 1$ for all $p\in\NN$, and that
\begin{equation}\label{Lemma206equ}
\forall\;c\in\NN_{>0}\;\forall\;p\in\NN:\;\;\;(M_p)^c\le M_{cp};
\end{equation}
for a proof see e.g. \cite[Lemma 2.0.6]{diploma}.

\begin{proposition}\label{multlemma65}
Let $\mathbf{M},\mathbf{N}\in\hyperlink{LCset}{\mathcal{LC}}$ be given and all relations below are understood as $t\rightarrow+\infty$.
\begin{itemize}
\item[$(i)$] $\omega_{\mathbf{M}}\check{\star}\omega_{\mathbf{N}}(t)=O(\omega_{\mathbf{M}}(t)),O(\omega_{\mathbf{N}}(t))$ is always valid and hence $\omega_{\mathbf{M}}(t),\omega_{\mathbf{N}}(t)=o(\omega_{\mathbf{M}}\check{\star}\omega_{\mathbf{N}}(t))$ fails.

\item[$(ii)$] $\omega_{\mathbf{M}}\check{\star}\omega_{\mathbf{N}}(t)=o(\omega_{\mathbf{M}}(t))$ holds if $\mathbf{M}$ satisfies \emph{(mg)} and $\omega_{\mathbf{M}}\check{\star}\omega_{\mathbf{N}}(t)=o(\omega_{\mathbf{N}}(t))$ if $\mathbf{N}$ satisfies \emph{(mg).} Consequently, $\omega_{\mathbf{M}}(t)=O(\omega_{\mathbf{M}}\check{\star}\omega_{\mathbf{N}}(t))$ fails if $\mathbf{M}$ satisfies \emph{(mg)} and $\omega_{\mathbf{N}}(t)=O(\omega_{\mathbf{M}}\check{\star}\omega_{\mathbf{N}}(t))$ fails if $\mathbf{N}$ satisfies \emph{(mg).}
\end{itemize}
\end{proposition}

\emph{Note:} $(i)$ does also follow by resp. is consistent with Theorem \ref{propBomanLemma4} and $(a)$ in Lemma \ref{wedgeproplemma}.

\demo{Proof}
First, by Theorem \ref{propBomanLemma4} it holds that
$$\forall\;t\in[0,+\infty):\;\;\;\omega_{\mathbf{M}\cdot\mathbf{N}}(t)=\omega_{\mathbf{M}}\check{\star}\omega_{\mathbf{N}}(t).$$

$(i)$ By \cite[Lemma 6.5 $(i)\Leftrightarrow(ii)$]{PTTvsmatrix} we see that $\omega_{\mathbf{M}\cdot\mathbf{N}}(t)=O(\omega_{\mathbf{M}}(t))$ as $t\rightarrow+\infty$ if and only if
$$\exists\;c\in\NN_{>0}\;\exists\;A\ge 1\;\forall\;p\in\NN:\;\;\;M_p\le A(M_{cp}N_{cp})^{1/c},$$
which is clear with $A:=1$ by \eqref{Lemma206equ} and since $N_q\ge 1$ for any $q\in\NN$. Similarly, relation $\omega_{\mathbf{M}\cdot\mathbf{N}}(t)=O(\omega_{\mathbf{N}}(t))$ is always true.\vspace{6pt}

$(ii)$ By \cite[Lemma 6.5 $(i)'\Leftrightarrow(ii)'$]{PTTvsmatrix} one has that $\omega_{\mathbf{M}\cdot\mathbf{N}}(t)=o(\omega_{\mathbf{M}}(t))$ as $t\rightarrow+\infty$ if and only if
$$\forall\;c\in\NN_{>0}\;\exists\;A\ge 1\;\forall\;p\in\NN:\;\;\;(M_{cp})^{1/c}\le A M_pN_p.$$
If $\mathbf{M}$ satisfies (mg), then this estimate is valid since $(M_{cp})^{1/c}\le B^pM_p\le AM_pN_p$ for some $A,B\ge 1$ and all $p\in\NN$ because $\lim_{p\rightarrow+\infty}(N_p)^{1/p}=+\infty$ (and $A$ is then depending on chosen $c$). Similarly, this is valid if $\mathbf{M}$ is replaced by $\mathbf{N}$.


\qed\enddemo

\begin{proposition}\label{divlemma65}
Let $\mathbf{M},\mathbf{N}\in\hyperlink{LCset}{\mathcal{LC}}$ be given such that $\mathbf{N}\hyperlink{mtriangle}{\vartriangleleft}\mathbf{M}$ and such that $\frac{\mathbf{M}}{\mathbf{N}}$ is log-convex. All relations below are understood as $t\rightarrow+\infty$.
\begin{itemize}
\item[$(i)$] $\omega_{\mathbf{M}}(t)=o(\omega_{\mathbf{M}}\widehat{\star}\omega_{\mathbf{N}}(t))$ holds if $\mathbf{M}$ satisfies \emph{(mg).} Consequently, $\omega_{\mathbf{M}}\widehat{\star}\omega_{\mathbf{N}}(t)=O(\omega_{\mathbf{M}}(t))$ is violated if $\mathbf{M}$ satisfies \emph{(mg).}

\item[$(ii)$] $\omega_{\mathbf{M}}(t)=O(\omega_{\mathbf{M}}\widehat{\star}\omega_{\mathbf{N}}(t))$ is valid. Consequently, $\omega_{\mathbf{M}}\widehat{\star}\omega_{\mathbf{N}}(t)=o(\omega_{\mathbf{M}}(t))$ fails.
\end{itemize}
\end{proposition}
\emph{Note:} $(ii)$ does also follow by resp. is consistent with Theorem \ref{propBomanLemma4inv} and $(a)$ in Lemma \ref{widehatproplemma}.

\demo{Proof}
First, by Theorem \ref{propBomanLemma4inv} it follows that
$$\forall\;t\in[0,+\infty):\;\;\;\omega_{\mathbf{M}}\widehat{\star}\omega_{\mathbf{N}}(t)=\omega_{\frac{\mathbf{M}}{\mathbf{N}}}(t),$$
and $\frac{\mathbf{M}}{\mathbf{N}}_{\iota}=+\infty$.

$(i)$ For $\omega_{\mathbf{M}}(t)=o(\omega_{\frac{\mathbf{M}}{\mathbf{N}}}(t))$ as $t\rightarrow+\infty$ one requires
$$\forall\;c\in\NN_{>0}\;\exists\;A\ge 1\;\forall\;p\in\NN:\;\;\;\left(\frac{M_{cp}}{N_{cp}}\right)^{1/c}\le AM_p.$$
Let now $c\in\NN_{>0}$ be given and fixed. If $\mathbf{M}$ satisfies (mg), then via iteration $(M_{cp})^{1/c}\le B^pM_p$ for some $B\ge 1$ and all $p\in\NN$; see \eqref{mgiterate}. Since $\lim_{p\rightarrow+\infty}(N_p)^{1/p}=+\infty$ it holds that $B^p\le AN_p$ for some $A\ge 1$ depending on $B$, and hence on $c$, and for all $p\in\NN$. Summarizing, $(M_{cp})^{1/c}\le AM_pN_p\le AM_p(N_{cp})^{1/c}$ for all $p\in\NN$ and the last estimate holds by \eqref{Lemma206equ} applied to $\mathbf{N}$.\vspace{6pt}

$(ii)$ For $\omega_{\mathbf{M}}(t)=O(\omega_{\frac{\mathbf{M}}{\mathbf{N}}}(t))$, one requires
$$\exists\;c\in\NN_{>0}\;\exists\;A\ge 1\;\forall\;p\in\NN:\;\;\;\frac{M_p}{N_p}\le A(M_{cp})^{1/c}.$$
Since $N_p\ge 1$ for all $p\in\NN$ we see that \eqref{Lemma206equ} implies this desired estimate with $A:=1$.


\qed\enddemo

\section{Inverse operations}\label{inversection}
The aim is to study how $\check{\star}$ and $\widehat{\star}$ are related to each other; indeed Theorems \ref{propBomanLemma4} and \ref{propBomanLemma4inv} together suggest that $\check{\star}$ and $\widehat{\star}$ are inverse (at least in the weight sequence setting).

\subsection{The weight sequence setting}
From Theorems \ref{propBomanLemma4} and \ref{propBomanLemma4inv}, we infer:

\begin{theorem}\label{propBomanLemma4comb}
Let log-convex sequences $\mathbf{M},\mathbf{N}\in\RR_{>0}^{\NN}$ be given with $\mathbf{M}_{\iota}=+\infty=\mathbf{N}_{\iota}$.
\begin{itemize}
\item[$(i)$] Then
    $$\forall\;t\in[0,+\infty):\;\;\;(\omega_{\mathbf{M}}\check{\star}\omega_{\mathbf{N}})\widehat{\star}\omega_{\mathbf{N}}(t)=\omega_{\mathbf{M}}(t),$$
    and
    $$\forall\;t\in[0,+\infty):\;\;\;(\omega_{\mathbf{M}}\check{\star}\omega_{\mathbf{N}})\widehat{\star}\omega_{\mathbf{M}}(t)=\omega_{\mathbf{N}}(t).$$

\item[$(ii)$] Assume that $\mathbf{N}\hyperlink{mtriangle}{\vartriangleleft}\mathbf{M}$ and that $\frac{\mathbf{M}}{\mathbf{N}}$ is log-convex. Then we obtain
$$\forall\;t\in[0,+\infty):\;\;\;\omega_{\mathbf{M}}(t)=\omega_{\mathbf{N}}\check{\star}(\omega_{\mathbf{M}}\widehat{\star}\omega_{\mathbf{N}})(t)=(\omega_{\mathbf{M}}\widehat{\star}\omega_{\mathbf{N}})\check{\star}\omega_{\mathbf{N}}(t).$$
\end{itemize}
\end{theorem}

\demo{Proof}
$(i)$ First, Theorem \ref{propBomanLemma4} implies $\omega_{\mathbf{M}\cdot\mathbf{N}}(t)=\omega_{\mathbf{M}}\check{\star}\omega_{\mathbf{N}}(t)$ for all $t\in[0,+\infty)$.

Then distinguish: In the first case apply Theorem \ref{propBomanLemma4inv} to $\mathbf{P}=\mathbf{M}\cdot\mathbf{N}$ and $\mathbf{Q}=\mathbf{N}$. Then clearly $\frac{\mathbf{P}}{\mathbf{Q}}=\mathbf{M}$ is log-convex, $\mathbf{Q}\hyperlink{mtriangle}{\vartriangleleft}\mathbf{P}$ precisely means $\lim_{p\rightarrow+\infty}\frac{1}{(M_p/M_0)^{1/p}}=0$, i.e. $\lim_{p\rightarrow+\infty}(M_p/M_0)^{1/p}=\mathbf{M}_{\iota}=+\infty$, and which holds by assumption. Finally, one gets $(\underline{C}_{\mathbf{Q}\preceq\mathbf{P}})^{-1}=\frac{\mathbf{P}}{\mathbf{Q}}_{\iota}=\mathbf{M}_{\iota}=+\infty$.

Therefore, Theorem \ref{propBomanLemma4inv} implies $\omega_{\mathbf{M}\cdot\mathbf{N}}\widehat{\star}\omega_{\mathbf{N}}(t)=\omega_{\mathbf{P}}\widehat{\star}\omega_{\mathbf{Q}}(t)=\omega_{\mathbf{M}}(t)$ for all $t\in[0,+\infty)$ and combining both identities gives the desired equation.

On the other hand, in the second case apply Theorem \ref{propBomanLemma4inv} to $\mathbf{P}=\mathbf{M}\cdot\mathbf{N}$ and $\mathbf{Q}=\mathbf{M}$ and the rest follows as before.\vspace{6pt}

$(ii)$ First, Theorem \ref{propBomanLemma4inv} gives
\begin{equation}\label{propBomanLemma4combequ1}
\forall\;t\in\left[0,\frac{\mathbf{M}}{\mathbf{N}}_{\iota}\right)=[0,+\infty):\:\:\;\omega_{\mathbf{M}}\widehat{\star}\omega_{\mathbf{N}}(t)=\omega_{\frac{\mathbf{M}}{\mathbf{N}}}(t).
\end{equation}
Note that by log-convexity for $\frac{\mathbf{M}}{\mathbf{N}}$ and by $\mathbf{N}\hyperlink{mtriangle}{\vartriangleleft}\mathbf{M}$ we get $\frac{\mathbf{M}}{\mathbf{N}}_{\iota}=+\infty$. Then apply Theorem \ref{propBomanLemma4} to $\mathbf{P}=\mathbf{N}$, $\mathbf{Q}=\frac{\mathbf{M}}{\mathbf{N}}$ and note again that $\frac{\mathbf{M}}{\mathbf{N}}$ is log-convex and $\frac{\mathbf{M}}{\mathbf{N}}_{\iota}=+\infty$. Hence
\begin{equation}\label{propBomanLemma4combequ}
\forall\;t\in[0,+\infty):\;\;\;\omega_{\mathbf{M}}(t)=\omega_{\mathbf{P}\cdot\mathbf{Q}}(t)=\omega_{\mathbf{N}}\check{\star}\omega_{\frac{\mathbf{M}}{\mathbf{N}}}(t),
\end{equation}
since $t\in[0,\mathbf{P}_{\iota}\cdot\mathbf{Q}_{\iota})=[0,\mathbf{N}_{\iota}\cdot\frac{\mathbf{M}}{\mathbf{N}}_{\iota})=[0,+\infty)$. Recall that for $\omega_{\mathbf{N}}\check{\star}\omega_{\frac{\mathbf{M}}{\mathbf{N}}}(t)$ naturally one has to consider the infimum over all $s\in\mathcal{I}_t=(t(\frac{\mathbf{M}}{\mathbf{N}}_{\iota})^{-1},\mathbf{N}_{\iota})=(0,+\infty)$ with $t\in[0,\mathbf{N}_{\iota}\cdot\frac{\mathbf{M}}{\mathbf{N}}_{\iota})=[0,+\infty)$, see \eqref{lowertransform} and for this note that $\left(\frac{\mathbf{M}}{\mathbf{N}}_{\iota}\right)^{-1}=0$. Finally, by combining \eqref{propBomanLemma4combequ1} and \eqref{propBomanLemma4combequ} the desired equality is verified.
\qed\enddemo

\subsection{The general weight function case}
We investigate how Theorem \ref{propBomanLemma4comb} transfers to abstractly given weight functions.

\begin{proposition}\label{inversethm}
Let $\omega,\sigma,\tau$ be weight functions. Assume that $\tau(0)=0$ and that $\sigma\widehat{\star}\tau$ is well defined; i.e. \eqref{uppetransfromwelldefequ} holds with $t_0=+\infty$.

Then $\omega\check{\star}(\sigma\widehat{\star}\tau)$ is a weight function such that $\omega\check{\star}(\sigma\widehat{\star}\tau)(0) =\omega(0)+\sigma(0)(\ge 0)$ and one has:
\begin{itemize}
\item[$(i)$]
$$\forall\;t>0:\;\;\;\sigma(t)+\inf_{u>0}\{\omega(u)-\tau(u)\}\le\omega\check{\star}(\sigma\widehat{\star}\tau)(t),$$

\item[$(ii)$]
$$\forall\;t\ge 0:\;\;\;\omega\check{\star}(\sigma\widehat{\star}\tau)(t)\le\min\{\omega(t)+\sigma\widehat{\star}\tau(1),\sigma\widehat{\star}\tau(t)+\omega(1)\},$$
which implies, in particular,
\begin{equation}\label{inversethmequ}
\limsup_{t\rightarrow+\infty}\frac{\omega\check{\star}(\sigma\widehat{\star}\tau)(t)}{\min\{\omega(t),\sigma\widehat{\star}\tau(t)\}}\le 1,
\end{equation}
hence $\sigma\widehat{\star}\tau,\omega\hyperlink{ompreceq}{\preceq}\omega\check{\star}(\sigma\widehat{\star}\tau)$.
\end{itemize}
\end{proposition}

\demo{Proof}
First, by assumption $\sigma\widehat{\star}\tau$ is a weight function; recall also Lemma \ref{uppertransformfinite}. Thus $\omega\check{\star}(\sigma\widehat{\star}\tau)$ is also a weight function and by Remark \ref{lowerLegendrepoint0rem} and $(b)$ in Lemma \ref{widehatproplemma} we get $$\omega\check{\star}(\sigma\widehat{\star}\tau)(0)=\omega(0)+\sigma\widehat{\star}\tau(0)=\omega(0)+\sigma(0)-\tau(0)=\omega(0)+\sigma(0).$$

$(i)$ For any $t>0$ it follows that
\begin{align*}
\omega\check{\star}(\sigma\widehat{\star}\tau)(t)&=\inf_{u>0}\{\omega(u)+\sigma\widehat{\star}\tau(t/u)\}=\inf_{u>0}\{\omega(u)+\sup_{s\ge 0}\{\sigma(s)-\tau(\frac{s}{t/u})\}\}
\\&
\underbrace{\ge}_{s:=t}\inf_{u>0}\{\omega(u)+\sigma(t)-\tau(u)\}.
\end{align*}
\emph{Note:} In general one only has $\inf_{u>0}\{\omega(u)-\tau(u)\}\in[-\infty,+\infty)$.\vspace{6pt}

$(ii)$ For any $t>0$ we estimate as follows, see also the proof of $(a)$ in Lemma \ref{wedgeproplemma}:
\begin{align*}
\omega\check{\star}(\sigma\widehat{\star}\tau)(t)&=\inf_{u>0}\{\omega(u)+\sigma\widehat{\star}\tau(t/u)\}\underbrace{\le}_{u:=t}\omega(t)+\sigma\widehat{\star}\tau(1).
\end{align*}

$\sigma\widehat{\star}\tau$ is a weight function and so, in particular, $\sigma\widehat{\star}\tau(1)<+\infty$. Since $\omega\check{\star}(\sigma\widehat{\star}\tau)(0)=\omega(0)+\sigma(0)=\omega(0)+\sigma\widehat{\star}\tau(0)\le\omega(0)+\sigma\widehat{\star}\tau(1)$ the above estimate is valid for all $t\ge 0$.

Analogously, the choice $u=1$ gives $\omega\check{\star}(\sigma\widehat{\star}\tau)(t)\le\omega(1)+\sigma\widehat{\star}\tau(t)$ for all $t\ge 0$ and for $t=0$ note that $\omega\check{\star}(\sigma\widehat{\star}\tau)(0)=\omega(0)+\sigma(0)\le\omega(1)+\sigma(0)=\omega(1)+\sigma\widehat{\star}\tau(0)$. Since $\omega$ and $\sigma\widehat{\star}\tau$ are weight functions, and so $\lim_{t\rightarrow+\infty}\omega(t)=\lim_{t\rightarrow+\infty}\sigma\widehat{\star}\tau(t)=+\infty$, we get \eqref{inversethmequ}.
\qed\enddemo

The previous result immediately gives:

\begin{corollary}\label{inversethmcor}
Let $\omega,\sigma,\tau$ be weight functions and assume that
\begin{itemize}
\item[$(i)$] $\tau(0)=0$;

\item[$(ii)$] $\sigma\widehat{\star}\tau$ is well defined;

\item[$(iii)$] $C_{\omega,\tau}:=\inf_{u>0}\{\omega(u)-\tau(u)\}>-\infty$; i.e.
$$\exists\;C>0\;\forall\;u\ge 0:\;\;\;\tau(u)\le\omega(u)+C.$$
\end{itemize}
Then it holds that
\begin{equation*}\label{inversethmcorequ}
\sigma\widehat{\star}\tau,\omega\hyperlink{ompreceq}{\preceq}\omega\check{\star}(\sigma\widehat{\star}\tau)\hyperlink{ompreceq}{\preceq}\sigma,
\end{equation*}
and, moreover that
\begin{equation}\label{inversethmcorequ1}
\forall\;t\ge 0:\;\;\;\sigma(t)\le\tau\check{\star}(\sigma\widehat{\star}\tau)(t)\le\min\{\tau(t)+\sigma\widehat{\star}\tau(1),\sigma\widehat{\star}\tau(t)+\tau(1)\}.
\end{equation}
\end{corollary}

\emph{Note:}

\begin{itemize}
\item[$(*)$] For \eqref{inversethmcorequ1} one considers $\omega=\tau$ and so $C_{\omega,\omega}=0$.

\item[$(*)$] In order to establish an equivalence in \eqref{inversethmcorequ1} between $\tau\check{\star}(\sigma\widehat{\star}\tau)$ and $\sigma$ one is interested in having $\sigma\hyperlink{ompreceq}{\preceq}\tau$ or $\sigma\hyperlink{ompreceq}{\preceq}\sigma\widehat{\star}\tau$ and Corollary \ref{inversethmcor} for $\omega=\tau$ should be compared with $(ii)$ in Theorem \ref{propBomanLemma4comb} with $\sigma=\omega_{\mathbf{M}}$ and $\tau=\omega_{\mathbf{N}}$.

\item[$(*)$] However, for the desired relations even in the weight sequence setting $\omega_{\mathbf{M}}\hyperlink{ompreceq}{\preceq}\omega_{\mathbf{M}}\widehat{\star}\omega_{\mathbf{N}}$ fails if $\mathbf{M}$ satisfies (mg); see $(i)$ in Proposition \ref{divlemma65}. On the other hand $\omega_{\mathbf{M}}\hyperlink{ompreceq}{\preceq}\omega_{\mathbf{N}}$, i.e. $\omega_{\mathbf{N}}(t)=O(\omega_{\mathbf{M}}(t))$, amounts to have
$$\exists\;c\in\NN_{>0}\;\exists\;A\ge 1\;\forall\;p\in\NN:\;\;\;M_p\le A(N_{cp})^{1/c},$$
see again \cite[Lemma 6.5 $(i)\Leftrightarrow(ii)$]{PTTvsmatrix}. In $(ii)$ in Theorem \ref{propBomanLemma4comb} the assumption $\mathbf{N}\hyperlink{mtriangle}{\vartriangleleft}\mathbf{M}$ is used and this gives that for all $h>0$ there exists $C_h\ge 1$ such that $N_p\le C_hh^pM_p$ for all $p\in\NN$. In view of \eqref{Lemma206equ} (for $\mathbf{N}$) both estimates can be valid simultaneously.
\end{itemize}

Concerning the analogue of $(i)$ in Theorem \ref{propBomanLemma4comb} for general weight functions we can show:

\begin{proposition}\label{inversethm1}
Let $\omega,\sigma,\tau$ be weight functions.
\begin{itemize}
\item[$(i)$] If $\tau(0)=0$, then (formally) $(\omega\check{\star}\sigma)\widehat{\star}\tau(t)\ge 0$ for all $t\ge 0$ and, moreover,
$$\forall\;t\ge 0:\;\;\;(\omega\check{\star}\sigma)\widehat{\star}\tau(t)\ge\omega\check{\star}\sigma(t)-\tau(1).$$

\item[$(ii)$] We get (formally)
$$\forall\;t>0:\;\;\;(\omega\check{\star}\sigma)\widehat{\star}\tau(t)\le\min\{\omega(t)+\sup_{u\ge 0}\{\sigma(u/t)-\tau(u/t)\},\sigma(t)+\sup_{u\ge 0}\{\omega(u/t)-\tau(u/t)\}\}.$$

\item[$(iii)$] Both $(\sigma\check{\star}\tau)\widehat{\star}\tau$ and $(\sigma\check{\star}\tau)\widehat{\star}\sigma$ are weight functions and
    $$\forall\;t\ge 0:\;\;\;\sigma\check{\star}\tau(t)-\tau(1)\le(\sigma\check{\star}\tau)\widehat{\star}\tau(t)\le\sigma(t),$$
    $$\forall\;t\ge 0:\;\;\;\sigma\check{\star}\tau(t)-\sigma(1)\le(\sigma\check{\star}\tau)\widehat{\star}\sigma(t)\le\tau(t).$$
\end{itemize}
\end{proposition}

\demo{Proof}
$(i)$ Let $t>0$, then
\begin{align*}
(\omega\check{\star}\sigma)\widehat{\star}\tau(t)&=\sup_{u\ge 0}\{\omega\check{\star}\sigma(u)-\tau(u/t)\}\underbrace{\ge}_{u=t}\omega\check{\star}\sigma(t)-\tau(1),
\end{align*}
and for $t=0$ one has $(\omega\check{\star}\sigma)\widehat{\star}\tau(0)=\omega(0)+\sigma(0)-\tau(0)=\omega(0)+\sigma(0)=\omega\check{\star}\sigma(0)(\ge 0)$ since $\tau(0)=0$. And since $\tau$ is non-decreasing, also $(\omega\check{\star}\sigma)\widehat{\star}\tau(0)\ge\omega\check{\star}\sigma(0)-\tau(1)$ follows.\vspace{6pt}

$(ii)$ Let $t>0$ be fixed and since $(\omega\check{\star}\sigma)\widehat{\star}\tau(t)=\sup_{u\ge 0}\{\inf_{s>0}\{\omega(s)+\sigma(u/s)\}-\tau(u/t)\}$ when choosing $s:=t$ in the infimum:
\begin{align*}
(\omega\check{\star}\sigma)\widehat{\star}\tau(t)\le\omega(t)+\sup_{u\ge 0}\{\sigma(u/t)-\tau(u/t)\}.
\end{align*}
On the other hand, $(\omega\check{\star}\sigma)\widehat{\star}\tau(t)=\sup_{u\ge 0}\{\inf_{v>0}\{\omega(u/v)+\sigma(v)\}-\tau(u/t)\}$ and choosing $v:=t$ in the infimum we infer the analogous estimate.\vspace{6pt}

$(iii)$ First, $\sigma\check{\star}\tau$ is a weight function and, second, like in $(i)$ one has $(\sigma\check{\star}\tau)\widehat{\star}\tau(0)=\sigma(0)\ge 0$ and $(\sigma\check{\star}\tau)\widehat{\star}\sigma(0)=\tau(0)\ge 0$. Therefore, $(\sigma\check{\star}\tau)\widehat{\star}\tau(t)\ge 0$ and $(\sigma\check{\star}\tau)\widehat{\star}\sigma(t)\ge 0$ (formally) for all $t\ge 0$.

More generally, in order to ensure that $(\omega\check{\star}\sigma)\widehat{\star}\tau$ is well defined, one requires (see \eqref{uppetransfromwelldefequ}):
\begin{equation}\label{uppertransformfiniteequtriple}
\forall\;t>0\;\exists\;D_t>0\;\forall\;u\ge 0:\;\;\;\omega\check{\star}\sigma(tu)\le\tau(u)+D_t.
\end{equation}
By \eqref{wedgeproplemmaequ} one has $\omega\check{\star}\sigma(tu)\le\min\{\omega(t)+\sigma(u),\omega(u)+\sigma(t)\}$ for all $t,u\ge 0$
and so \eqref{uppertransformfiniteequtriple} holds for $(\sigma\check{\star}\tau)\widehat{\star}\tau$ with $D_t:=\sigma(t)$ and for $(\sigma\check{\star}\tau)\widehat{\star}\sigma$ with $D_t:=\tau(t)$.

Finally, the desired estimates follow by $(i)$ and $(ii)$ and note that by $(\sigma\check{\star}\tau)\widehat{\star}\tau(0)=\sigma(0)$, $(\sigma\check{\star}\tau)\widehat{\star}\sigma(0)=\tau(0)$ the second inequalities also hold for $t=0$ (even with equality).
\qed\enddemo

\emph{Note:} When using the estimates in $(iii)$ in Proposition \ref{inversethm1}, in order to ensure $(\sigma\check{\star}\tau)\widehat{\star}\tau\hyperlink{sim}{\sim}\sigma$ resp. $(\sigma\check{\star}\tau)\widehat{\star}\sigma\hyperlink{sim}{\sim}\tau$ it suffices to have $\sigma(t)=O(\tau\check{\star}\sigma(t))$ resp. $\tau(t)=O(\tau\check{\star}\sigma(t))$ as $t\rightarrow+\infty$. But, by $(ii)$ in Proposition \ref{multlemma65}, this is in general even violated in the weight sequence setting.\vspace{6pt}

On the other hand, we can apply Proposition \ref{inversethm1} to show that the identity in $(b)$ in Lemma \ref{widehatproplemma} does not follow even if $\sigma\widehat{\star}\tau$ is a weight function; see also the comments in Section \ref{arbuppertrafosectprelim}.

\begin{corollary}\label{inversethm1cor}
For given weight functions $\omega$ and $\omega_1$, even if $\omega\widehat{\star}\omega_1$ is again a weight function, in general $\lim_{t\rightarrow 0^{+}}\omega\widehat{\star}\omega_1(t)\neq\omega(0)-\omega_1(0)=\omega\widehat{\star}\omega_1(0)$.
\end{corollary}

\demo{Proof}
Let $\sigma$, $\tau$ be weight functions such that $\sigma(0)<\underline{\sigma}:=\lim_{t\rightarrow 0^{+}}\sigma(t)$ and $\tau(t)=0$ for all $t\in[0,1]$ \emph{(``normalization'').} Then $(iii)$ in Proposition \ref{inversethm1} implies that $(\sigma\check{\star}\tau)\widehat{\star}\tau$ is a weight function and such that for all $t>0$:
\begin{align*}
(\sigma\check{\star}\tau)\widehat{\star}\tau(t)&\ge\sigma\check{\star}\tau(t)=\inf_{s>0}\{\sigma(s)+\tau(t/s)\}\ge\inf_{s>0}\{\sigma(s)\}=\underline{\sigma}>\sigma(0)=\sigma(0)+\tau(0)-\tau(0)
\\&
=(\sigma\check{\star}\tau)\widehat{\star}\tau(0).
\end{align*}
\qed\enddemo

\appendix
\section{On non-standard cases in the weight sequence setting}\label{nonstandardappendix}
In order to present a complete approach from the point of view of abstractly given (log-convex) sequences of positive real numbers and also motivated by the general formulation of \cite[Lemma 4]{Boman98} we are interested in non-standard cases and thus we even allow $\mathbf{M}_{\iota}<+\infty$ and/or $\mathbf{N}_{\iota}<+\infty$; i.e. $\mathbf{M}_{\iota}\in(0,+\infty)$ and/or $\mathbf{N}_{\iota}\in(0,+\infty)$. For the definition of weighted settings usually such sequences are not considered since they yield trivial and/or restricting cases. Note that if $\mathbf{M}_{\iota}<+\infty$ resp. $\mathbf{N}_{\iota}<+\infty$, then $\omega_{\mathbf{M}}$ resp. $\omega_{\mathbf{N}}$ is not a weight function in the sense of Definition \ref{weightfctdef}; recall $(ii)$ in Lemma \ref{assofunctsectlemma}. But, on the other hand, from an abstract point of view it is natural to study the operations $\check{\star}$, $\widehat{\star}$ between $\omega_{\mathbf{M}}$ and $\omega_{\mathbf{N}}$ for arbitrary $\mathbf{M},\mathbf{N}\in\RR_{>0}^{\NN}$. Recall that, apart from \cite{Boman98} also in \cite[Sect. 3.1]{Komatsu73} formally the identity $\mathbf{M}_{\iota}=+\infty$ is not ensured; see the comments in \cite[Sect. 2.5]{regularnew}.\vspace{6pt}

In order to proceed we are using from now on the conventions $0^0:=1$, $\frac{0}{0}:=0$, $\frac{c}{+\infty}=0$ for any $c>0$, $+\infty-(+\infty)=0$, $0\cdot(-\infty):=0=:0\cdot(+\infty)$ and $p\cdot(-\infty):=-\infty$ for any $p\in\NN_{>0}$.

\subsection{Non-standard cases for the generalized lower Legendre conjugate}\label{nonstandardlowersection}
Let $\mathbf{M},\mathbf{N}\in\RR_{>0}^{\NN}$ be given, then in general $\mathcal{I}_t\neq(0,+\infty)$ with $\mathcal{I}_t:=(t\mathbf{N}_{\iota}^{-1},\mathbf{M}_{\iota})$, $t\in[0,\mathbf{M}_{\iota}\cdot\mathbf{N}_{\iota})$ (see \eqref{Itset}) and for the next comments recall $(ii)$ in Lemma \ref{assofunctsectlemma} and Lemma \ref{wedgeproplemma}:

\begin{itemize}
\item[$(i)$] $\mathcal{I}_t\neq\emptyset$ is clear for $t=0$, for $t>0$ when $\mathbf{N}_{\iota}=+\infty$ (using $\frac{1}{+\infty}=0$) and also otherwise since $t<\mathbf{M}_{\iota}\cdot\mathbf{N}_{\iota}$.

    For each $s\in\mathcal{I}_t$ we have $\omega_{\mathbf{M}}(s)<+\infty$ because $s<\mathbf{M}_{\iota}$ and $\omega_{\mathbf{N}}(t/s)<+\infty$ because $\frac{t}{s}<\mathbf{N}_{\iota}$. Consequently, $\omega_{\mathbf{M}}\check{\star}\omega_{\mathbf{N}}(t)<+\infty$ for all $t\in[0,\mathbf{M}_{\iota}\cdot\mathbf{N}_{\iota})$.

Clearly, $\omega_{\mathbf{M}}\check{\star}\omega_{\mathbf{N}}$ is non-decreasing and one has $$\omega_{\mathbf{M}}\check{\star}\omega_{\mathbf{N}}(0)=\inf_{s\in\mathcal{I}_0}\{\omega_{\mathbf{M}}(s)\}=\inf_{s\in(0,\mathbf{M}_{\iota})}\{\omega_{\mathbf{M}}(s)\}=0;$$
compare this formula with Remark \ref{lowerLegendrepoint0rem}: Recall that $\omega_{\mathbf{M}}$ is continuous on $[0,\mathbf{M}_{\iota})$ on hence the obstruction mentioned in this remark cannot occur in the weight sequence setting.

\item[$(ii)$] We can consider in the infimum also the value $s=0$: Since $\omega_{\mathbf{M}}(0)=0$ and by using the conventions $\frac{1}{0}=+\infty$ and $\omega_{\mathbf{N}}(+\infty)=+\infty$ this will not change the value of $\omega_{\mathbf{M}}\check{\star}\omega_{\mathbf{N}}(t)$. Note that $\omega_{\mathbf{N}}(+\infty)=+\infty$ is justified by continuity and $\lim_{u\rightarrow+\infty}\omega_{\mathbf{N}}(u)=+\infty$ if $\mathbf{N}_{\iota}=+\infty$ and by $\omega_{\mathbf{N}}(u)=+\infty$ for all $u>\mathbf{N}_{\iota}$ if $\mathbf{N}_{\iota}<+\infty$; recall $(b)$ in Section \ref{assofunctsect}.

\item[$(iii)$] Assume $t>0$ and $\mathcal{I}_t\neq(0,+\infty)$, i.e. $\mathbf{N}_{\iota}\in(0,+\infty)$ or $\mathbf{M}_{\iota}\in(0,+\infty)$. Then for any $s\notin\overline{\mathcal{I}_t}$ we either have $\omega_{\mathbf{M}}(s)=+\infty$, if $s>\mathbf{M}_{\iota}$, or $\omega_{\mathbf{N}}(t/s)=+\infty$, if $s<t\mathbf{N}_{\iota}^{-1}$. Thus formally we can even allow $s\in[0,+\infty)$, because if either $\omega_{\mathbf{M}}(s)=+\infty$ or $\omega_{\mathbf{N}}(t/s)=+\infty$, this will not effect the definition of the value $\omega_{\mathbf{M}}\check{\star}\omega_{\mathbf{N}}(t)$.

\item[$(iv)$] On the other hand, if $\mathbf{M}_{\iota}\in(0,+\infty)$ and $\mathbf{N}_{\iota}\in(0,+\infty)$, then fix $t>\mathbf{M}_{\iota}\cdot\mathbf{N}_{\iota}$. For any $s\in\mathcal{I}_t$ one has $\frac{t}{s}>\mathbf{N}_{\iota}$ and $s>\frac{t}{\mathbf{N}_{\iota}}>\mathbf{M}_{\iota}$ and so $\omega_{\mathbf{M}}(s)=+\infty$ and $\omega_{\mathbf{N}}(t/s)=+\infty$. And if $s\ge 0$ with $s\notin\overline{\mathcal{I}_t}$, and such values exist since $\mathcal{I}_t\neq(0,+\infty)$, then as before either $\omega_{\mathbf{M}}(s)=+\infty$ or $\omega_{\mathbf{N}}(t/s)=+\infty$.
\end{itemize}

Summarizing, we have shown the following:

\begin{lemma}\label{sequencelowerlemma1}
Let $\mathbf{M}$, $\mathbf{N}$ be given with $\mathbf{M}_{\iota}>0$, $\mathbf{N}_{\iota}>0$. Then it holds that
\begin{itemize}
\item[$(a)$] $$\forall\;t\in[0,\mathbf{M}_{\iota}\cdot\mathbf{N}_{\iota}):\;\;\;\omega_{\mathbf{M}}\check{\star}\omega_{\mathbf{N}}(t)=\inf_{s\in[0,+\infty)}\{\omega_{\mathbf{M}}(s)+\omega_{\mathbf{N}}(t/s)\};$$

\item[$(b)$]
$$\forall\;t\in[0,\mathbf{M}_{\iota}\cdot\mathbf{N}_{\iota}):\;\;\;0\le\omega_{\mathbf{M}}\check{\star}\omega_{\mathbf{N}}(t)<+\infty;$$

\item[$(c)$] if $\mathbf{M}_{\iota}\in(0,+\infty)$, $\mathbf{N}_{\iota}\in(0,+\infty)$, then
$$\forall\;t\in(\mathbf{M}_{\iota}\cdot\mathbf{N}_{\iota},+\infty):\;\;\;\omega_{\mathbf{M}}\check{\star}\omega_{\mathbf{N}}(t)=+\infty.$$
\end{itemize}
\end{lemma}
Consequently, also in non-standard cases when using the above conventions we can replace $\mathcal{I}_t$ by $[0,+\infty)$ and hence the dependence on $t$ in the definition disappears. And \eqref{wedgeproplemmaequ} transfers into $\omega_{\mathbf{M}}\check{\star}\omega_{\mathbf{N}}(tu)\le\omega_{\mathbf{M}}(t)+\omega_{\mathbf{N}}(u)(<+\infty)$ for all $t\in[0,\mathbf{M}_{\iota})$ and $u\in[0,\mathbf{N}_{\iota})$.

Next we verify \emph{commutativity} in all cases.

\begin{lemma}\label{sequencelowercommlemma}
Let $\mathbf{M}$, $\mathbf{N}$ be given with $\mathbf{M}_{\iota}>0$, $\mathbf{N}_{\iota}>0$. Then
\begin{equation}\label{lowercommassoweightfct}
\forall\;t\in[0,\mathbf{M}_{\iota}\cdot\mathbf{N}_{\iota}):\;\;\;\omega_{\mathbf{M}}\check{\star}\omega_{\mathbf{N}}(t)=\omega_{\mathbf{N}}\check{\star}\omega_{\mathbf{M}}(t).
\end{equation}
\end{lemma}

\demo{Proof}
First, let $t=0$. Then, as seen before in $(i)$ we get $\omega_{\mathbf{M}}\check{\star}\omega_{\mathbf{N}}(0)=\inf_{s\in(0,\mathbf{M}_{\iota})}\{\omega_{\mathbf{M}}(s)\}=0=\inf_{s\in(0,\mathbf{N}_{\iota})}\{\omega_{\mathbf{N}}(s)\}=\omega_{\mathbf{N}}\check{\star}\omega_{\mathbf{M}}(0)$.

So let $t\in(0,\mathbf{M}_{\iota}\cdot\mathbf{N}_{\iota})$. We apply the substitution $u:=\frac{t}{s}$ in \eqref{lowertransform} and distinguish:

\begin{itemize}
\item[$(*)$] If $\mathbf{M}_{\iota}=+\infty=\mathbf{N}_{\iota}$, then this situation corresponds to $(c)$ in Lemma \ref{wedgeproplemma}.

\item[$(*)$] If $\mathbf{M}_{\iota}\in(0,+\infty)$ and $\mathbf{N}_{\iota}=+\infty$, then for any $t>0$ and $s\in\mathcal{I}_t=(0,\mathbf{M}_{\iota})$ we get $u\in(t\mathbf{M}_{\iota}^{-1},+\infty)$. This gives $\omega_{\mathbf{M}}(t/u)<+\infty$ and $\omega_{\mathbf{N}}(u)<+\infty$ for all such $u$. Moreover,  $(t\mathbf{M}_{\iota}^{-1},+\infty)$ corresponds to the natural interval under consideration for the definition of $\omega_{\mathbf{N}}\check{\star}\omega_{\mathbf{M}}(t)$ and so commutativity follows.

\item[$(*)$] If $\mathbf{M}_{\iota}=+\infty$ and $\mathbf{N}_{\iota}\in(0,+\infty)$, then for any $t>0$ and $s\in\mathcal{I}_t=(t\mathbf{N}_{\iota}^{-1},+\infty)$ we get $u\in(0,\mathbf{N}_{\iota})$ which gives $\omega_{\mathbf{M}}(t/u)<+\infty$ and $\omega_{\mathbf{N}}(u)<+\infty$ for all such $u$ and $(t\mathbf{M}_{\iota}^{-1},\mathbf{N}_{\iota})=(0,\mathbf{N}_{\iota})$ corresponds to the interval under consideration for defining $\omega_{\mathbf{N}}\check{\star}\omega_{\mathbf{M}}(t)$. Consequently, commutativity follows in this case, too.

\item[$(*)$] If $\mathbf{M}_{\iota}\in(0,+\infty)$ and $\mathbf{N}_{\iota}\in(0,+\infty)$, then $\mathbf{M}_{\iota}\cdot\mathbf{N}_{\iota}\in(0,+\infty)$. For any $t\in(0,\mathbf{M}_{\iota}\cdot\mathbf{N}_{\iota})$ and any $s\in\mathcal{I}_t=(t\mathbf{N}_{\iota}^{-1},\mathbf{M}_{\iota})$ we get $u\in(t\mathbf{M}_{\iota}^{-1},\mathbf{N}_{\iota})$. This corresponds to the interval under consideration for $\omega_{\mathbf{N}}\check{\star}\omega_{\mathbf{M}}(t)$ and $\omega_{\mathbf{M}}(t/u)<+\infty$, $\omega_{\mathbf{N}}(u)<+\infty$ for all such $u$. Therefore, commutativity is valid as well.
\end{itemize}
\qed\enddemo

Finally, consider the ``most extreme cases'' $\mathbf{M}_{\iota}=0$ or $\mathbf{N}_{\iota}=0$ and note that for log-convex sequences $\mathbf{M},\mathbf{N}\in\RR_{>0}^{\NN}$ this is excluded. But when allowing formally the log-convex sequence given by $M_0\ge 0$ and $M_p=0$ for all $p\in\NN_{>0}$, then $\mathbf{M}_{\iota}=0$. This corresponds to the regularization approach studied in \cite[Sect. 4.1]{regularnew}; see also $(d)$ in Section \ref{assofunctsect}.

\begin{lemma}\label{exoticlowertrafolemma}
Let $\mathbf{M},\mathbf{N}\in\RR_{>0}^{\NN}$ be given and assume that either $\mathbf{M}_{\iota}=0$ or $\mathbf{N}_{\iota}=0$. Then
$$\omega_{\mathbf{M}}\check{\star}\omega_{\mathbf{N}}(0)=0,\hspace{15pt}\omega_{\mathbf{M}}\check{\star}\omega_{\mathbf{N}}(t)=+\infty,\;\;\;\forall\;t>0.$$
\end{lemma}

\demo{Proof}
We distinguish between different situations:

\begin{itemize}
\item[$(\ast)$] $\mathbf{M}_{\iota}=0$ and $\mathbf{N}_{\iota}>0$. If $\mathbf{N}_{\iota}\in(0,+\infty)$, then $\mathcal{I}_t=\emptyset$ for any $t>0$ and we only consider $t=0$ and get $\mathcal{I}_0=\{0\}$. The convention $\frac{0}{0}=0$ gives $\omega_{\mathbf{M}}\check{\star}\omega_{\mathbf{N}}(0)=0$. Formally, we can put $\omega_{\mathbf{M}}\check{\star}\omega_{\mathbf{N}}(t)=+\infty$ for all $t>0$ when considering only $s=0$ in the infimum and this is then also justified by the conventions $t/0=+\infty$, $\omega_{\mathbf{N}}(+\infty)=+\infty$ and since $\omega_{\mathbf{M}}(s)=+\infty$ for all $s>0$; see $(d)$ in Section \ref{assofunctsect}.

    If $\mathbf{N}_{\iota}=+\infty$, then $\frac{1}{+\infty}=0$ gives $I_t=\{0\}$ for any $t\ge 0$ and so the conventions $\frac{0}{0}=0$ and $\frac{1}{0}=+\infty$ imply again $\omega_{\mathbf{M}}\check{\star}\omega_{\mathbf{N}}(0)=0$ and $\omega_{\mathbf{M}}\check{\star}\omega_{\mathbf{N}}(t)=+\infty$ for all $t>0$.

\item[$(\ast)$] $\mathbf{M}_{\iota}>0$ and $\mathbf{N}_{\iota}=0$. If $\mathbf{M}_{\iota}\in(0,+\infty)$, then for $t=0$ the conventions $\frac{1}{0}=+\infty$ and $0=0\cdot(+\infty)$ give $\mathcal{I}_0=(0,\mathbf{M}_{\iota})$ and for any $t>0$ the interval $\mathcal{I}_t$ is empty when using $\mathbf{N}_{\iota}^{-1}=+\infty$. Therefore, again $\omega_{\mathbf{M}}\check{\star}\omega_{\mathbf{N}}(0)=0$ and we can put $\omega_{\mathbf{M}}\check{\star}\omega_{\mathbf{N}}(t)=+\infty$ for all $t>0$ which is justified by the fact that here $\omega_{\mathbf{N}}(t/s)=+\infty$ for all $s>0$ if $t>0$ is arbitrary but fixed.

    If $\mathbf{M}_{\iota}=+\infty$, then for $t=0$ one has $\mathcal{I}_0=(0,+\infty)$ and for any $t>0$ let us set $I_t=\{+\infty\}$. Hence $\omega_{\mathbf{M}}\check{\star}\omega_{\mathbf{N}}(0)=0$ and using $\omega_{\mathbf{M}}(+\infty)=+\infty$ and $\frac{1}{+\infty}=0$ gives $\omega_{\mathbf{M}}\check{\star}\omega_{\mathbf{N}}(t)=+\infty$ for any $t>0$.

 \item[$(\ast)$] $\mathbf{N}_{\iota}=0=\mathbf{M}_{\iota}$. Then $\omega_{\mathbf{M}}(s)=+\infty=\omega_{\mathbf{N}}(t/s)$ for all $t,s>0$ and $I_0=\{0\}$ by the convention $0\cdot(+\infty)=0$ whereas $I_t=\emptyset$ for any $t>0$. Therefore, $\omega_{\mathbf{M}}\check{\star}\omega_{\mathbf{N}}(0)=0$ (using $\frac{0}{0}=0$) and $\omega_{\mathbf{M}}\check{\star}\omega_{\mathbf{N}}(t)=+\infty$ for all $t>0$ when formally again only considering $s=0$; see the first case before.
\end{itemize}
\qed\enddemo

Next we verify that Theorem \ref{propBomanLemma4} holds even if either $\mathbf{M}_{\iota}<+\infty$ or $\mathbf{N}_{\iota}<+\infty$. First, concerning \eqref{propBomanLemma4equ} we point out that naturally one considers $t\in[0,\mathbf{M}_{\iota}\cdot\mathbf{N}_{\iota})$ if $\mathbf{M}_{\iota},\mathbf{N}_{\iota}\in(0,+\infty)$. Let $\mathbf{M},\mathbf{N}\in\RR_{>0}^{\NN}$ be log-convex:

\begin{itemize}
\item[$(i)$] If $\mathbf{M}_{\iota}\cdot\mathbf{N}_{\iota}\in(0,+\infty)$, i.e. $\mathbf{M}_{\iota}\in(0,+\infty)$ and $\mathbf{N}_{\iota}\in(0,+\infty)$, then \eqref{propBomanLemma4equ} can be extended to all $t\in[0,+\infty)$ when we set both sides to be equal to $+\infty$ if $t>\mathbf{M}_{\iota}\cdot\mathbf{N}_{\iota}$. And for $t=\mathbf{M}_{\iota}\cdot\mathbf{N}_{\iota}$ we consider $\lim_{t\rightarrow\mathbf{M}_{\iota}\cdot\mathbf{N}_{\iota}}\omega_{\mathbf{M}\cdot\mathbf{N}}(t)$.

    This is justified by $(ii)$ in Lemma \ref{assofunctsectlemma} applied to $\mathbf{M}\cdot\mathbf{N}$ and by $(c)$ in Lemma \ref{sequencelowerlemma1}.

\item[$(ii)$] And by taking into account Lemma \ref{exoticlowertrafolemma} this identity formally even holds if either $\mathbf{M}_{\iota}=0$ or $\mathbf{N}_{\iota}=0$ on whole $[0,+\infty)$: take on both sides for $t=0$ the value $0$ and $+\infty$ for all $t>0$. Therefore recall the convention $0\cdot(+\infty)=0$, $(c)$ in Section \ref{preliminarysection} and $(d)$ in Section \ref{assofunctsect}. The sequence for which $(\cdot)_{\iota}=0$ is equal to $0$ for all $p\in\NN_{>0}$ and so $M_p\cdot N_p=0$ for all $p\in\NN_{>0}$ too; recall that both sequences are assumed to be log-convex. Moreover, the convention $0\cdot(+\infty)=0$ which is used in the proof of Lemma \ref{exoticlowertrafolemma} also yields that \eqref{propBomanLemma4equ} is only required to hold for $t=0$ (and for this value one has the equality $0=0$).
\end{itemize}

Let $\mathbf{M},\mathbf{N}\in\RR_{>0}^{\NN}$ be log-convex and that either $\mathbf{M}_{\iota}\in(0,+\infty)$ or $\mathbf{N}_{\iota}\in(0,+\infty)$, then the proof of Theorem \ref{propBomanLemma4} changes as follows:

\begin{itemize}
\item[$(a)$] $\mathbf{M}_{\iota}\in(0,+\infty)$ and $\mathbf{N}_{\iota}=+\infty$; i.e. $\mathcal{I}_t=(0,\mathbf{M}_{\iota})$ for all $t\ge 0$. Then in \eqref{stdecompose} let $t\ge 0$ and $s\in(0,\mathbf{M}_{\iota})$ and consider the supremum over all $p\in\NN$ and the infimum over all $s\in\mathcal{I}_t$. For the converse, we follow the proof of Theorem \ref{propBomanLemma4} with $t\in[\mu_1\nu_1,+\infty)$, $s_t\in[\mu_1,\mathbf{M}_{\iota}]$ resp. $t\in[0,\mu_1\nu_1)$, $s_t\in(0,\mu_1)$ and use Lemma \ref{lemma1}.\vspace{6pt}

\item[$(b)$] $\mathbf{M}_{\iota}=+\infty$ and $\mathbf{N}_{\iota}\in(0,+\infty)$; i.e. $\mathcal{I}_t=(t\mathbf{N}_{\iota}^{-1},+\infty)$ for all $t\ge 0$. Then in \eqref{stdecompose} let $t\ge 0$ and $s\in\mathcal{I}_t$ and consider the supremum over all $p\in\NN$ and the infimum over all $s\in\mathcal{I}_t$. For the converse, we follow the proof of Theorem \ref{propBomanLemma4} with $t\in[\mu_1\nu_1,+\infty)$, $s_t\in[\mu_1,+\infty)$ and so $\frac{t}{s_t}\in[\nu_1,\mathbf{N}_{\iota}]$ and $s_t\in[t\mathbf{N}_{\iota}^{-1},+\infty)$. When $t\in[0,\mu_1\nu_1)$, $s_t\in(0,\mu_1)$ with $\frac{t}{s_t}\in[0,\nu_1)$, then apply again Lemma \ref{lemma1}.\vspace{6pt}

\item[$(c)$] $\mathbf{M}_{\iota}\in(0,+\infty)$ and $\mathbf{N}_{\iota}\in(0,+\infty)$; i.e. $\mathcal{I}_t=(t\mathbf{N}_{\iota}^{-1},\mathbf{M}_{\iota})$. Then in \eqref{stdecompose} let $0\le t<\mathbf{M}_{\iota}\cdot\mathbf{N}_{\iota}$ and $s\in\mathcal{I}_t$ and consider the supremum over all $p\in\NN$ and the infimum over all $s\in\mathcal{I}_t$. For the converse, we follow the proof of Theorem \ref{propBomanLemma4} and consider $t\in[\mu_1\nu_1,\mathbf{M}_{\iota}\mathbf{N}_{\iota})$, $s_t\in[\mu_1,\mathbf{M}_{\iota}]$ and $\frac{t}{s_t}\in[\nu_1,\mathbf{N}_{\iota}]$ and so $s_t\in[t\mathbf{N}_{\iota}^{-1},\mathbf{M}_{\iota}]$. When $t\in[0,\mu_1\nu_1)$, $s_t\in(0,\mu_1)$ with $\frac{t}{s_t}\in[0,\nu_1)$, then use Lemma \ref{lemma1}.
\end{itemize}

\begin{remark}\label{stabrem}
\emph{For all sequences in the above proof for which $(\cdot)_{\iota}\in(0,+\infty)$ we have to consider the particular closed interval to which the corresponding sequence of quotients belongs: In this situation it is not excluded that the sequence of quotients eventually stabilizes. However, then Lemma \ref{lemma1} still can be applied by taking into account Remark \ref{lemma1stabrem}.}
\end{remark}

\subsection{Non-standard cases for the generalized upper Legendre conjugate}\label{nonstandardsection}
Let $\mathbf{M},\mathbf{N}\in\RR_{>0}^{\NN}$ be log-convex with $\mathbf{M}_{\iota},\mathbf{N}_{\iota}\in(0,+\infty)$ and hence \emph{violating the basic assumption $(a)$} in Definition \ref{uppertransformdef}. The aim is to study \eqref{propBomanLemma4invequ} in this situation and we distinguish:

\begin{itemize}
\item[$(a)$] Let $0<\mathbf{N}_{\iota}<\mathbf{M}_{\iota}=+\infty$. Relation $\mathbf{N}\hyperlink{mtriangle}{\vartriangleleft}\mathbf{M}$ is trivial and in \eqref{uppertransform}, for any given $t\in(0,+\infty)$, we can restrict to all $s\in[0,\mathbf{N}_{\iota}t)$ since $\omega_{\mathbf{M}}(s)<+\infty$ for any $s\ge 0$ and so for $s>\mathbf{N}_{\iota}t$ one has $\omega_{\mathbf{M}}(s)-\omega_{\mathbf{N}}(s/t)=-\infty$; i.e. all such $s$ are not effecting the supremum. Therefore, write
\begin{equation}\label{uppertransform1}
\omega_{\mathbf{M}}\widehat{\star}\omega_{\mathbf{N}}(t)=\sup_{s\in[0,\mathbf{N}_{\iota}t)}\{\omega_{\mathbf{M}}(s)-\omega_{\mathbf{N}}(s/t)\},\;\;\;t\in(0,+\infty).
\end{equation}
One also has $\omega_{\mathbf{M}}\widehat{\star}\omega_{\mathbf{N}}(0)=\omega_{\mathbf{M}}(0)-\omega_{\mathbf{N}}(0)=0$ (again by using $\frac{0}{0}=0$ and $\frac{1}{0}=+\infty$). $\omega_{\mathbf{M}}(s)\le D_t+\omega_{\mathbf{N}}(s/t)$ is valid with e.g. $D_t:=\omega_{\mathbf{M}}(\mathbf{N}_{\iota}t)$ and so $\omega_{\mathbf{M}}\widehat{\star}\omega_{\mathbf{N}}(t)<+\infty$ for all $t\in(0,+\infty)$.

By analogous arguments \eqref{widehatproplemmaequ} transfers as follows:
\begin{equation}\label{widehatproplemmaequnonstandard}
\forall\;t,u\in(\mathbf{N}_{\iota}^{-1},+\infty):\;\;\;\omega_{\mathbf{M}}\widehat{\star}\omega_{\mathbf{N}}(tu)\ge\max\{\omega_{\mathbf{M}}(t)-\omega_{\mathbf{N}}(1/u),\omega_{\mathbf{M}}(u)-\omega_{\mathbf{N}}(1/t)\};
\end{equation}
indeed when $\mathbf{N}_{\iota}=+\infty$, then \eqref{widehatproplemmaequnonstandard} corresponds to \eqref{widehatproplemmaequ}. However, it is not clear that one can consider any $t,u>0$ in \eqref{widehatproplemmaequnonstandard}.

Finally, if $\frac{\mathbf{M}}{\mathbf{N}}$ is log-convex, then equality in \eqref{propBomanLemma4invequ} holds for all $t\in[0,+\infty)$ since the proof of Theorem \ref{propBomanLemma4inv} can be repeated without any changes except restricting to all $s\in[0,\mathbf{N}_{\iota}t)$.

\item[$(b)$] Let $0<\mathbf{M}_{\iota}<\mathbf{N}_{\iota}=+\infty$. Then $\mathbf{N}\hyperlink{mtriangle}{\vartriangleleft}\mathbf{M}$ and even $\mathbf{N}\hyperlink{preceq}{\preceq}\mathbf{M}$ fails. Moreover, $\frac{\mathbf{M}}{\mathbf{N}}_{\iota}=0$, $\omega_{\mathbf{M}}(s)=+\infty$ for all $s>\mathbf{M}_{\iota}$ and $\omega_{\mathbf{N}}(s/t)\in[0,+\infty)$ for all $s\ge 0$, $t>0$. Therefore, the definition yields $\omega_{\mathbf{M}}\widehat{\star}\omega_{\mathbf{N}}(t)=+\infty$ for any $t>0$ and put $\omega_{\mathbf{M}}\widehat{\star}\omega_{\mathbf{N}}(0)=0$ when using the conventions $\frac{0}{0}=0$, $\frac{1}{0}=+\infty$ and $+\infty-(+\infty)=0$. This definition is consistent with \eqref{widehatproplemmaequnonstandard} resp. with \eqref{widehatproplemmaequ} in $(a)$ in Lemma \ref{widehatproplemma} and one obtains equality in \eqref{propBomanLemma4invequ} for all $t\in[0,+\infty)$ by taking into account for the right-hand side there comment $(d)$ in Section \ref{assofunctsect} which is applied to the sequence $\frac{\mathbf{M}}{\mathbf{N}}$.

\item[$(c)$] Let $\mathbf{M}_{\iota},\mathbf{N}_{\iota}\in(0,+\infty)$. Relation $\mathbf{N}\hyperlink{mtriangle}{\vartriangleleft}\mathbf{M}$ is violated but $\mathbf{N}\hyperlink{preceq}{\preceq}\mathbf{M}$ holds. (Indeed, recall that $\mathbf{M}$ and $\mathbf{N}$ are even equivalent by $\lim_{p\rightarrow+\infty}(M_p/M_0)^{1/p}=\mathbf{M}_{\iota}$ and similarly for $\mathbf{N}$.)

     Let now $t\in(0,+\infty)$ be fixed.

     When $s\in[0,\mathbf{M}_{\iota})$ and $s\in[0,\mathbf{N}_{\iota}t)$, i.e. $s\in[0,\min\{\mathbf{M}_{\iota},\mathbf{N}_{\iota}t\})$, then $\omega_{\mathbf{M}}(s)<+\infty$ and $\omega_{\mathbf{N}}(s/t)<+\infty$.

     When $\mathbf{N}_{\iota}t<\mathbf{M}_{\iota}\Leftrightarrow t\in(0,\mathbf{M}_{\iota}(\mathbf{N}_{\iota})^{-1})$ then for $s\in(\mathbf{N}_{\iota}t,\mathbf{M}_{\iota})$ we get $\omega_{\mathbf{M}}(s)-\omega_{\mathbf{N}}(s/t)=-\infty$ and so all such $s$ are not effecting the supremum in the definition.

     When $\mathbf{M}_{\iota}<\mathbf{N}_{\iota}t\Leftrightarrow t\in(\mathbf{M}_{\iota}(\mathbf{N}_{\iota})^{-1},+\infty)$, then for $s\in(\mathbf{M}_{\iota},\mathbf{N}_{\iota}t)$ we get $\omega_{\mathbf{M}}(s)-\omega_{\mathbf{N}}(s/t)=+\infty$. Finally, for $s>\max\{\mathbf{M}_{\iota},\mathbf{N}_{\iota}t\}$ one has $\omega_{\mathbf{M}}(s)-\omega_{\mathbf{N}}(s/t)=0$ when using the convention $+\infty-(+\infty)=0$.

Summarizing, in accordance with Theorem \ref{propBomanLemma4inv} and \eqref{uppertransform1} we put
\begin{equation}\label{uppertransform2}
\omega_{\mathbf{M}}\widehat{\star}\omega_{\mathbf{N}}(t)=\sup_{s\in[0,\mathbf{N}_{\iota}t)}\{\omega_{\mathbf{M}}(s)-\omega_{\mathbf{N}}(s/t)\},\;\;\;t\in(0,(\underline{C}_{\mathbf{N}\preceq\mathbf{M}})^{-1})=\left(0,\frac{\mathbf{M}_{\iota}}{\mathbf{N}_{\iota}}\right)=\left(0,\frac{\mathbf{M}}{\mathbf{N}}_{\iota}\right);
\end{equation}
recall also \eqref{MNiotarelation}. Moreover, one has $\omega_{\mathbf{M}}\widehat{\star}\omega_{\mathbf{N}}(\frac{\mathbf{M}_{\iota}}{\mathbf{N}_{\iota}})=\lim_{t\rightarrow\frac{\mathbf{M}_{\iota}}{\mathbf{N}_{\iota}}}\omega_{\mathbf{M}}\widehat{\star}\omega_{\mathbf{N}}(t)$, $\omega_{\mathbf{M}}\widehat{\star}\omega_{\mathbf{N}}(t)=+\infty$ for $t\in(\frac{\mathbf{M}_{\iota}}{\mathbf{N}_{\iota}},+\infty)$ and $\omega_{\mathbf{M}}\widehat{\star}\omega_{\mathbf{N}}(0)=0$. If $\frac{\mathbf{M}}{\mathbf{N}}$ is log-convex, then equality in \eqref{propBomanLemma4invequ} formally holds for all $t\in[0,+\infty)$ since the proof of Theorem \ref{propBomanLemma4inv} can be repeated when restricting to all $s\in[0,\mathbf{N}_{\iota}t)$ and $t\in(0,\frac{\mathbf{M}_{\iota}}{\mathbf{N}_{\iota}})$ and note that for all $t\in(\frac{\mathbf{M}_{\iota}}{\mathbf{N}_{\iota}},+\infty)$ on both sides in \eqref{propBomanLemma4invequ} we get $+\infty$ by taking into account the fact that in this case $\frac{\mathbf{M}}{\mathbf{N}}_{\iota}=\frac{\mathbf{M}_{\iota}}{\mathbf{N}_{\iota}}$.

Concerning \eqref{widehatproplemmaequ} in $(a)$ in Lemma \ref{widehatproplemma} note that under the assumption $1<\mathbf{M}_{\iota}\cdot\mathbf{N}_{\iota}$ we get (cf. \eqref{uppertransform2})
\begin{equation*}\label{widehatproplemmaequnonstandard1}
	\forall\;t,u:\;\mathbf{N}_{\iota}^{-1}<t,u,\;\;tu<\frac{\mathbf{M}_{\iota}}{\mathbf{N}_{\iota}}:\;\;\;\omega_{\mathbf{M}}\widehat{\star}\omega_{\mathbf{N}}(tu)\ge\max\{\omega_{\mathbf{M}}(t)-\omega_{\mathbf{N}}(1/u),\omega_{\mathbf{M}}(u)-\omega_{\mathbf{N}}(1/t)\}.
	\end{equation*}
\end{itemize}

Note that the comments from Remark \ref{stabrem} apply to the cases $(a)$ and $(c)$ when following the proof in Theorem \ref{propBomanLemma4inv}.\vspace{6pt}

Finally, for the sake of completeness let us again comment on the extreme cases when either $\mathbf{M}_{\iota}=0$ or $\mathbf{N}_{\iota}=0$. Recall that in this situation the particular sequence(s) cannot be log-convex. We distinguish and for each case we use the conventions $\frac{0}{0}=0$, $\frac{1}{0}=+\infty$, and $+\infty-(+\infty)=0$:

\begin{itemize}
\item[$(*)$] Let $0=\mathbf{M}_{\iota}<\mathbf{N}_{\iota}$. Then $\omega_{\mathbf{M}}(s)=+\infty$ for any $s>0$ and $\omega_{\mathbf{M}}(0)=0$ but $\omega_{\mathbf{N}}(s/t)<+\infty$ at least for all $s$ sufficiently small (depending on given $t>0$). Then we get
    $$\forall\;t\in(0,+\infty):\;\;\;\omega_{\mathbf{M}}\widehat{\star}\omega_{\mathbf{N}}(t)=+\infty,$$
 and set $\omega_{\mathbf{M}}\widehat{\star}\omega_{\mathbf{N}}(0):=0$ which is justified by the conventions above. Therefore, $\omega_{\mathbf{M}}\widehat{\star}\omega_{\mathbf{N}}=\omega_{\mathbf{M}}$ and in view of $(d)$ in Section \ref{assofunctsect} this is consistent with \eqref{widehatproplemmaequ}.

\item[$(*)$] Let $0=\mathbf{N}_{\iota}<\mathbf{M}_{\iota}$. Since $\omega_{\mathbf{N}}(u)=+\infty$ for any $u>0$ and $\omega_{\mathbf{N}}(0)=0$, by the above conventions we get $\omega_{\mathbf{M}}\widehat{\star}\omega_{\mathbf{N}}(t)=\omega_{\mathbf{M}}(0)-\omega_{\mathbf{N}}(0)=0$ for all $t\in(0,+\infty)$ and also $\omega_{\mathbf{M}}\widehat{\star}\omega_{\mathbf{N}}(0)=0$.

    Summarizing, $\omega_{\mathbf{M}}\widehat{\star}\omega_{\mathbf{N}}=0$ holds and \eqref{widehatproplemmaequ} transfers into the (formal) ``estimate'' $0\ge-\infty$ for all $t,u>0$.

\item[$(*)$] Let $\mathbf{M}_{\iota}=0=\mathbf{N}_{\iota}$. Here $\omega_{\mathbf{N}}(s/t)=+\infty$ and $\omega_{\mathbf{M}}(s)=+\infty$ for any $s,t>0$. The convention $+\infty-(+\infty)=0$ gives again $\omega_{\mathbf{M}}\widehat{\star}\omega_{\mathbf{N}}(t)=\omega_{\mathbf{M}}(0)-\omega_{\mathbf{N}}(0)=0$ for all $t\in(0,+\infty)$. And for $t=0$ the above conventions  yield $\omega_{\mathbf{M}}\widehat{\star}\omega_{\mathbf{N}}(0)=0$. Thus $\omega_{\mathbf{M}}\widehat{\star}\omega_{\mathbf{N}}=0$ is verified and \eqref{widehatproplemmaequ} turns into the (formal) equality $0=0=+\infty-(+\infty)$ for all $t,u>0$.
\end{itemize}

\subsection{Inverse operations for non-standard cases}\label{invnonstandardsection}
Finally, we transfer Theorem \ref{propBomanLemma4comb} to some non-standard situations:

\begin{theorem}\label{propBomanLemma4combvar}
Let $\mathbf{M},\mathbf{N}\in\RR_{>0}^{\NN}$ be log-convex and hence $\mathbf{M}_{\iota},\mathbf{N}_{\iota}\in(0,+\infty]$.
\begin{itemize}
\item[$(i)$] If $\mathbf{N}_{\iota}=+\infty$, then
    $$\forall\;t\in[0,\mathbf{M}_{\iota}):\;\;\;(\omega_{\mathbf{M}}\check{\star}\omega_{\mathbf{N}})\widehat{\star}\omega_{\mathbf{N}}(t)=\omega_{\mathbf{M}}(t),$$
    and if $\mathbf{M}_{\iota}=+\infty$, then
    $$\forall\;t\in[0,\mathbf{N}_{\iota}):\;\;\;(\omega_{\mathbf{M}}\check{\star}\omega_{\mathbf{N}})\widehat{\star}\omega_{\mathbf{M}}(t)=\omega_{\mathbf{N}}(t).$$

\item[$(ii)$] Assume that $\mathbf{M}_{\iota}=+\infty$, $\frac{\mathbf{M}}{\mathbf{N}}$ is log-convex and that $\mathbf{N}\hyperlink{preceq}{\preceq}\mathbf{M}$. Then one gets
$$\forall\;t\in[0,+\infty):\;\;\;\omega_{\mathbf{M}}(t)=\omega_{\mathbf{N}}\check{\star}(\omega_{\mathbf{M}}\widehat{\star}\omega_{\mathbf{N}})(t)=(\omega_{\mathbf{M}}\widehat{\star}\omega_{\mathbf{N}})\check{\star}\omega_{\mathbf{N}}(t).$$
\end{itemize}
Consequently, in $(i)$ for the equality one has to restrict to a finite interval in the case $\mathbf{M}_{\iota}<+\infty$ resp. $\mathbf{N}_{\iota}<+\infty$ whereas the equality in $(ii)$ is the same as the one in $(ii)$ in Theorem \ref{propBomanLemma4comb}; hence this identity is even valid under non-standard (weaker) assumptions.
\end{theorem}

\demo{Proof}
$(i)$ First, by the comments in Section \ref{nonstandardlowersection} concerning non-standard cases for Theorem \ref{propBomanLemma4} we get
$$t\in[0,\mathbf{M}_{\iota}\cdot\mathbf{N}_{\iota})=[0,+\infty):\;\;\;\omega_{\mathbf{M}\cdot\mathbf{N}}(t)=\omega_{\mathbf{M}}\check{\star}\omega_{\mathbf{N}}(t).$$ Then distinguish: When $\mathbf{N}_{\iota}=+\infty$, then apply Theorem \ref{propBomanLemma4inv} to $\mathbf{P}=\mathbf{M}\cdot\mathbf{N}$ and $\mathbf{Q}=\mathbf{N}$, get $(\underline{C}_{\mathbf{Q}\preceq\mathbf{P}})^{-1}=\frac{\mathbf{P}}{\mathbf{Q}}_{\iota}=\mathbf{M}_{\iota}$ and Theorem \ref{propBomanLemma4inv} implies $\omega_{\mathbf{M}\cdot\mathbf{N}}\widehat{\star}\omega_{\mathbf{N}}(t)=\omega_{\mathbf{P}}\widehat{\star}\omega_{\mathbf{Q}}(t)=\omega_{\mathbf{M}}(t)$ for all $t\in[0,\mathbf{M}_{\iota})$. Combining both identities the desired equation is verified and the case $\mathbf{M}_{\iota}=+\infty$ is similar with $\mathbf{P}=\mathbf{M}\cdot\mathbf{N}$, $\mathbf{Q}=\mathbf{M}$.\vspace{6pt}

$(ii)$ \emph{Case I:} $\mathbf{N}_{\iota}=+\infty$. First, by the assumptions Theorem \ref{propBomanLemma4inv} gives
\begin{equation}\label{propBomanLemma4combvarequ1}
\forall\;t\in\left[0,\frac{\mathbf{M}}{\mathbf{N}}_{\iota}\right):\:\:\;\omega_{\mathbf{M}}\widehat{\star}\omega_{\mathbf{N}}(t)=\omega_{\frac{\mathbf{M}}{\mathbf{N}}}(t).
\end{equation}
Then apply the comments in Section \ref{nonstandardlowersection} concerning non-standard cases for Theorem \ref{propBomanLemma4} to $\mathbf{P}=\mathbf{N}$ and $\mathbf{Q}=\frac{\mathbf{M}}{\mathbf{N}}$ (note that $\frac{\mathbf{M}}{\mathbf{N}}$ is log-convex) and hence
\begin{equation}\label{propBomanLemma4combvarequ}
\forall\;t\in[0,\mathbf{P}_{\iota}\cdot\mathbf{Q}_{\iota})=\left[0,\mathbf{N}_{\iota}\cdot\frac{\mathbf{M}}{\mathbf{N}}_{\iota}\right)=[0,+\infty):\;\;\;\omega_{\mathbf{M}}(t)=\omega_{\mathbf{P}\cdot\mathbf{Q}}(t)=\omega_{\mathbf{N}}\check{\star}\omega_{\frac{\mathbf{M}}{\mathbf{N}}}(t).
\end{equation}
Next recall that for $\omega_{\mathbf{N}}\check{\star}\omega_{\frac{\mathbf{M}}{\mathbf{N}}}(t)$ naturally one has to consider the infimum over all $s\in\mathcal{I}_t=(t(\frac{\mathbf{M}}{\mathbf{N}}_{\iota})^{-1},\mathbf{N}_{\iota})=(t(\frac{\mathbf{M}}{\mathbf{N}}_{\iota})^{-1},+\infty)$ with $t\in[0,\mathbf{N}_{\iota}\cdot\frac{\mathbf{M}}{\mathbf{N}}_{\iota})=[0,+\infty)$; see \eqref{lowertransform}. Since $\frac{\mathbf{M}}{\mathbf{N}}$ is log-convex, $\mathbf{N}\hyperlink{preceq}{\preceq}\mathbf{M}$ implies $\left(\frac{\mathbf{M}}{\mathbf{N}}_{\iota}\right)^{-1}>0$.

Thus for all $s$ under consideration in the definition of $\omega_{\mathbf{N}}\check{\star}\omega_{\frac{\mathbf{M}}{\mathbf{N}}}(t)$ one has $0\le\frac{t}{s}<\frac{\mathbf{M}}{\mathbf{N}}_{\iota}$. By combining \eqref{propBomanLemma4combvarequ1} and \eqref{propBomanLemma4combvarequ} the desired equality is shown.\vspace{6pt}

\emph{Case II:} $\mathbf{N}_{\iota}<+\infty$; i.e. $\mathbf{N}_{\iota}\in(0,+\infty)$. Because $\mathbf{M}_{\iota}=+\infty$ one infers automatically $\mathbf{N}\hyperlink{mtriangle}{\vartriangleleft}\mathbf{M}$ and comment $(a)$ in Section \ref{nonstandardsection} gives that
$$\forall\;t\in\left[0,\frac{\mathbf{M}}{\mathbf{N}}_{\iota}\right)=[0,+\infty):\;\;\;\omega_{\mathbf{M}}\widehat{\star}\omega_{\mathbf{N}}(t)=\omega_{\frac{\mathbf{M}}{\mathbf{N}}}(t).$$
Apply again the comments in Section \ref{nonstandardlowersection} concerning non-standard cases for Theorem \ref{propBomanLemma4} to $\mathbf{P}=\mathbf{N}$ and $\mathbf{Q}=\frac{\mathbf{M}}{\mathbf{N}}$; then get \eqref{propBomanLemma4combvarequ} and note that $\mathbf{N}_{\iota}\cdot\frac{\mathbf{M}}{\mathbf{N}}_{\iota}=+\infty$ because $\mathbf{N}\hyperlink{mtriangle}{\vartriangleleft}\mathbf{M}$ and log-convexity for $\frac{\mathbf{M}}{\mathbf{N}}$ yields $\frac{\mathbf{M}}{\mathbf{N}}_{\iota}=+\infty$. In $\check{\star}$ one considers the infimum over all $s\in\mathcal{I}_t=(t\left(\frac{\mathbf{M}}{\mathbf{N}}_{\iota}\right)^{-1},\mathbf{N}_{\iota})=(0,\mathbf{N}_{\iota})$ with $t\in[0,\mathbf{N}_{\iota}\cdot\frac{\mathbf{M}}{\mathbf{N}}_{\iota})=[0,+\infty)$. Combining the information we are done.
\qed\enddemo

\bibliographystyle{plain}
\bibliography{Bibliography}
\end{document}